\documentclass[sn-mathphys,Numbered]{sn-jnl}% Math and Physical Sciences Reference Style
\usepackage{graphicx}%
\usepackage{amsmath,amssymb,amsfonts,bm}%
\usepackage{multirow}
\usepackage{amsthm}%
\usepackage{mathrsfs}%
\usepackage[title]{appendix}%
\usepackage{xcolor}%
\usepackage{textcomp}%
\usepackage{manyfoot}%
\usepackage{booktabs}%
\usepackage{algpseudocode}%
\usepackage{listings}%
\usepackage{bigints}
\usepackage{outlines}
\usepackage{geometry}
\usepackage{subfigure}
\usepackage{siunitx}
\usepackage{float}
\usepackage{comment}
\usepackage[percent]{overpic}

\usepackage[onehalfspacing]{setspace}
\usepackage[utf8]{inputenc}
\usepackage{graphicx}
\usepackage[ruled,vlined]{algorithm2e}
\usepackage{amsmath}
\usepackage{enumitem}
\usepackage{color,soul}
\usepackage{xfrac}

\newcommand{\bb}{\boldsymbol}
\usepackage{nicefrac}
\usepackage{mathtools}

\usepackage{hyperref}
\begin{document}

\title[Article Title]{A Jacobian-free Newton-Krylov method for high-order cell-centred finite volume solid mechanics}

\author[1,2]{\fnm{Ivan} \sur{Batisti\'{c}}}
\author*[1,3]{\fnm{Pablo} \sur{Castrillo}}\email{pablo.castrillo@ucd.ie}
\author[1]{\fnm{Philip} \sur{Cardiff}}

\affil*[1]{\orgdiv{School of Mechanical and Materials Engineering}, \orgname{University College Dublin}, \orgaddress{\country{Ireland}}}
\affil[2]{\orgdiv{Faculty of Mechanical Engineering and Naval Architecture}, \orgname{University of Zagreb}, \orgaddress{\country{Croatia}}}
\affil[3]{\orgdiv{Instituto de Estructuras y Transporte, Facultad de Ingeniería}, \orgname{Universidad de la República}, \orgaddress{\country{Uruguay}}}

\abstract
{
%\coffeestainA{0.2}{0.85}{-25}{5cm}{10.3cm}
This work extends the application of Jacobian-free Newton–Krylov (JFNK) methods to higher-order cell-centred finite-volume formulations for solid mechanics. While conventional schemes are typically limited to second-order accuracy, we present third- and fourth-order formulations employing local least-squares reconstructions for gradient evaluation and Gaussian quadrature at cell faces. These schemes enable accurate resolution of complex stress and deformation fields in linear and nonlinear solids while retaining the flexibility of finite-volume methods.
A key contribution is a JFNK solution strategy for these higher-order schemes, eliminating the need to assemble complex Jacobian matrices. A compact-stencil approximate Jacobian is used as a preconditioner, providing efficiency gains similar to second-order frameworks. To enhance robustness on irregular meshes, an $\alpha$-stabilisation scheme is incorporated, damping high-frequency error modes without compromising formal accuracy.
The proposed methodology is benchmarked across a suite of two- and three-dimensional test problems involving elastic and nonlinear materials, where key performance metrics, including accuracy, computational cost, memory usage, and robustness, are systematically evaluated.
Results confirm that the higher-order formulations deliver substantial accuracy improvements over second-order schemes, while the JFNK approach achieves strong performance with only minimal modifications to existing segregated frameworks.
These findings underscore the potential of combining higher-order finite-volume methods with JFNK solvers to advance the state of the art in computational solid mechanics.
The implementations are openly released in the solids4foam toolbox for OpenFOAM, supporting further exploration and adoption by the community.
}


\maketitle

\newpage
%%%%%%%%%%%%%%%%%%%%%%%%%%%%%%%%%%%%%%%%%%%%%%%%%%%%%%%%%%%%%%%%%%
\section{Introduction}\label{sec:intro}
High-order finite-volume discretisations are generally defined as schemes that achieve at least third-order accuracy \cite{wang2013}. Over the past two decades, such methods have attracted growing interest within the computational fluid dynamics (CFD) community, leading to continuous development and refinement \cite{Tsoutsanis2025, Nishikawa2025, Costa2022, Costa2023, Antoniadis2022}. In contrast, the application of high-order finite volume discretisations to solid mechanics remains relatively underexplored  \cite{Demirdzic2016, Castrillo2022, Castrillo2024}, despite several recent advances. The untapped potential of these formulations is significant: for a given level of accuracy, high-order methods can be computationally more efficient and reduce data movement, an increasingly important factor for large-scale simulations on high-performance computing systems, where memory bandwidth is often the main bottleneck. As noted in \cite{Tzanio2021}, these two criteria, computational efficiency and reduced data movement, are guiding principles in the development of next-generation discretisations. Within the solid mechanics community, the development of such methods lags decades behind the finite-element counterparts \cite{Bathe1996}. Their greatest promise lies in multiphysics simulations, where solving the coupled problem within a unified finite volume framework is particularly advantageous.

The first high-order, cell-centred finite volume method for solid mechanics was introduced by Demirdžić \cite{Demirdzic2016}, who also pioneered the general finite volume approach for solids \cite{Demirdzic1988, Cardiff2021}. In his work, a fourth-order discretisation for two-dimensional structured hexahedral meshes was proposed, and it was shown that shear locking is not inherent to high-order finite volume formulations but is a limitation of conventional second-order schemes. However, Demirdžić’s formulation was not easily generalisable to irregular meshes, and he suggested adopting a Moving Least-Squares (MLS) reconstruction, similar to those already used in CFD \cite{Felgueroso2007, Ramirez2014}. This step was realised later by Castrillo et. al. \cite{Castrillo2022, Castrillo2023}, where it was demonstrated that higher-order accuracy can be achieved through gradient reconstruction using least-squares interpolation, specifically the Local Regression Estimator (LRE) technique, and face-flux integration via Gaussian quadrature.

A major distinction between these studies lies in their solution strategies. Demirdžić \cite{Demirdzic2016} employed a semi-implicit, segregated framework with Picard iterations, whereas Castrillo et. al. \cite{Castrillo2022,Castrillo2024} adopted a fully implicit, block-coupled Newton–Raphson approach. While segregated solvers are simpler to implement, they typically exhibit only linear convergence and limited robustness. The Newton–Raphson method, in contrast, achieves quadratic convergence but requires the explicit formation and storage of the Jacobian matrix, which is both computationally expensive and memory-intensive. Moreover, constructing an exact Jacobian for nonlinear material models is often non-trivial, further complicating implementation. 

To overcome these limitations, Jacobian-free Newton–Krylov (JFNK) methods have emerged as a powerful alternative. JFNK algorithms retain the quadratic convergence properties of Newton’s method while avoiding explicit Jacobian construction by estimating matrix-vector products through finite differencing of the residual. These methods were first introduced in the 1980s and early 1990s for solving ordinary and partial differential equations \citep{Gear1983, Chan1984, Brown1986, Brown1990}, and subsequently applied to the Navier–Stokes equations with various preconditioners and linear solvers  \cite{Qin2000, Geuzaine2001, Pernice2001, Nejat2008}. Despite their potential, JFNK methods are still not widely adopted in CFD, although some codes provide them as built-in options \cite{Nakashima2025}. Of particular relevance is that recent developments in high-order finite volume schemes increasingly rely on JFNK formulations, due to their natural compatibility with compact-stencil preconditioners and flexible treatment of nonlinearities \cite{Nishikawa2025}. Notably, the preconditioning matrix can be constructed from a lower-order spatial discretisation, which provides a good approximation to the true Jacobian while maintaining computational efficiency and robustness \cite{Mchugh1994, Qin2000, Geuzaine2001, Pernice2001, Knoll2004, Vaassen2008, Lucas2010, Nejat2008, Nishikawa2020}. 

The application of JFNK methods to cell-centred finite volume solid mechanics was first demonstrated by the present authors in previous work \cite{Cardiff2025}. There, it was shown that a compact-stencil approximate Jacobian is highly effective as a preconditioner and that the resulting JFNK solution procedure, when combined with a second-order finite volume discretisation, offers a feasible, efficient, and robust framework applicable to static, dynamic, linear, and nonlinear solid mechanics problems. However, the combination of JFNK with high-order finite volume formulations for solids has remained unexplored prior to the current study.

The present work addresses this gap by extending the JFNK framework to third- and fourth-order cell-centred finite volume formulations for both linear and nonlinear solid mechanics. High-order accuracy is achieved through least-squares-based gradient reconstruction and high-order flux integration, following the approaches validated in prior studies \cite{Castrillo2022, Castrillo2024}. We show that even in this higher-order setting, a compact-stencil Jacobian derived from a lower-order scheme remains an effective preconditioner for the Krylov solver, thereby preserving the efficiency and fast convergence of the JFNK approach. To our knowledge, this work provides the first demonstration of a JFNK solver successfully applied to high-order (beyond second-order) finite-volume solid mechanics formulations.

Another key contribution of this study is its implementation within a conventional open-source framework, which distinguishes it from earlier custom-code prototypes.
Specifically, we build our solver using the OpenFOAM-based solids4foam toolbox \cite{Weller1998, Cardiff2025a}, maintaining its support for arbitrary polyhedral cells by decomposing cells and faces into tetrahedral and triangles.
This \emph{conventional} implementation means that the high-order JFNK solver can handle general unstructured meshes and can be directly compared against standard second-order solvers.
The implementation is streamlined to be more consistent with the OpenFOAM paradigm, removing the additional boundary unknowns required in \cite{Castrillo2022, Castrillo2024} which required different treatment of boundary conditions. Finally, we introduce the $\alpha$-stabilisation scheme \cite{Nishikawa2010, Nishikawa2011, Nishikawa2025a, Antoniadis2022}, and show its effectiveness in stabilising solutions on unstructured (irregular) grids and suppressing high-frequency oscillations (zero-energy modes). The concept of $\alpha$-stabilisation originates from CFD but is shown here to be naturally compatible with high-order solid mechanics finite volume discretisations. Importantly, it scales consistently with the order of accuracy, unlike conventional second-order stabilisations (e.g. second order Rhie–Chow-type stabilisation \cite{Cardiff2025a}) used in segregated solvers. The aim of this paper is therefore twofold: to contribute to the relatively sparse literature on high-order cell-centred finite-volume methods for solids, and to explore the powerful synergy between such methods and Jacobian-free Newton-Krylov solution strategies.

The remainder of this paper is structured as follows: Section~\ref{sec:math_model} presents the governing equations, followed by the formulation of the high-order finite volume discretisation in Section~\ref{sec:HO-FVD}. Section~\ref{sec:Sol_Alg} summarises the JFNK solution algorithm. Section~\ref{sec:test_cases} reports the performance of the proposed solver across a series of benchmark problems, analysing accuracy, computational cost, robustness and memory requirements. Finally, Section~\ref{sec:conclusion} provides concluding remarks and outlines future research directions.
%
%%%%%%%%%%%%%%%%%%%%%%%%%%%%%%%%%%%%%%%%%%%%%%%%%%%%%%%%%%%%%%%%%%%%%%%%%%%%%%%
%
\section{Mathematical Model}
\label{sec:math_model}
%
%%%%%%%%%%%%%%%%%%%%%%%%%%%%%%%%%%%%%%%%%%%%%%%%%%%%%%%%%%%%%%%%%%%%%%%%%%%%%%%
%
For an arbitrary body of volume $\Omega$ bounded by surface $\Gamma$ with outward facing unit normal $\bb{n}$
the strong integral form of linear momentum in \emph{total} Lagrangian form  in the initial reference configuration is:
\begin{eqnarray} \label{eqn:momentum_TL}
     \int_{\Omega_o} \rho_o \frac{\partial^2 \bb{u} }{\partial t^2} \text{d}\Omega_o
     =
     \oint_{\Gamma_o} \left( \bb{n}_o \cdot  \bb{P}^{\text{T}}\right)  \text{d}\Gamma_o
     + \int_{\Omega_o}  \bb{f}_b \, \text{d}\Omega_o,
\end{eqnarray}
where $\rho$ is density, $\bb{u}$ is the displacement vector, $\bb{P}$ is the first Piola–Kirchhoff stress tensor,
and $\bb{f}_b$ is a body force per unit volume, e.g., $\rho \bb{g}$, where $\bb{g}$ is gravity.
Subscript $o$ is used to indicate quantities in the initial reference configuration.
Through Nanson's formula it is possible to relate the first Piola–Kirchhoff stress tensor with the Cauchy stress tensor, $\bb{\sigma}$, in the current configuration \cite{Gurtin1981, Gurtin2010}:
\begin{eqnarray} \label{eqn:nanson}
\bb{P}=J \bb{\sigma} \bb{F}^{-\text{T}},
\end{eqnarray}
where $\bb{F}$ is the deformation gradient tensor defined as $\bb{I} + (\bb{\nabla}_o \bb{u})^{\text{T}}$ and $J$ is its determinant $J = \text{det}(\bb{F})$.
\\
This work also considers linear geometry formulation, i.e. small strain assumption:
\begin{eqnarray} \label{eqn:momentum_lingeom}
    \int_{\Omega} \rho \frac{\partial^2 \bb{u} }{\partial t^2} \, \text{d}\Omega
     =
     \oint_{\Gamma} \bb{n} \cdot \bb{\sigma}_\text{s} \,  \text{d}\Gamma
     + \int_{\Omega}  \bb{f}_b \, \text{d}\Omega,
\end{eqnarray}
where $\bb{\sigma}_\text{s}$ is the engineering (small strain) stress tensor which coincide with Piola stress tensor in the limit of small strains and rotations.

The constitutive relation for Cauchy stress tensor in Eqs.~\eqref{eqn:nanson}~and~\eqref{eqn:momentum_lingeom} is given by a chosen mechanical law.
Mechanical laws considered in this work (Hooke's law, Mooney–Rivlin and neo-Hookean) are briefly outlined in Appendix \ref{app:mechLaws}.
The governing equations are complemented by boundary conditions, with three types considered here: prescribed displacement, prescribed traction, and symmetry.
%
%%%%%%%%%%%%%%%%%%%%%%%%%%%%%%%%%%%%%%%%%%%%%%%%%%%%%%%%%%%%%%%%%%%%%%%%%%%%%%%
%
\section{High-order Finite Volume Discretisation}
\label{sec:HO-FVD}
%
%%%%%%%%%%%%%%%%%%%%%%%%%%%%%%%%%%%%%%%%%%%%%%%%%%%%%%%%%%%%%%%%%%%%%%%%%%%%%%%
%
In this work, a finite volume discretisation is employed to approximate the strong integral form of the governing equations.
The computational nodes are located at the cell centroids, where the solution is treated as point-valued rather than cell-averaged.
This type of discretisation is referred to as the deconvolution finite volume method \cite{Nishikawa2025, Nishikawa2021}.
%
%------------------------------------------------------------------------------
\subsection{Solution Domain Discretisation}
\begin{figure}[h]
 	\centering
	\subfigure[Control volume in 2D]
 	{
 		\label{fig:cell2D}
    	\includegraphics[scale=1]{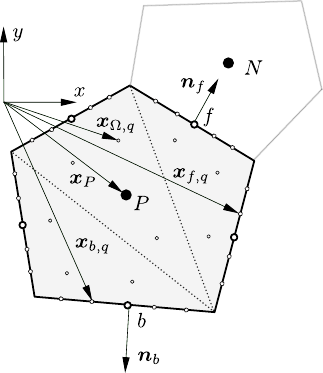}
    }
    \quad
 	\subfigure[Control volume in 3D]
 	{
 		\label{fig:cell3D}
    	\includegraphics[scale=1]{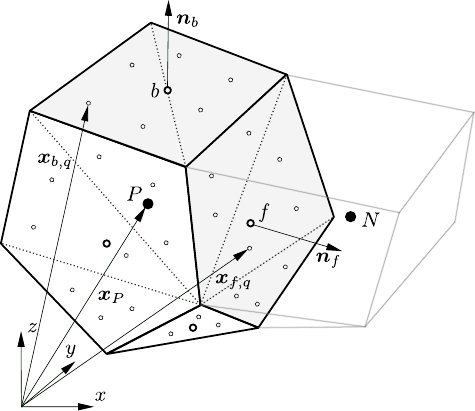}
    }
 	\caption{Representative convex polyhedral cell $P$ and neighbouring cell $N$. Quantities calculated at internal and boundary faces are denoted by the subscripts $f$ and $b$, respectively.}
 	\label{fig:cell}
\end{figure}

%------------------------------------------------------------------------------
%
The spatial domain is partitioned into a finite set of contiguous convex polyhedral control volumes, each denoted by $P$.
Representative control volumes in 2D and 3D settings are shown in Figs.~\ref{fig:cell2D}~and~\ref{fig:cell3D}, respectively.
The computational node associated with each cell $P$ is positioned at the cell centroid $\bb{x}_P$, the cell volume is denoted by $\Omega_P$,
and the centroid of a neighboring control volume $N$ is denoted by $\bb{x}_N$.
Each control volume is bounded by a set of polygonal faces, which are categorised as follows:
\begin{itemize}
\item[•] Internal faces — shared between adjacent control volumes. The centroid of an internal face is denoted by $\bb{x}_f$, its outward unit normal vector by $\bb{n}_f$, and its face area by $\Gamma_f$. 
\item[• ] Boundary faces — located on the boundary of the spatial domain. The centroid of a boundary face is denoted by $\bb{x}_b$, its outward unit normal vector by $\bb{n}_b$, and its face area by $\Gamma_b$. 
\end{itemize}
Accurate and robust flux integration is required over these arbitrary convex polyhedral volumes.
To facilitate this, in 3D cases, each polygonal face $f$ is subdivided into triangular subfaces,
on which quadrature points $\bb{x}_{f,g}$ are defined for numerical flux evaluation.
The fan triangulation method is adopted for face decomposition to minimise the number of integration points.
%For robustness, centre-point decomposition is used whenever the fan triangulation produces degenerate or near-zero-area triangles.
Although the computational mesh is always three-dimensional,
two-dimensional problems are treated as planar, where the surface integration reduces to line integration along edges,
therefore, no subdivision is required.
Each cell volume $P$ is decomposed into tetrahedral
(in 3D) or triangles (in 2D) subelements, with volume quadrature points $\bb{x}_{\Omega,q}$ defined within each subelement.
For all geometric entities, Gaussian quadrature rules are employed.
Specifically, Gauss–Legendre quadrature is used for line integrals,
Dunavant’s symmetric Gaussian quadrature rules \cite{Dunavant1985} are applied for triangular subelements,
and Shunn and Ham’s symmetric quadrature rules \cite{Shunn2012} are adopted for volume integration over tetrahedral subelements.
A key advantage of this geometric framework is its unified treatment of arbitrary convex polyhedral topologies, where all control
volumes are handled consistently using the same reconstruction and integration procedures.

Before proceeding with the discretisation of the volume and surface integrals, we introduce the following sets of cell faces.
The set of faces of a cell $P$ is denoted $\mathcal{F}_P=\mathcal{F}^{\text{int}}_P \cup \mathcal{F}^{\text{bnd}}_P$, 
where $\mathcal{F}_P^{\text{int}}$ is the set of internal faces, and $\mathcal{F}_P^{\text{bnd}}$ is the set of boundary faces.
The boundary-face set $\mathcal{F}_P^{\text{bnd}}$ is further partitioned into three mutually disjoint subsets,
$\mathcal{F}_P^{\mathrm{bnd}} := \mathcal{F}_P^{\text{disp}} \cup \mathcal{F}_P^{\text{sym}} \cup \mathcal{F}_P^{\text{trac}}$,
representing boundary faces where displacement ($\mathcal{F}_P^{\text{disp}}$), traction ($\mathcal{F}_P^{\text{trac}}$),
 or symmetry ($\mathcal{F}_P^{\text{sym}}$) conditions are prescribed. For convenience, we additionally define the set of non-traction boundary faces as $\mathcal{F}_P^{\text{non-trac}} := \mathcal{F}_P^{\text{disp}} \cup \mathcal{F}_P^{\text{sym}}$.
%
%------------------------------------------------------------------------------
\subsection{Surface Integrals}
\label{sec:surf_int}
%------------------------------------------------------------------------------
%
The surface integral term in Eq.~\eqref{eqn:momentum_lingeom} (or Eq.~\eqref{eqn:momentum_TL}) is discretised by integrating the stress over the quadrature points of each cell face:
\begin{eqnarray}\label{eq:divStressDiscret}
\begin{split}
&\oint_{\Gamma_P} \bb{n} \cdot \bb{\sigma}_\text{s} \; \text{d}\Gamma_P
=
\sum_{f \, \in \,\mathcal{F}_P} \int_{\Gamma_f} \bb{n} \cdot \bb{\sigma}_\text{s} \; \text{d}\Gamma_f
\approx
\sum_{f \,\in\, \mathcal{F}_P^{\text{int}}} \bb{n}_{f} \cdot \left [ \sum_{q=1}^{q=N_{f,q}}\alpha_q \; (\bb{\sigma}_\text{s})_{f,q} \right]\Gamma_f \\
& + \sum_{b\, \in \,\mathcal{F}_P^{\text{non-trac}}}  \bb{n}_{b} \cdot \left [ \sum_{q=1}^{q=N_{f,q}}\alpha_q \; (\bb{\sigma}_\text{s})_{b,q} \right]\Gamma_b
+ \sum_{b\, \in \, \mathcal{F}_P^{\text{trac}}}   \bb{n}_{b} \cdot \left[ \sum_{q=1}^{q=N_{f,q}}\bb{T}_{b,q} \,\Gamma_b\right],
\end{split}
\end{eqnarray}
where $\Gamma_P$ indicates the surface of cell $P$,
vector $\bb{T}_{b,g}$ represents the prescribed traction at the traction boundary quadrature point,
$N_{f,q}$ is the number of face quadrature points and
$\alpha_q$ is quadrature weight. Subscript ${f,q}$ indicates a quantity calculated at quadrature point location $\bb{x}_{f,q}$.
For non-triangular faces, face quadrature weights are scaled by the ratio of each triangle’s area to the total face area,
computed as the sum of all sub-triangle areas. Using a consistent definition of face area is essential here,
as inconsistency can degrade overall accuracy.
Additionally, a planar-face assumption is made by using a single outward unit normal \(\bb{n}_f\) for all quadrature points.
For non-planar faces, this geometric approximation can affects the accuracy and may limit the achievable convergence rate.

The stress tensor at a quadrature point $\bb{x}_{f,q}$, denoted $(\bb{\sigma}_\text{s})_{f,q}$, is computed as a function of the displacement gradient according to the adopted mechanical constitutive law.
For the case of small deformation engineering stress, it is given by:
\begin{eqnarray}
\bb{\sigma}_\text{s}(\bb{x}_{f,q}) =   \mu (\nabla_o \bb{u})_{f,q} + \mu (\nabla_o \bb{u})^T_{f,q} +\lambda \text{tr}\left[(\nabla_o \bb{u})_{f,q} \right]\bb{I},
\end{eqnarray}
where $\mu$ and $\lambda$ are the Lamé parameters.
Expressions for other constitutive laws are provided in Appendix \ref{app:mechLaws}.
The computation of the displacement gradient $(\nabla_o \bb{u})_{f,q}$ is shown in Section \ref{sec:grad}.
%
%------------------------------------------------------------------------------
\subsection{Volume Integrals}
\label{sec:vol_int}
\begin{figure}[h]
 	\centering
    \includegraphics[scale=1]{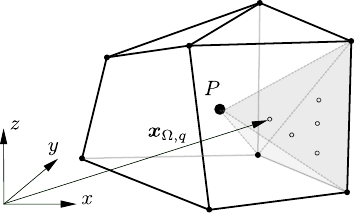}
 	\caption{Volume integration, cell decomposition into tetrahedral elements.}
 	\label{fig:volumeTerm}
\end{figure}

%------------------------------------------------------------------------------
%
The body force volume integral is discretised by integrating over the cell volume quadrature points:
\begin{eqnarray}\label{eq:bodyTerm}
\int_{\Omega_P} \bb{f}_b \,\text{d}\Omega_P
\approx
 \sum_{q=1}^{q=N_{\Omega,q}}\beta_q \, \left(\bb{f}_b\right)_{\Omega,q}  \Omega_P,
\end{eqnarray}
where $N_{\Omega,q}$ is the number of cell quadrature points and $\beta_q$ are the corresponding quadrature weights.
Subscript ${\Omega,q}$ is used to denote a quantity calculated at quadrature point location $\bb{x}_{\Omega,q}$.
Similarly, the inertia term (e.g. left-hand side term of Eq.~\eqref{eqn:momentum_lingeom}) is discretised as:
\begin{eqnarray}\label{eq:inertiaTerm}
\int_{\Omega_P} \rho \frac{\partial^2 \bb{u}}{\partial t^2} \text{d}\Omega_P
\approx
 \sum_{q=1}^{q=N_{\Omega,q}}\beta_q \; \bb{\rho}_{\,\Omega,q}  \left(\frac{\partial^2 \bb{u}}{\partial t^2} \right)_{\Omega,q}  \Omega_P.
\end{eqnarray}
To evaluate these integrals over arbitrary polyhedral cells, the control volume is partitioned into a set of non-overlapping tetrahedra using the face subdivisions already defined for surface integration.
Each tetrahedron is formed by connecting the cell centroid P with one of the triangular faces, as illustrated in Fig.~\ref{fig:volumeTerm}.
The cell partition is resolved at the quadrature weight calculation level, where the integration weights are scaled by the ratio of the host tetrahedron’s volume to the total cell volume. As in the computation of face quadrature weights, it is crucial to use a consistent cell volume \cite{Nishikawa2025}.
%
%------------------------------------------------------------------------------
\subsection{Stabilisation Scheme}
\label{sec:alpha}
%------------------------------------------------------------------------------
%
The evaluation of the surface term (diffusive flux) requires estimates of the solution gradient at the quadrature points $\bb{x}_{f,q}$ of a face.
The most common approach is to average the gradients reconstructed from each side of a quadrature point, either using arithmetic averaging \cite{Nishikawa2011,Nishikawa2025a} or interpolation with geometrically weighted coefficients \cite{Cardiff2025,Cardiff2017,Cardiff2018}.
In the present work, the solution gradient is evaluated directly at the quadrature points, eliminating the need for such interpolation.
However, as in the aforementioned averaging-based approaches, the discretisation of the diffusive flux may lead to numerical instabilities, manifested as zero-energy (checkerboard-type) modes.
To suppress these high-frequency oscillations, an appropriate stabilisation (damping) must be introduced.

A widely used stabilisation technique in the finite-volume framework is the Rhie–Chow type stabilisation \cite{Rhie1983}, which is second-order accurate \cite{Cardiff2025}.
Since this limits the overall convergence rate to second order, it is not suitable for the present high-order discretisation.
Instead, we employ the $\alpha$-stabilisation scheme \cite{Nishikawa2010}, which preserves the designed high-order convergence of the method \cite{Nishikawa2011, Nishikawa2025a}.
This stabilisation is introduced in the form of an additional traction term at each face $\bb{t}_{\alpha, f}$, defined as a function of the solution jump between the values extrapolated from the centroids of the two cells sharing a common face:
\begin{eqnarray}\label{eq:alpha}
\bb{t}_{\alpha, f} = \alpha \frac{\bar{K}}{|\bb{d}_{PN} \cdot \bb{n}_f|}(\bb{u}_{f,N}^*-\bb{u}_{f,P}^*),
\end{eqnarray}
where $\alpha$ is the scale factor and $\bar{K}$ a stiffness-type parameter that gives the stabilisation an appropriate scale and dimension.
In the current work, $\bar{K}$ is calculated using Lamé parameters $\bar{K} = 2\mu +\lambda$, following previous works \cite{Tukovic2013,Cardiff2017,Cardiff2018} where this value is used for Rhie-Chow stabilisation.
For hyperelastic constitutive laws, the values of $\mu$ and $\lambda$ correspond to the linearised (small-strain) moduli obtained by expanding the hyperelastic model about the undeformed configuration, ensuring that $\bar{K}$ is consistent with the model’s tangent stiffness in the small-strain limit
The vector $\bb{d}_{PN}=\bb{x}_N-\bb{x}_P$ connects the centroid of cell $P$ with that of its neighboring cell $N$,
while the product $\bb{d}_{PN} \cdot \bb{n}_f$ serves as a skewness measure, increasing the stabilisation on highly skewed cells \cite{Nishikawa2011}.
The quantities $\bb{u}_{f,N}^*$ and $\bb{u}_{f,P}^*$ denote the displacements at the face center $f$,
 extrapolated from the adjacent cell centers using the displacement gradients evaluated at $N$ and $P$, respectively:
\begin{eqnarray}\label{eq:alphaJump}
\begin{split}
\bb{u}_{f,N}^* &= \bb{u}_N + \bb{d}_{Nf} \cdot (\nabla \bb{u})_N + \frac{1}{2}\bb{d}_{Nf}^2 \!:\! (\nabla \nabla \bb{u})_N + \frac{1}{6}\bb{d}_{Nf}^3 \!::\! (\nabla \nabla \nabla \bb{u})_N, \\
\bb{u}_{f,P}^*& = \bb{u}_P + \bb{d}_{Pf} \cdot (\nabla \bb{u})_P + \frac{1}{2}\bb{d}_{Pf}^2 \!:\! (\nabla \nabla \bb{u})_P + \frac{1}{6}\bb{d}_{Pf}^3 \!::\! (\nabla \nabla \nabla \bb{u})_P,
\end{split}
\end{eqnarray}
where $\bb{d}_{Pf} = \bb{x}_f - \bb{x}_P$ and $\bb{d}_{Nf} = \bb{x}_f - \bb{x}_N$ represent the vectors from the cell centers to the face center.
The squared and cubed forms of the vector denote outer (tensor) products ($\bb{d}_{Nf}^3=\bb{d}_{Nf}\otimes \bb{d}_{Nf} \otimes  \bb{d}_{Nf}$), while the symbols $:$ and $\!::\!$ represent double and triple tensor contractions, respectively.
The second and third derivatives in Eq.~\eqref{eq:alphaJump} are employed only for third- and fourth-order discretisations, respectively, and their computation is described in Section \ref{sec:ho_scheme}.

The stabilisation contribution is added to the surface integral in Eq.~\eqref{eq:divStressDiscret}, for internal faces and boundary faces on which displacement is prescribed:
\begin{eqnarray}
\sum_{f \,\in\, \mathcal{F}_P^{\text{int}} \cup \mathcal{F}_P^{\text{disp}}} \bb{t}_{\alpha, f} \Gamma_f =
\sum_{f \,\in\, \mathcal{F}_P^{\text{int}}}  \alpha \frac{\bar{K}}{|\bb{d}_{PN} \cdot \bb{n}_f|}(\bb{u}_{f,N}^*-\bb{u}_{f,P}^*)\Gamma_f
+
\sum_{b \,\in\, \mathcal{F}_P^{\text{disp}}}  \alpha \frac{\bar{K}}{|\bb{d}_{Pb} \cdot \bb{n}_b|}(\bb{u}_{b}-\bb{u}_{f,P}^*)\Gamma_b,
\end{eqnarray}
where $\bb{u}_{b}$ is the prescribed displacement and $\bb{d}_{Pb}=\bb{x}_b - \bb{x}_P$. Note that the stabilisation term is integrated using a single quadrature point per face, as there is no requirement for high accuracy in evaluating the stabilisation contribution; it is only important that the added term converges with the desired order.
%
%------------------------------------------------------------------------------
\subsection{High order interpolation scheme}
\label{sec:ho_scheme}
%------------------------------------------------------------------------------
%
For an arbitrary point in space, denoted by $\tilde{\bb{x}}$ (for example, a face quadrature point), the displacement field can be reconstructed using a pointwise interpolation of the surrounding cell-centre values:
\begin{eqnarray}\label{eq:interp}
(\bb{u})_{\tilde{\bb{x}}} = \sum_{N \in \mathcal{S}}c_N \bb{u}_N,
\end{eqnarray}
where $c_N$ are scalar interpolation coefficients associated with the neighbouring cell-centre (computational) nodes $\bb{u}_N=(u_{x,N}, u_{y,N}, u_{z,N})^{\mathrm{T}}$, and $\mathcal{S}$ denotes the interpolation stencil, i.e. the set of computational nodes used in the interpolation, with $|\mathcal{S}|$ being its size.
The same interpolation can be used to evaluate the first and higher spatial derivatives of the displacement field:
\begin{eqnarray}\label{eq:derInterp}
\left(\frac{\partial \bb{u}}{\partial x}\right)_{\tilde{\bb{x}}}=\sum_{N \in \mathcal{S}}c_{x_N} \bb{u}_N, \qquad
\left(\frac{\partial \bb{u}}{\partial y}\right)_{\tilde{\bb{x}}}=\sum_{N \in \mathcal{S}}c_{y_N} \bb{u}_N, \qquad
\left(\frac{\partial \bb{u}}{\partial z}\right)_{\tilde{\bb{x}}}=\sum_{N \in \mathcal{S}}c_{z_N} \bb{u}_N,
\end{eqnarray}
where $c_{x_N}$, $c_{y_N}$ and $c_{z_N}$ are the interpolation coefficients corresponding to the spatial derivatives in the $x-$, $y-$, and $z-$directions, respectively.
The procedure for computing these interpolation coefficients is described in the following section.
%
%------------------------------------------------------------------------------
\subsubsection{Least-Squares Taylor Reconstruction}
%------------------------------------------------------------------------------
%
The interpolation coefficients in Eqs.~\eqref{eq:interp}~and~\eqref{eq:derInterp} are determined using a local least-squares reconstruction method, commonly referred to as the Local Regression Estimator (LRE) in \cite{Castrillo2022, Castrillo2023, Castrillo2024}. A truncated Taylor expansion of the component $\alpha$ of the displacement field, denoted by $\bar{u}_\alpha$, in the vicinity of $\tilde{\boldsymbol{x}}$ is written as:
\begin{eqnarray}
\begin{split}
\bar{u}_\alpha(\bb{x}) =
(u_\alpha)_{\tilde{\bb{x}}}
 + \left(\frac{\partial u_\alpha}{\partial x}\right)_{\tilde{\bb{x}}}(x-\tilde{x})&
 + \left(\frac{\partial u_\alpha}{\partial y}\right)_{\tilde{\bb{x}}}(y-\tilde{y})
 + \left(\frac{\partial u_\alpha}{\partial z}\right)_{\tilde{\bb{x}}}(z-\tilde{z})
 + \frac{1}{2}\left(\frac{\partial^2 u_\alpha}{\partial x^2}\right)_{\tilde{\bb{x}}}(x-\tilde{x})^2 \\
 +  \left(\frac{\partial^2 u_\alpha}{\partial x \partial y}\right)_{\tilde{\bb{x}}}&(x-\tilde{x})(y-\tilde{y})
 +  \left(\frac{\partial^2 u_\alpha}{\partial x \partial z}\right)_{\tilde{\bb{x}}}(x-\tilde{x})(z-\tilde{z})
 + \dots
\end{split}
\label{eq:te}
\end{eqnarray}
%
%The expansion is truncated at polynomial order $p$, resulting in $N_t$ terms in the polynomial basis.
%Number of terms $N_t$ depends on spatial dimensions $d$ and polynomial order $p$, $N_t = \dfrac{(d+p)!}{(d!\,p!)}$.
The expansion is truncated at polynomial order $p$, resulting in $N_t$ terms in the polynomial basis.
The number of terms $N_t$ depends on the spatial dimension $d$ and the polynomial order $p$, and is given by
 $N_t = \dfrac{(d+p)!}{(d!\,p!)}$ \cite{Khelladi2011}.
%The expansion is truncated at polynomial order $p$, which corresponds to $N_t$ terms, where $N_t = \dfrac{(d+p)!}{(d!\,p!)}$ in $d$ spatial dimensions. 
Although the interpolation formulation supports arbitrary polynomial order $p$, in this work we restrict attention to orders up to $p = 3$.
The truncated expansion can then be written compactly using a polynomial basis $\mathbf{q}$ and a vector of unknown coefficients $\bar{\mathbf{a}}$ representing the derivatives:
\begin{eqnarray}
\bar{u}_\alpha(\bb{x}) = \mathbf{q}^{\mathrm{T}}(\bb{x}-\tilde{\bb{x}})\bar{\mathbf{a}}(\tilde{\bb{x}}),
\end{eqnarray}
where
\begin{eqnarray}\label{eq:qAnda}
\begin{split}
&\mathbf{q}^{\mathrm{T}}(\bb{x}-\tilde{\bb{x}}) = \left[1, (x-\tilde{x}), (y-\tilde{y}), (z-\tilde{z}), \frac{1}{2}(x-\tilde{x})^2, \dots \right],\\
& \bar{\mathbf{a}}^{\mathrm{T}}(\tilde{\bb{x}}) = \left[(u_\alpha)_{\tilde{\bb{x}}}, \left(\frac{\partial u_\alpha}{\partial x}\right)_{\tilde{\bb{x}}},
 \left(\frac{\partial u_\alpha}{\partial y}\right)_{\tilde{\bb{x}}},
 \left(\frac{\partial u_\alpha}{\partial z}\right)_{\tilde{\bb{x}}},
 \left(\frac{\partial^2 u_\alpha}{\partial x^2}\right)_{\tilde{\bb{x}}},
 \dots  \right].
\end{split}
\end{eqnarray}
The parameter vector $\bar{\mathbf{a}}$ is calculated by minimising the weighted sum of squared differences between the Taylor approximation and the cell-centre values:
\begin{eqnarray}\label{eq:ls}
\mathcal{R} = \frac{1}{2} \sum_{N \in \mathcal{S}} w_N \left [ \bar{u}_\alpha(\bb{x}_N) - u_{\alpha,N}\right]^2 =\frac{1}{2} \sum_{N \in \mathcal{S}} w_N \left [ \mathbf{q}^{\mathrm{T}}(\bb{x}_N-\tilde{\bb{x}})\bar{\mathbf{a}}(\tilde{\bb{x}}) - u_{\alpha,N}\right]^2.
\end{eqnarray}
Minimizing Eq.~\eqref{eq:ls} with respect to $\bar{\mathbf{a}}(\tilde{\boldsymbol{x}})$ leads to the following equation:
\begin{eqnarray}\label{eq:maQw}
\mathbf{{M}}(\tilde{\bb{x}}) \mathbf{\bar{a}}(\tilde{\bb{x}}) = \mathbf{{Q}}(\tilde{\bb{x}}) \mathbf{{W}}(\tilde{\bb{x}}) \mathbf{u}_\alpha,
\end{eqnarray}
where $\mathbf{M}(\tilde{\bb{x}})=\mathbf{Q}(\tilde{\bb{x}})\mathbf{W}(\tilde{\bb{x}})\mathbf{Q}^{\mathrm{T}}(\tilde{\bb{x}})$, $\mathbf{Q}(\tilde{\bb{x}})$ is the $N_t \times |\mathcal{S}|$ matrix whose $N-$th column is $\mathbf{q}^{\mathrm{T}}(\bb{x}_N-\tilde{\bb{x}})$, $\mathbf{W}(\tilde{\bb{x}})$ is a diagonal weighting matrix of size $|\mathcal{S}| \times |\mathcal{S}|$ and $\mathbf{u}^\text{T}_\alpha=\left[u_{\alpha,1}, u_{\alpha,2},\dots,u_{\alpha,|\mathcal{S}|} \right] $.
Although Eq.~\eqref{eq:maQw} represents the normal-equation form of the least-squares problem, it is not solved directly; instead, a QR factorization via Householder transformations \cite{strang2012} is used to obtain the coefficients in a more robust way.
Defining $\bar{\mathbf{A}}(\tilde{\bb{x}})=\mathbf{M}(\tilde{\bb{x}})^{-1}\mathbf{Q}(\tilde{\bb{x}})\mathbf{W}(\tilde{\bb{x}})$, the system becomes:
\begin{eqnarray}\label{eq:lsA}
\mathbf{\bar{a}}(\tilde{\bb{x}}) = \mathbf{\bar{A}}(\tilde{\bb{x}})\mathbf{u}_\alpha.
\end{eqnarray}
The rows of $\bar{\mathbf{A}}(\tilde{\boldsymbol{x}})$ correspond to the interpolation coefficients used in
Eqs.~\eqref{eq:interp}~and~\eqref{eq:derInterp} depending on order they were defined in Eq.~\eqref{eq:qAnda}.
The interpolation coefficients from matrix $\bar{\mathbf{A}}(\tilde{\boldsymbol{x}})$ depend only on the mesh (geometry), so they can be computed once (or whenever the mesh moves) and stored.
The inversion of ${\mathbf{M}}(\tilde{\boldsymbol{x}})$ may be ill-conditioned, for poor spatial distribution of stencil nodes.
This issue is mitigated by scaling the Taylor basis, ensuring sufficient stencil size and restricting the weights in $\mathbf{W}(\tilde{\boldsymbol{x}})$ to positive values, see \cite{Castrillo2022, Castrillo2023}.
The weight function in Eq.~\eqref{eq:ls} is chosen here as a radially symmetric exponential function:
\begin{eqnarray}\label{eq:weightFunc}
w_N = \dfrac{\exp\left( -\tilde{d}^2k^2\right) -\exp\left(-k^2\right)}{1-\exp\left(-k^2\right)},
\end{eqnarray}
where $d=|\bb{\tilde{x}}-\bb{x}_N|$, $d_{\mathrm{s}}=2r_{\mathrm{s}}$ is the stencil diameter, $\tilde{d}=d/d_{\mathrm{s}}$ and $r_{\mathrm{s}}$ is the stencil radius, defined as the maximum distance $|\bb{\tilde{x}}-\bb{x}_N|$ in stencil $\mathcal{S}$. The weighting function in Eq.~\eqref{eq:weightFunc} has been previously used in the context of the Navier–Stokes equations
\cite{Khelladi2011, Ramirez2014} and solid mechanics problems \cite{Castrillo2022, Castrillo2023, Castrillo2024}.
Following these studies, the shape parameter $k = 6$ is adopted herein for both 2D and 3D cases, as it was shown to provide optimal performance \cite{Castrillo2022, Castrillo2023, Castrillo2024}.
%
%------------------------------------------------------------------------------
\subsubsection{Interpolation Stencil}
\label{sec:stencil}
%------------------------------------------------------------------------------
%
The minimum stencil size $|\mathcal{S}|_{\text{min}}$ is determined by the dimension of the polynomial basis and the interpolation order, i.e. $|\mathcal{S}|_{\text{min}} = N_t$:
\begin{eqnarray}
\begin{split}
\text{for 2D:}\quad & |\mathcal{S}|_{\text{min}} = \frac{(p+1)(p+2)}{2}, \\
\text{for 3D:}\quad & |\mathcal{S}|_{\text{min}} =\frac{(p+1)(p+2)(p+3)}{6}.
\end{split}
\end{eqnarray}
The actual stencil size $|S|$ used for computations is set as:
\begin{eqnarray}
|\mathcal{S}| = |\mathcal{S}|_{\text{min}} + n^+
\end{eqnarray}
where $n^{+}$ is a user-defined number of surplus nodes. The choice of $n^{+}$ represents a trade-off: an excessively large stencil introduces numerical dissipation and increases computational cost, whereas a small stencil may lead to unstable interpolation \cite{Khelladi2011}.
Based on those findings presented in \cite{Castrillo2022,Castrillo2024}, $n^{+}$ is set to $10$ for all two-dimensional cases, and to $45$, $55$, and $65$ for $p = 1$, $p = 2$, and $p = 3$ in three-dimensional cases, respectively. These values are consistent with those reported in the literature. In particular, \cite{Tsoutsanis2023} noted that within the $k-$exact least-squares framework, the stencil size is typically chosen in the range $1.5\,|\mathcal{S}|_{\text{min}} \leq |\mathcal{S}| \leq 3\,|\mathcal{S}|_{\text{min}}$.

\begin{figure}[h]
 	\centering
	\subfigure[Stencil $\mathcal{S}_f$ for internal face $f$ quadrature  \newline  points $\bb{x}_{f,q}$]
 	{
 		\label{fig:faceStencil}
    	\includegraphics[scale=1]{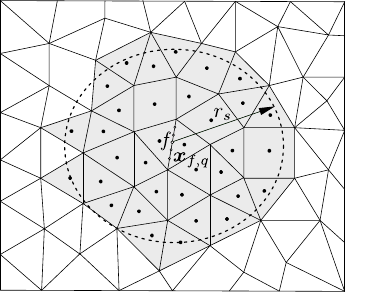}
    }\hspace*{-0.1cm}
 	\subfigure[Stencil for boundary face $b$ quadrature points $\bb{x}_{b,q}$ in the case of \newline
 	displacement $\mathcal{S}^{\text{disp}}_b$ and symmetry $\mathcal{S}^{\text{symm}}_b$ boundary type. ]
 	{
 		\label{fig:faceBoundaryStencil}
    	\includegraphics[scale=1]{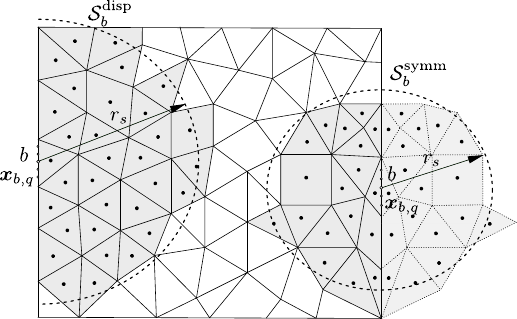}
    }
 	\caption{Interpolation stencils for internal $f$ and boundary $b$ face quadrature points, $\bb{x}_{f,q}$ and $\bb{x}_{b,q}$, respectively.}
 	\label{fig:stencil}
\end{figure}

The stencil is constructed using a nearest-neighbour approach, where the set $\mathcal{S}$ comprises the cells closest to the reference point (face or cell centre), enclosed within a sphere of radius $r_s$. To reduce computational cost, all quadrature points on a given face share the same stencil, and likewise, all quadrature points within a cell use a common stencil. The construction of the internal-face stencil $\mathcal{S}_f$ is illustrated in Fig.~\ref{fig:faceStencil}.

Stencil construction is adapted for boundary faces according to the imposed boundary condition (see Fig.~\ref{fig:faceBoundaryStencil}):
\begin{itemize}
\item Traction boundary – no stencil is defined, since interpolation is not applied.
\item Displacement boundary – the stencil, denoted $\mathcal{S}_b^{\text{disp}}$, is constructed in the same manner as for internal faces.
\item Symmetry boundary – the stencil, denoted $\mathcal{S}_b^{\text{symm}} := \mathcal{S}_b^{\text{symm,p}} \cup \mathcal{S}_b^{\text{symm,g}}$, is divided equally: half of the nodes correspond to the physical domain $\mathcal{S}_b^{\text{symm,p}}$, while the remaining half are mirrored ghost nodes $\mathcal{S}_b^{\text{symm,g}}$ that mirror both the position and state variables.
\end{itemize}

As an alternative, the stencil can be constructed by selecting layers of neighbouring cells around the reference point; however, this approach can lead to less predictable stencil sizes. More advanced stencil construction strategies, such as those that account for cell distribution and moment-matrix conditioning, may also be employed \cite{Khelladi2011,Tsoutsanis2023}, but a detailed assessment of these approaches in the context of solid mechanics is required and is left for future work.
%
%
%------------------------------------------------------------------------------
\subsubsection{Boundary Conditions}
\label{sec:bc}
%------------------------------------------------------------------------------
%
Boundary conditions are imposed through modifications of the local polynomial reconstruction, except for traction boundaries,
where the known tractions at the quadrature points are directly integrated and added to the surface term in Eq.~\eqref{eq:divStressDiscret}.

At displacement boundaries, the prescribed displacement is enforced in a weak sense by introducing a penalty contribution into the polynomial reconstruction.
This is achieved by extending the stencil of each boundary quadrature point with an additional ghost nodes, co-located with the quadrature point itself and assigned the prescribed displacement value $\mathbf{u}_{b,q}$.
The inclusion of this ghost node effectively augments the reconstruction matrices $\mathbf{W}$ and $\mathbf{Q}$.
Since the Taylor expansion in Eq.~\eqref{eq:te} is centred at the same location, the associated polynomial basis vector simplifies to
$\mathbf{q}^{\mathrm{T}} = [1,\,0,\,0,\,\dots,\, 0]$, and its weight in $\mathbf{W}$ is typically set to unity.
Increasing this weight enhances the local satisfaction of the prescribed displacement but does not affect the overall convergence rate.

At symmetry boundaries, only the half of the stencil located inside the physical domain is available, and the remaining nodes are generated by reflecting this interior stencil across the symmetry plane. The polynomial basis vectors for these mirrored nodes are computed using their reflected coordinates. This approach guarantees enforcement of symmetry constraints; zero tangential component of boundary traction $(\bb{I}-2\bb{n}_b\otimes\bb{n}_b)\cdot \bb{T}_{b,q}= 0$, and zero normal component of displacement $\bb{n}_b \cdot \bb{u}_b = 0$.

Achieving high-order convergence near boundaries requires careful attention to geometric representation. When curved surfaces are approximated by planar faces, the resulting geometric discretisation is second-order accurate, which inherently limits the overall convergence rate of the method. High-order accuracy can be recovered either by introducing isoparametric boundary representations \cite{Coppeans2026} or by maintaining flat faces while introducing a correction \cite{costa2021, Costa2022, Gooch2002}. Curved boundary treatment has not yet been tackled in literature (for finite volume solid mechanics) and is therefore left for future work, as it presents a challenging problem deserving dedicated study.
%
%------------------------------------------------------------------------------
\subsubsection{Computing Gradients}
\label{sec:grad}
%------------------------------------------------------------------------------
%
For each quadrature point on internal faces, Eq.~\eqref{eq:lsA} is evaluated and the resulting interpolation coefficients from the matrix $\mathbf{\bar{A}}$ are stored for computation of the displacement gradient:
\begin{eqnarray}\label{eq:faceQuadGrad}
(\nabla \bb{u})_{f,q} = \sum_{N \,\in\,  \mathcal{S}_f}  \bb{c}_{\bb{x}_N} \otimes \bb{u}_N,
\end{eqnarray}
where $\bb{c}_{\bb{x}_N}^{\mathrm{T}} = \left[ \,c_{x_N}, \; c_{y_N},\; c_{z_N} \right]$ is the interpolation vector (see Eq.~\eqref{eq:derInterp}).
The same procedure is applied at displacement boundaries using the stencil $\mathcal{S}_b^{\text{disp}}$  extended by the ghost node:
\begin{eqnarray}
(\nabla \bb{u})_{b,q} = \left(\sum_{N \,\in\,  \mathcal{S}_b^{\text{disp}}}  \bb{c}_{\bb{x}_N} \otimes \bb{u}_N \right)+ \bb{c}_{\bb{b}_N} \otimes \bb{u}_{b,q},
\end{eqnarray}
where $\bb{u}_{b,q}$ is prescribed displacement at the boundary quadrature point. Coefficient vector $\bb{c}_{\bb{b}_N}$ is extracted from the corresponding location from the same matrix $\mathbf{\bar{A}}$ as for $\bb{c}_{\bb{x}_N}$ coefficients.
For a symmetry boundary, the contribution of the mirrored (ghost) nodes is also included, giving:
\begin{eqnarray}
(\nabla \bb{u})_{b,q} =
\sum_{N \,\in\,  \mathcal{S}_b^{\text{symm,p}}}  \bb{c}_{\bb{x}_N} \otimes \bb{u}_N +
\sum_{N \,\in\,  \mathcal{S}_b^{\text{symm,g}}}  \bb{c}_{\bb{x}_N} \otimes ( \bb{R}_b \cdot \bb{u}_N),
\end{eqnarray}
where $ \bb{R}_b = (\bb{I}-2\bb{n}_b \otimes \bb{n}_b)$ is the reflection tensor \cite{Demirdzic2022}.
Note that only first derivatives are stored at face quadrature points, as higher derivatives are not required by the constitutive equations.
%In the present solution procedure, the displacements $\bb{u}_N$ are known, making the gradient evaluation straightforward.
%However, if $\bb{u}_N$ are treated as unknowns, the interpolation coefficients must be incorporated into the left-hand-side system matrix, as %shown in \cite{Castrillo2022,Castrillo2024} and Appendix \ref{app:P}.

The stabilisation scheme requires the evaluation of first and higher derivatives at cell centres when computing Eq.~\eqref{eq:alphaJump}.
For this purpose, each cell is associated with a stencil $\mathcal{S}_P$, constructed in the same manner as the face stencil.
At each cell centre, the gradient is computed according to Eq.~\eqref{eq:faceQuadGrad}, while the second derivative is evaluated for $p = 2$, and the third derivative tensor for $p = 3$.
The formulation and computation of the second and the third derivative tensors are described in Appendix \ref{app:derivatives}.
%
%%%%%%%%%%%%%%%%%%%%%%%%%%%%%%%%%%%%%%%%%%%%%%%%%%%%%%%%%%%%%%%%%%%%%%%%%%%%%%%
%
\section{Solution Algorithm}
\label{sec:Sol_Alg}
%
%%%%%%%%%%%%%%%%%%%%%%%%%%%%%%%%%%%%%%%%%%%%%%%%%%%%%%%%%%%%%%%%%%%%%%%%%%%%%%%
%
The nonlinear system arising from the finite volume discretisation is solved using the Jacobian-free Newton–Krylov (JFNK) algorithm.
The subsequent subsection describes the numerical formulation and implementation aspects of this approach.
%
%------------------------------------------------------------------------------
\subsection{Jacobian-free Newton-Krylov Algorithm}
\label{sec:JFNK_alg}
%------------------------------------------------------------------------------
%
The linear momentum balance (Eqs.~\eqref{eqn:momentum_TL}~and~\eqref{eqn:momentum_lingeom}) can be expressed in residual form as:
\begin{eqnarray}
\bb{R}(\bb{u}) = 0,
\end{eqnarray}
where $\bb{u}$ is the primary unknown field (displacement) and $\bb{R}$ represents the \emph{residual} (imbalance) of the equation.
For the implicit discretisation of the momentum balance, the residual vector is a function of the unknown displacement field at iteration $k+1$:
\begin{eqnarray}
\bb{R}(\bb{u}_{k+1}) = 0.
\end{eqnarray}
The residual can be approximated using a first-order Taylor expansion about the current iteration $k$:
\begin{eqnarray}
 	\bb{R}(\bb{u}_{k+1}) \approx \bb{R}(\bb{u}_{k}) \;+\;  \bb{R}'(\bb{u}_{k}) (\bb{u}_{k+1} - \bb{u}_{k}),
\end{eqnarray}
where $\bb{R}'$ is the Jacobian matrix $\bb{J}\equiv\bb{R}'\equiv\partial \bb{R} / \partial \bb{u}$.
Imposing $\bb{R}(\bb{u}_{k+1}) = \bb{0}$ yields the standard Newton update equation:
\begin{eqnarray}
\bb{J}(\bb{u}_k) \delta \bb{u} = -\bb{R}(\bb{u}_k),
\end{eqnarray}
where $\delta \bb{u}=\bb{u}_{k+1} - \bb{u}_{k}$.
In the Jacobian-free Newton–Krylov approach, the Jacobian matrix is not formed explicitly.
Instead, a Krylov subspace method (e.g., GMRES) is used to solve the linear system iteratively, requiring only the action of the Jacobian on a arbitrary vector $\bb{v}$:
\begin{eqnarray}
\bb{J} \bb{v} \approx \frac{\bb{R}(\bb{u}+\epsilon \bb{v})-\bb{R}(\bb{u})}{\epsilon},
\end{eqnarray}
where $\epsilon$ is a small scalar perturbation.
To mitigate the effect of ill-conditioned Jacobian matrices on convergence and robustness, a left preconditioning strategy is employed:
\begin{eqnarray}\label{eq:jfnk}
 	\bb{P}^{-1}\bb{J}  \bb{v}\approx 	\frac{\bb{P}^{-1}(\bb{R}(\bb{u} + \epsilon  \bb{v}) - \bb{R}(\bb{u}))}{\epsilon},
\end{eqnarray}
where $\bb{P}$ denotes the preconditioning matrix.
As the preconditioning matrix $\bb{P}$, we adopt a compact-stencil approximate Jacobian, the one that is used in the second-order cell-centre semi-implicit discretisation with segregated solution procedure \cite{Cardiff2018, Tukovic2018}. Following the classification of Knoll and Keyes \cite{Knoll2004}, this preconditioning choice may be considered ``physics-based".

The approximate Jacobian, for cell $P$, corresponds to the discretised diffusion term using a central-difference scheme:
\begin{eqnarray}
\begin{split}
 	\bb{P}_P = \oint_{\Gamma_P} \bar{K} \, \bb{n} \cdot \bb{\nabla} \bb{u} \; \text{d}\Gamma_P \approx
 		\sum_{f \in \mathcal{F}^{\text{int}}_P}  \bar{K} &
 		\left|\bb{\Delta}_{f} \right| \frac{ \bb{u}_{N} - \bb{u}_P}{\left|\bb{d}_{PN}\right|}    \Gamma_{f}
 	+  \sum_{b \in \mathcal{F}^{\text{disp}}_P}  \bar{K}
 		\left|\bb{\Delta}_{b} \right| \frac{ \bar{\bb{u}}_{b}  - \bb{u}_P}{\left|\bb{d}_{Pb}\right|}
 		{\Gamma}_{b}    \\
 	&	+ \sum_{b \in \mathcal{F}^{\text{symm}}_P}  \bar{K}
 		\left|\bb{\Delta}_{b} \right| \frac{ \bb{R}_{b} \cdot \bb{u}_{P} - \bb{u}_P}{\left|\bb{d}_{Pb}\right|} {\Gamma}_{b},
\end{split}
\end{eqnarray}
where $\bb{\Delta}_{f}=\bb{d_{PN}}/(\bb{d}_{PN} \cdot \bb{n}_f)$ is \emph{over-relaxed orthogonal} vector \cite{Demirdzic2015}.
Note that the non-orthogonal correction, which is normally required to preserve accuracy when discretising the diffusion term on skewed meshes, is omitted here.
Including it would enlarge the stencil (molecule) and compromise the compactness achieved with the central-difference approximation.
As will be shown later, this simplification does not degrade robustness or accuracy, since the overall accuracy of the discretisation is governed entirely by the residual evaluation itself.

For the JFNK solver, the PETSc library (version 3.22.2) \cite{PETSc} is used, where the interfaces for evaluating the approximate Jacobian (matrix $\bb{P}$) and the high-order residual ($\boldsymbol{R}(\bb{u})$) are implemented within the \texttt{solids4foam} toolbox \cite{Cardiff2018, Tukovic2018, Cardiff2025a}, which is built on OpenFOAM-v2312 \cite{Weller1998}.
An extended discussion with additional implementation details is provided in the authors’ previous work \cite{Cardiff2025}, where the JFNK solver is coupled with a second-order residual evaluation.
%
%%%%%%%%%%%%%%%%%%%%%%%%%%%%%%%%%%%%%%%%%%%%%%%%%%%%%%%%%%%%%%%%%%%%%%%%%%%%%%%
%
\section{Test Cases}
\label{sec:test_cases}
%
%%%%%%%%%%%%%%%%%%%%%%%%%%%%%%%%%%%%%%%%%%%%%%%%%%%%%%%%%%%%%%%%%%%%%%%%%%%%%%%
%
The performance of the proposed high-order solver is verified using several benchmark cases.
These are employed to assess both the solution accuracy and the formal order of convergence, and to provide a direct comparison against a standard second-order discretisation.
The comparison is evaluated in terms of solution accuracy, computational cost, and memory requirements.
In addition, the influence of the stabilisation factor, different choices of preconditioning matrices, and other numerical aspects of the method are examined.

The benchmark cases are selected to cover a broad range of features, including both two- and three-dimensional problems, small- and large-strain regimes,
and various material models.
The convergence behaviour is studied on a sequence of mesh refinements to detect potential oscillatory convergence trends, which may arise in high-order finite-volume formulations \cite{Nishikawa2025a}.

In all simulations, the convergence criterion corresponds to a six-order-of-magnitude reduction in the residual, and a line-search procedure is employed driving the Newton-step solution update \cite{Cardiff2025}.
The number of quadrature points per tetrahedral, triangular, or line element is kept to the minimum required for exact integration of the displacement gradient $\nabla \bb{u}$, from which the stresses are computed.
For example, triangular elements with $p=1$ and $p=2$ interpolation require one integration point, whereas $p=3$ uses three.
The interpolation stencil size for all cases is used as specified in Section~\ref{sec:stencil}.
Different mesh types are considered: the open-source mesher Gmsh \cite{geuzaine2009gmsh} was employed to generate tetrahedral meshes. 
For the structured hexahedral meshes, the native OpenFOAM utilities \texttt{blockMesh} and \texttt{extrudeMesh} were used. 
Polyhedral meshes were subsequently produced by converting a tetrahedral mesh using the \texttt{polyDualMesh} OpenFOAM utility.
%
% Aim of each test case
%
% MMS 3D - to show convergence behaviour
% MMS 2D - to show behavior in 2D, p=2 behaves like p=1 also in 2D
%
% Cantilever - 1. machine precistion reached if underlying solution is of same order as interpolation function. 
%              2. shear locking effect
%
%------------------------------------------------------------------------------
\subsection{Testing the Accuracy and Order of Accuracy}
%------------------------------------------------------------------------------
%
The accuracy and convergence behaviour of the method are evaluated for interpolation orders $p = 1$, $p = 2$, and $p = 3$.
For cases where an analytical solution is available, the solution error is quantified using the $L_2$ and $L_\infty$ norms,
where the $L_2$ norm represents the root mean square error:
\begin{eqnarray}
L_2 = \sqrt{\frac{1}{N_{\text{c}}}\sum_{i=1}^{N_{\text{c}}} \Delta \phi_i^2},
\end{eqnarray}
and $L_\infty$ is infinity norm representing the maximum absolute error:
\begin{eqnarray}
L_{\infty} = \max_{1 \leq i \leq N_{\text{c}}} |\Delta \phi_i|,
\end{eqnarray}
where $\Delta \phi_i$ is the difference between expected and predicted solutions at computational nodes and $N_{\text{c}}$ is the overall number of computational nodes, i.e. cells.
Both norms are calculated for displacement magnitude $|\boldsymbol{u}|=\sqrt{u_x^2+u_y^2+u_z^2}$ and von Mises stress:
\begin{eqnarray}
\sigma_\text{eq}=\sqrt{\frac{1}{2}\left[  (\sigma_{xx}-\sigma_{yy})^2+(\sigma_{yy}-\sigma_{zz})^2+(\sigma_{zz}-\sigma_{xx})^2 \right]+3\left(\sigma_{xy}^2+\sigma_{yz}^2+\sigma_{zx}^2\right)}.
\end{eqnarray}
%
%------------------------------------------------------------------------------
\subsubsection{Case 1: Order Verification via the Manufactured Solution Procedure}
\label{case:mms}
%------------------------------------------------------------------------------
%
The first test case consists of a  $0.2\,\text{m} \times 0.2$ m square or $0.2\,\text{m} \times 0.2\,\text{m} \times 0.2$ m cube with linear elastic (Young's modulus of $E = 200$ GPa and Poisson's ratio of $\nu = 0.3$) properties.
A manufactured solution for displacement, prescribed on domain boundaries, is employed in the form or trigonometric function in 3D \citep{Mazzanti2024} and mixed form in 2D \citep{Castrillo2022}:
\begin{itemize}
\item[•] \textbf{2D}
\begin{eqnarray}
	\bb{u} =
	\begin{pmatrix}
 e^{x^2}\sin(y)\\
\ln(3+y)\cos(x)+\sin(y)\\
 0
	\end{pmatrix},
\end{eqnarray}
\item[•] \textbf{3D}
 \begin{eqnarray}
	\bb{u} =
	\begin{pmatrix}
	a_x \sin(4\pi x) \sin(2\pi y) \sin(\pi z) \\
	a_y \sin(4 \pi x) \sin(2 \pi y) \sin(\pi z) \\
	a_z \sin(4 \pi x) \sin(2 \pi y) \sin(\pi z)
	\end{pmatrix},
\end{eqnarray}
where $a_x = 2\,\mu$m, $a_y = 4\,\mu$m, and $a_z = 6\,\mu$m.
\end{itemize}
The Cartesian coordinates are given by $x$, $y$ and $z$.
The corresponding manufactured body force term ($\bb{f}_b$ in Eq.~\eqref{eqn:momentum_lingeom}) can be obtained by the manufactured solution procedure.

Both manufactured-solution cases were examined using four mesh types: (i) \emph{regular} hexahedral meshes, (ii) \emph{regular} tetrahedral meshes, (iii) \emph{regular} polyhedral meshes, and (iv) \emph{irregular} tetrahedral meshes. In two dimensions, these correspond to quadrilateral, triangular, and polygonal meshes. For reference, the finest 2D quadrilateral mesh contains $1\;024$ cells, while the finest irregular triangular mesh contains $2\;048$ cells. In 3D, the finest regular hexahedral mesh consists of $91\;125$ cells, whereas the finest regular tetrahedral mesh contains $511\;104$ cells.

\begin{figure}[H]
 	\centering
	\subfigure[Regular quadrilateral mesh]
 	{
 	    \begin{overpic}[scale=0.8]{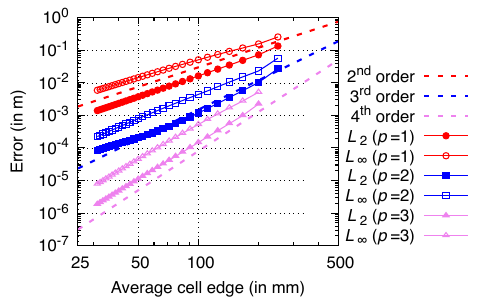}
  			\put(50,13){\includegraphics[scale=0.016]{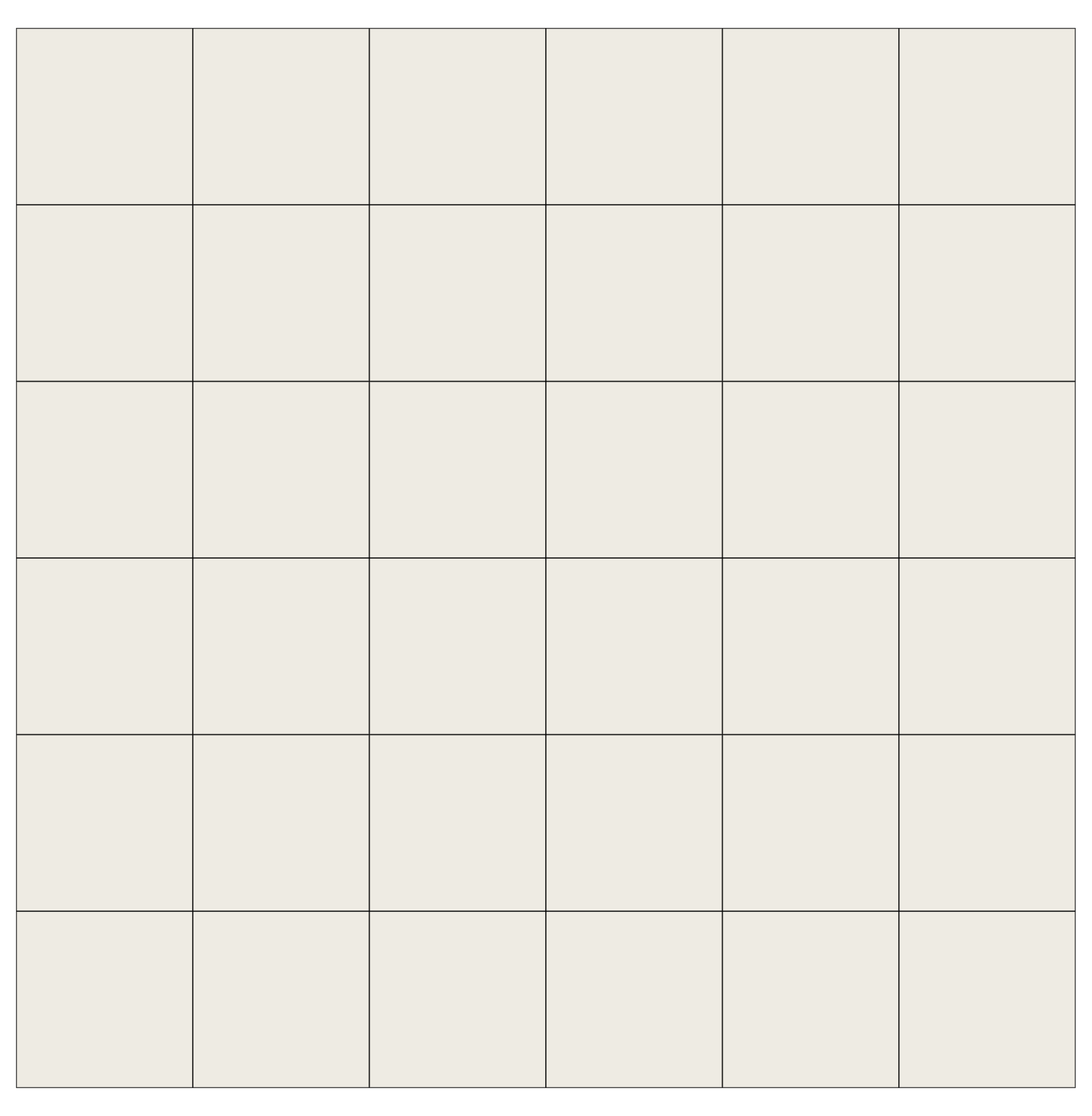}}
    	\end{overpic}
    }
 	\subfigure[Regular triangular mesh]
 	{
 	 	\begin{overpic}[scale=0.8]{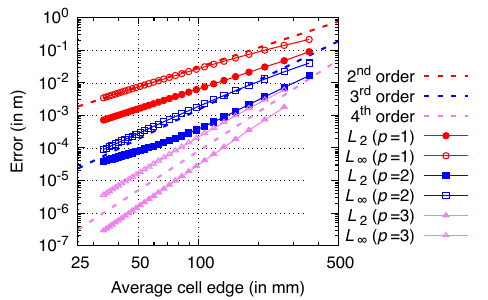}
  			\put(50,13){\includegraphics[scale=0.016]{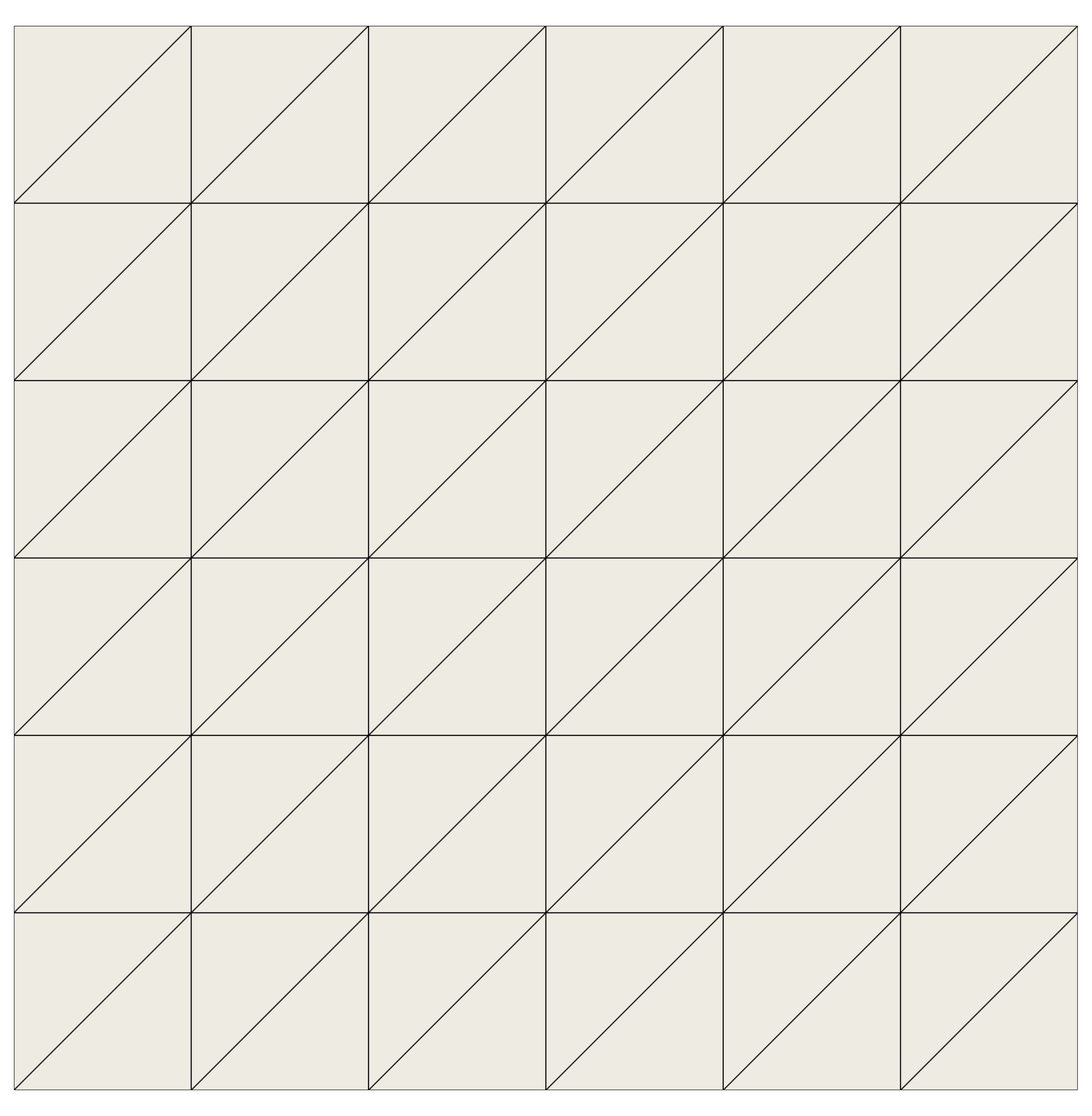}}
    	\end{overpic}
    }
    \subfigure[Irregular triangular mesh]
 	{
 	 	\begin{overpic}[scale=0.8]{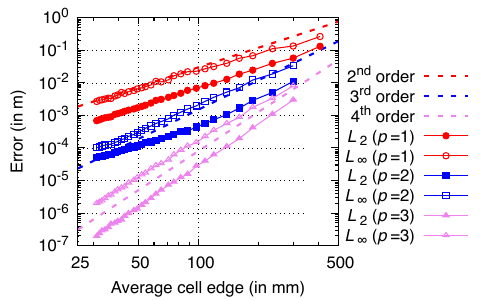}
  			\put(50,13){\includegraphics[scale=0.016]{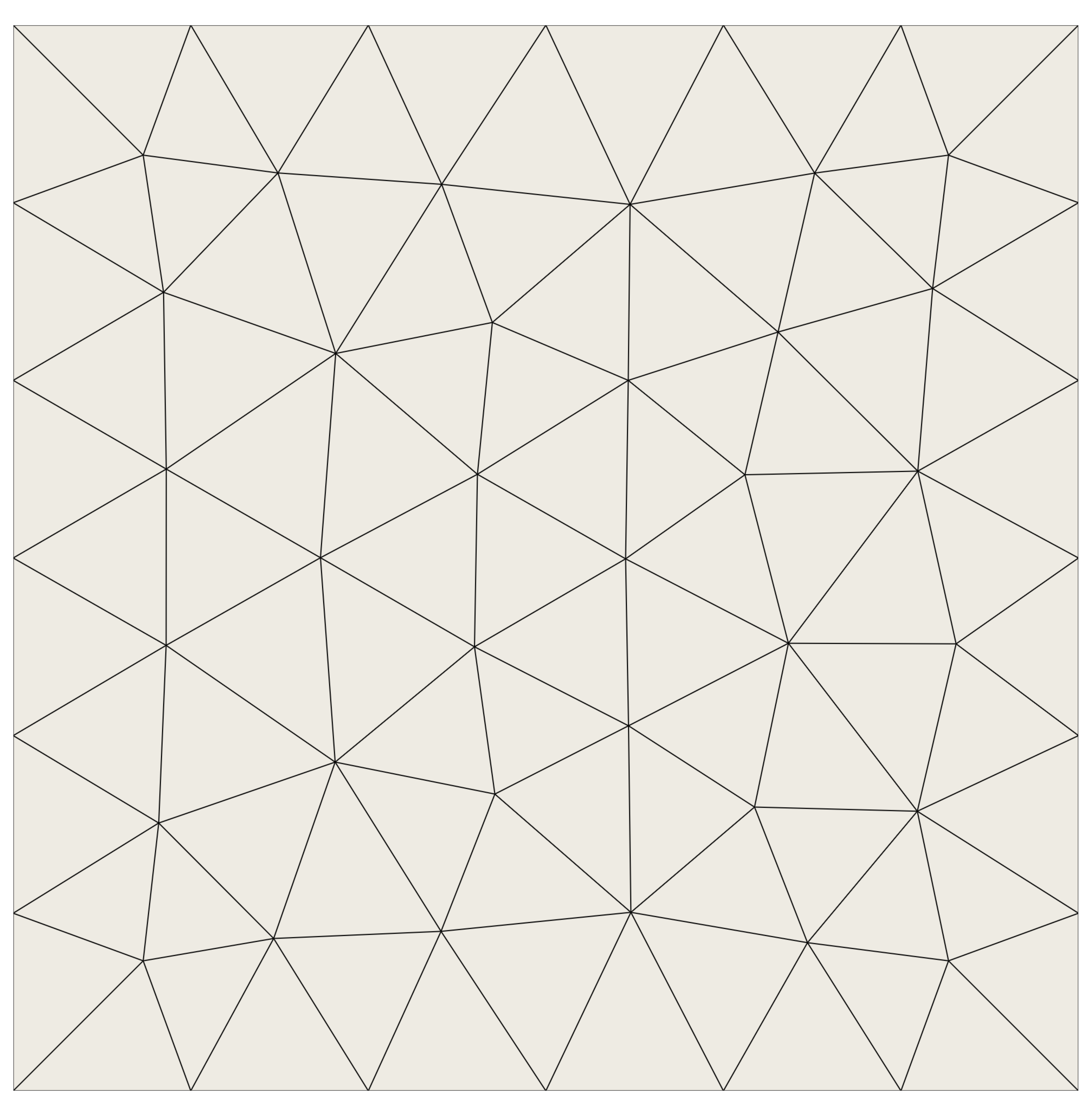}}
    	\end{overpic}
    }
    \subfigure[Regular polygonal mesh]
 	{
 	 	\begin{overpic}[scale=0.8]{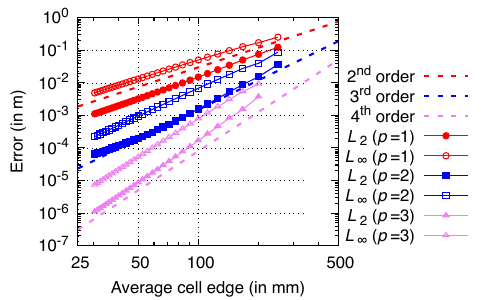}
  			\put(50,13){\includegraphics[scale=0.016]{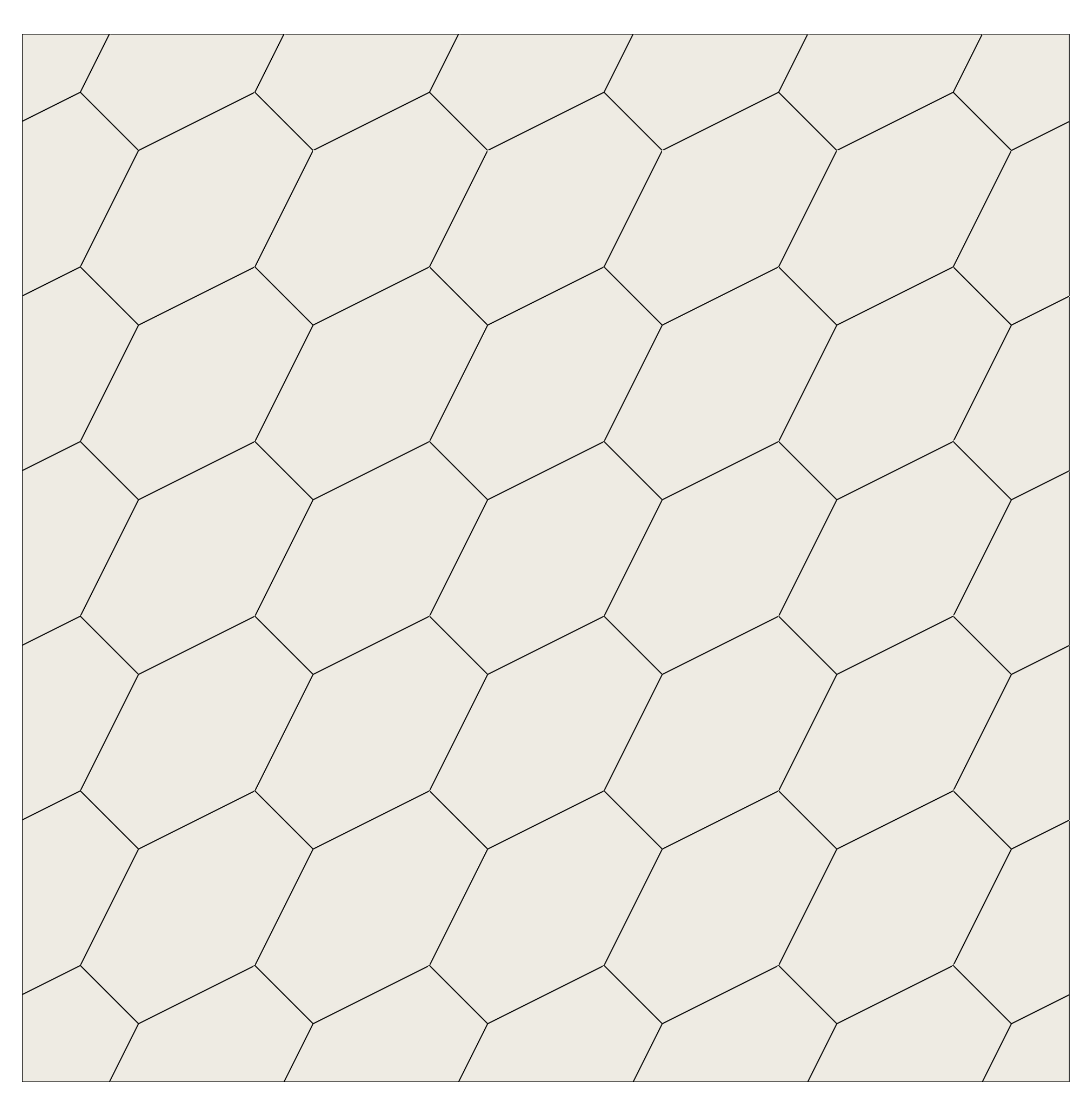}}
    	\end{overpic}
    }
 	\caption{Manufactured solution square (2D case): error convergence for displacement magnitude.}
 	\label{fig:mms2D-disp}
 \end{figure}

\begin{figure}[H]
 	\centering
	\subfigure[Regular quadrilateral mesh]
 	{
 	    \begin{overpic}[scale=0.8]{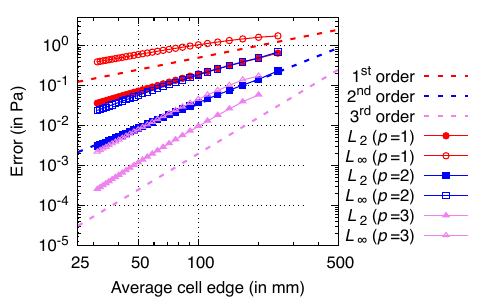}
  			\put(50,13){\includegraphics[scale=0.016]{figures/mms_2D_hex_struct-mesh}}
    	\end{overpic}
    }
 	\subfigure[Regular triangular mesh]
 	{
 	 	\begin{overpic}[scale=0.8]{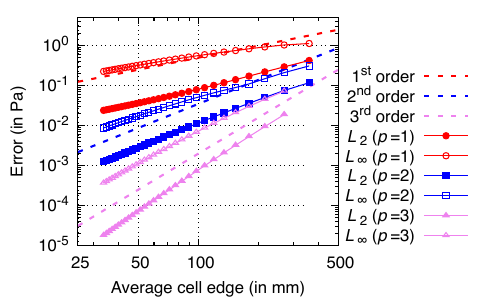}
  			\put(50,13){\includegraphics[scale=0.016]{figures/mms_2D_tet_struct-mesh}}
    	\end{overpic}
    }
    \subfigure[Irregular triangular mesh]
 	{
 	 	\begin{overpic}[scale=0.8]{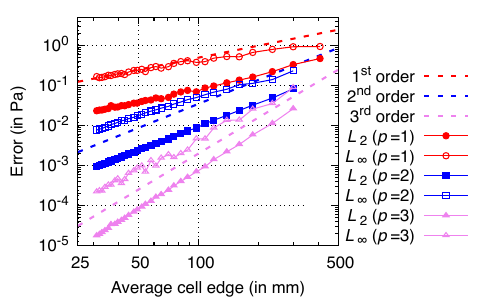}
  			\put(50,13){\includegraphics[scale=0.016]{figures/mms_2D_tet_unstruct_v2-mesh}}
    	\end{overpic}
    }
    \subfigure[Regular polygonal mesh]
 	{
 	 	\begin{overpic}[scale=0.8]{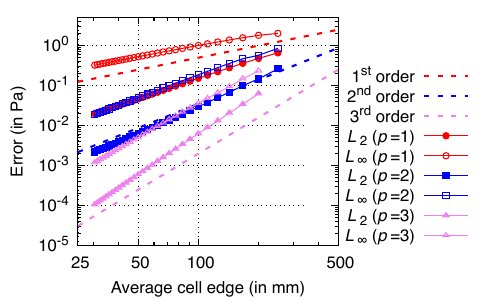}
  			\put(50,13){\includegraphics[scale=0.016]{figures/mms_2D_poly_struct-mesh}}
    	\end{overpic}
    }
 	\caption{Manufactured solution square (2D case): error convergence for stress magnitude.}
 	\label{fig:mms2D-stress}
\end{figure}

Figs.~\ref{fig:mms2D-disp}~and~\ref{fig:mms2D-stress} show the $L_2$ and $L_\infty$ norms of the displacement and stress magnitude errors plotted against the average cell spacing. For visual reference, an image of the corresponding mesh is included in the lower-right corner of each plot. Across all mesh types, the convergence behaviour of both displacement and stress is smooth. The stress errors exhibit the expected orders of accuracy: first-order for $p=1$, second-order for $p=2$, and third-order for $p=3$. Regular hexahedral and polyhedral meshes exhibit slightly higher than expected convergence rates, attributed to error cancellation effects commonly observed on regular grids. 
 Theoretical expectations for displacement errors are reached for second- and fourth-order discretisation, i.e, for $p=1$ and $p=3$. For 
$p=2$, the displacement error exhibits a transition from an initial third-order to second-order convergence as the mesh is refined. Despite the convergence rate degradation, the $p=2$ discretisation consistently yields displacement errors approximately an order of magnitude smaller. Importantly, this reduction in displacement convergence rate does not degrade the stress convergence rate, which remains uniformly second-order. Such behaviour has been reported by several authors \cite{Nishikawa2025a}, although, to the best of the author’s knowledge, a definitive explanation for why odd-order schemes often converge as if they were of the next lower even order has not yet been established.

\begin{figure}[H]
 	\centering
	\subfigure[Hexahedral mesh]
 	{
 	    \begin{overpic}[scale=0.8]{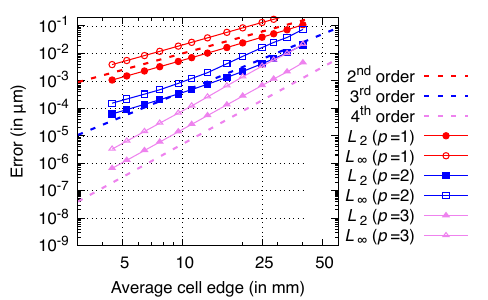}
  			\put(48,13){\includegraphics[scale=0.017]{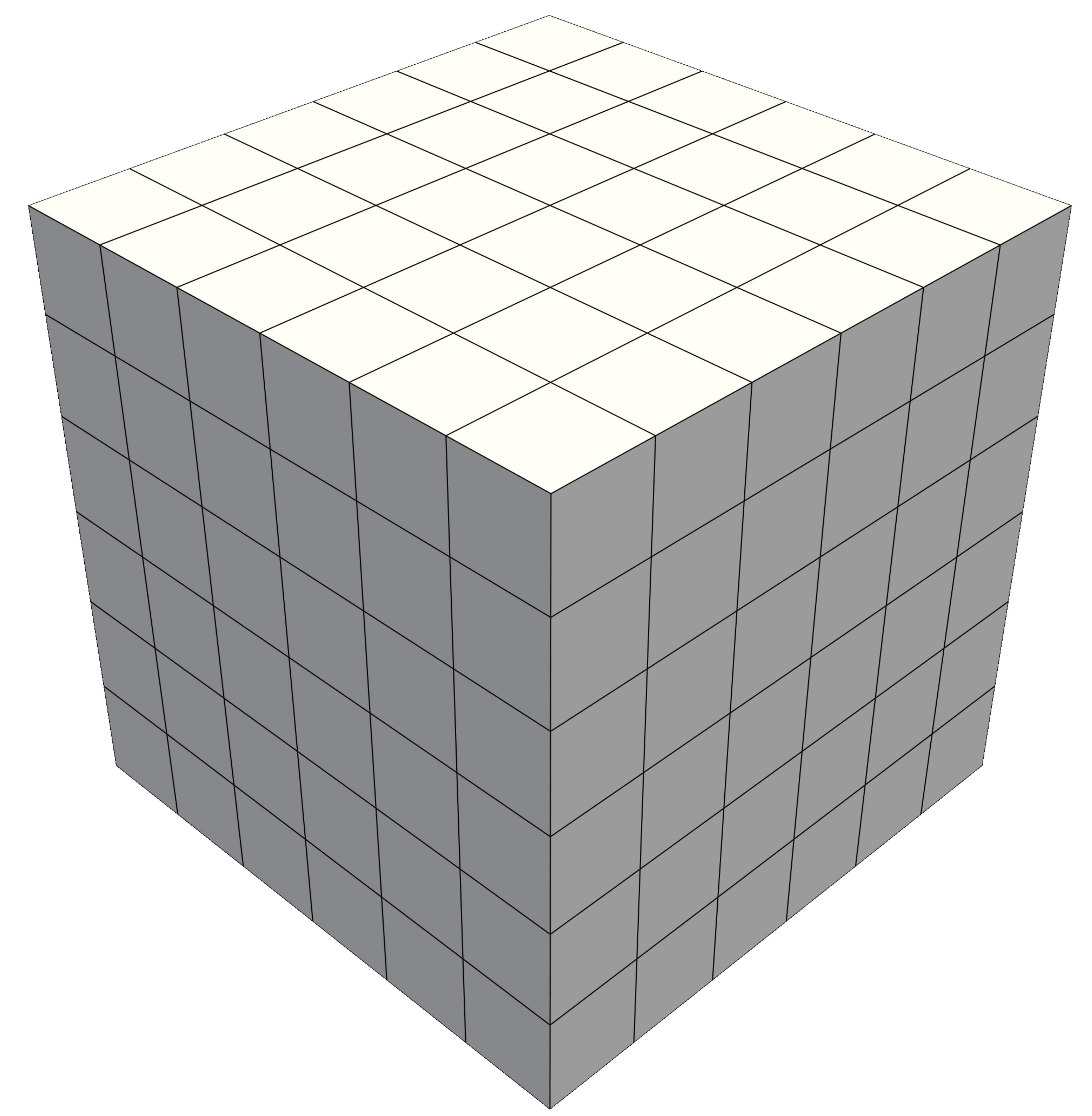}}
    	\end{overpic}
    }
 	\subfigure[Tetrahedral structured mesh]
 	{
 	 	\begin{overpic}[scale=0.8]{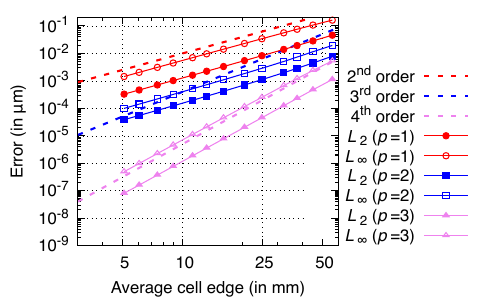}
  			\put(48,13){\includegraphics[scale=0.017]{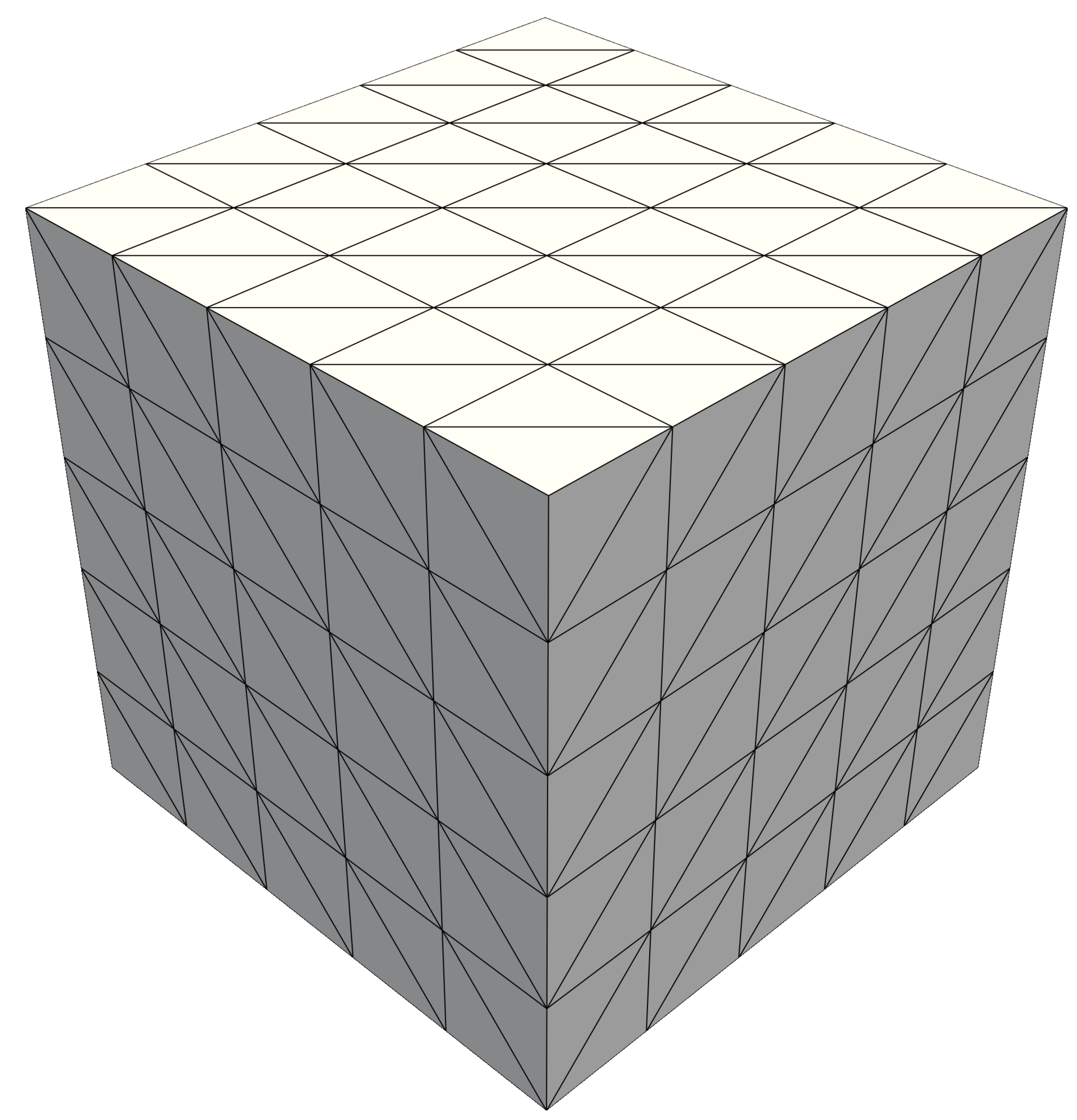}}
    	\end{overpic}
    }
    \subfigure[Tetrahedral unstructured mesh]
 	{
 	 	\begin{overpic}[scale=0.8]{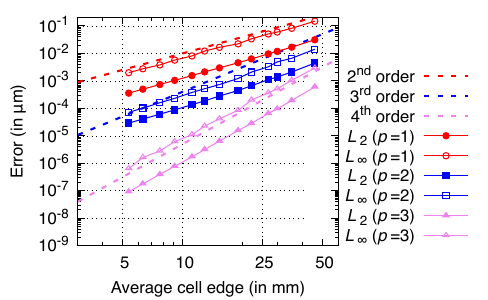}
  			\put(48,13){\includegraphics[scale=0.017]{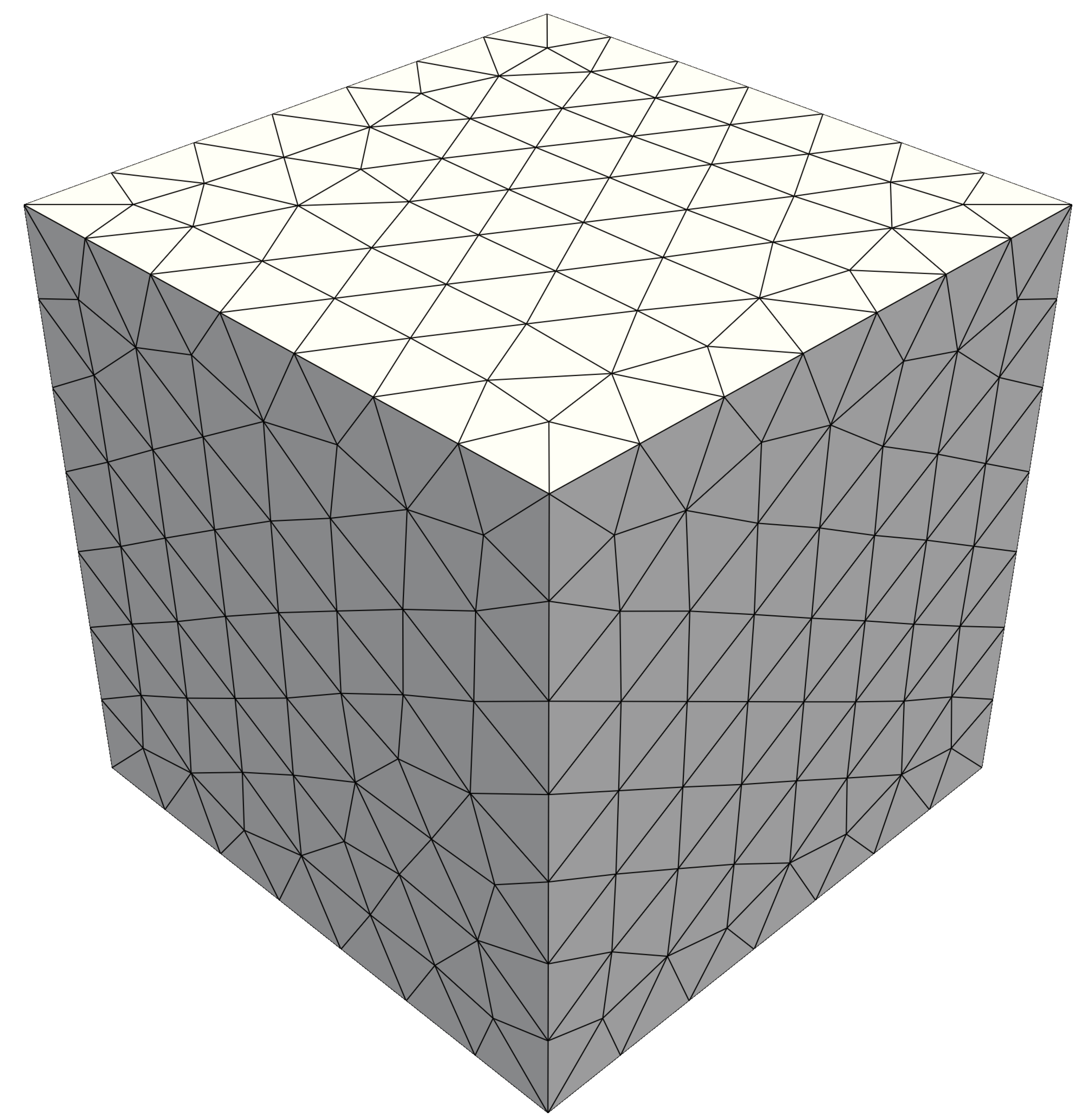}}
    	\end{overpic}
    }
    \subfigure[Polyhedral structured mesh]
 	{
 	 	\begin{overpic}[scale=0.8]{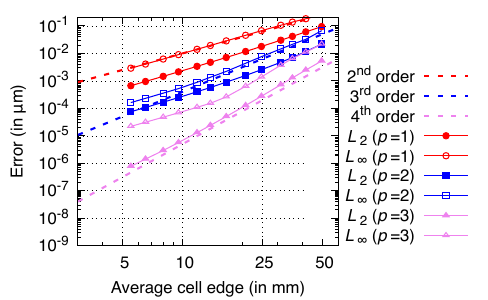}
  			\put(48,13){\includegraphics[scale=0.017]{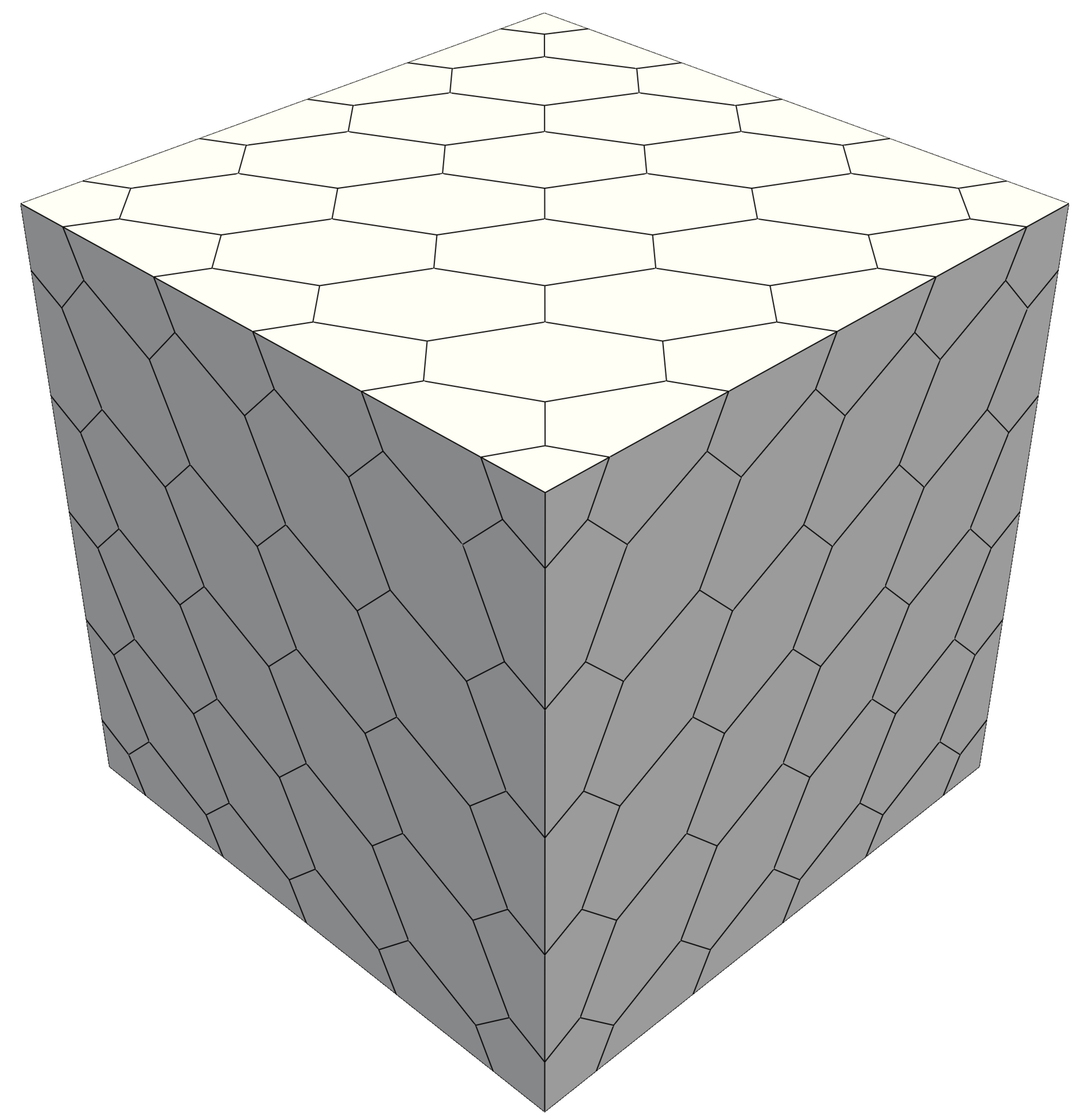}}
    	\end{overpic}
    }
 	\caption{Manufactured solution cube (3D case): error convergence for displacement magnitude.}
 	\label{fig:mms3D-disp}
 \end{figure}

\begin{figure}[H]
 	\centering
	\subfigure[Hexahedral mesh]
 	{
 	    \begin{overpic}[scale=0.8]{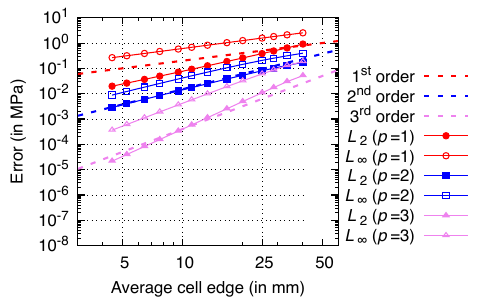}
  			\put(48,13){\includegraphics[scale=0.017]{figures/mms_3D_hex_struct-mesh}}
    	\end{overpic}
    }
 	\subfigure[Tetrahedral structured mesh]
 	{
 	 	\begin{overpic}[scale=0.8]{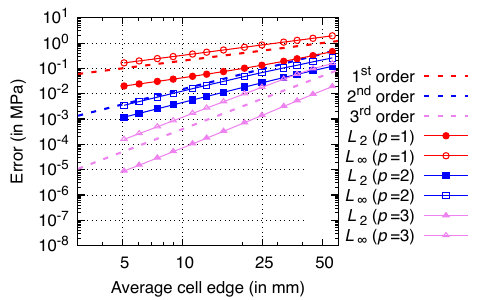}
  			\put(48,13){\includegraphics[scale=0.017]{figures/mms_3D_tet_struct-mesh}}
    	\end{overpic}
    }
    \subfigure[Tetrahedral unstructured mesh]
 	{
 	 	\begin{overpic}[scale=0.8]{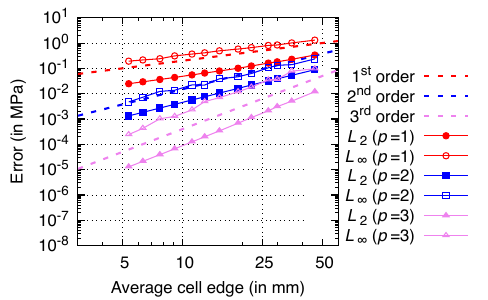}
  			\put(48,13){\includegraphics[scale=0.017]{figures/mms_3D_tet_unstruct_v2-mesh}}
    	\end{overpic}
    }
    \subfigure[Polyhedral structured mesh]
 	{
 	 	\begin{overpic}[scale=0.8]{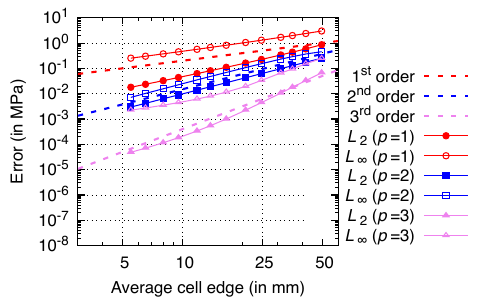}
  			\put(48,13){\includegraphics[scale=0.017]{figures/mms_3D_poly_struct-mesh}}
    	\end{overpic}
    }
 	\caption{Manufactured solution cube (3D case): error convergence for stress magnitude.}
 	\label{fig:mms3D-stress}
\end{figure}

Figs.~\ref{fig:mms3D-disp}~and~\ref{fig:mms3D-stress} present the corresponding 3D results: convergence of the $L_2$ and $L_\infty$ norms of the displacement and stress magnitude errors.  The trends closely mirror those observed in the 2D case: smooth convergence, expected order of accuracy for the stress field, and the same characteristic reduction of displacement convergence from third- to second-order for $p=2$. These results confirm that the observed behaviour is not specific to mesh dimensionality or topology, but rather intrinsic to the reconstruction order.
 % line width 2, Ambient 0.15. mesh number 6 is shown at diagrams
 %
 %
 %
%
%
%
%
%------------------------------------------------------------------------------
\subsubsection{Case 2: Cantilever Beam}
\begin{figure}[H]
	\centering
	\subfigure[Case geometry (in m) and boundary conditions]
 	{
 		\label{fig:cantilever-geometry}
    	\includegraphics[scale=1]{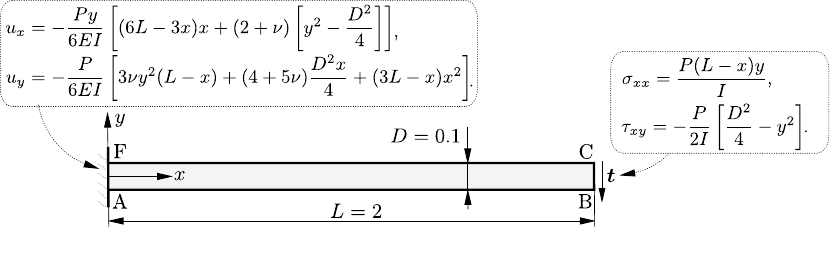}
    }
    \subfigure[Tetrahedral unstructured mesh]
 	{
 		\label{fig:cantilever-tetstruct-mesh}
        \includegraphics[scale=0.1]{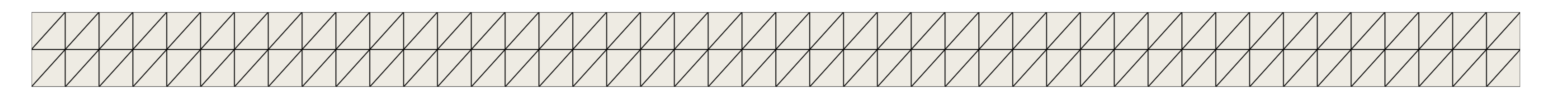}
    }
    \subfigure[Tetrahedral unstructured mesh]
 	{
 		\label{fig:cantilever-tetunstruct-mesh}
        \includegraphics[scale=0.1]{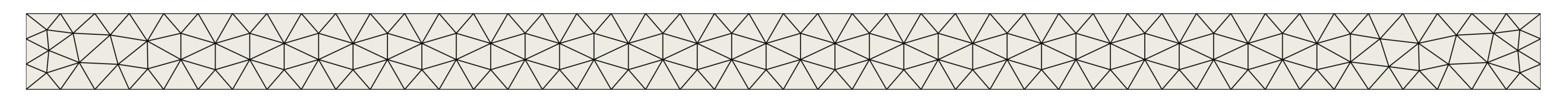}
    }  
	\label{fig:cantilever}
	\caption{Cantilever beam case: case geometry (in m), boundary conditions and the coarsest triangular meshes used.}
\end{figure}

%------------------------------------------------------------------------------
%
The test case geometry, shown in Fig.~\ref{fig:cantilever-geometry}, is a rectangular 2D beam with dimensions of $L\times D=2\,\text{m}\times0.1\,\text{m}$ (unit depth), a Young’s modulus of $E = 200$ GPa, and a Poisson’s ratio of $\nu = 0.3$. Beam edges AB and CF are traction-free, and plane strain conditions are assumed. The edge AF is constrained by the analytical displacement, while the edge BC is subjected to the load $P=1 \cdot 10^{5}$ N. This setup enables quantification of the difference between the predicted displacement and the analytical solution as a measure of convergence across the entire domain, rather than only at the beam end \cite{augarde2008}. Convergence is assessed using 15 successively refined regular and irregular triangular meshes, shown at Fig.~\ref{fig:cantilever-tetstruct-mesh} and \ref{fig:cantilever-tetunstruct-mesh}, respectively. The finest irregular triangular mesh consists of $12\;290$ cells while the finest regular triangular mesh of $10\;368$ cells.

\begin{figure}[H]
 	\centering
 	\subfigure[Tetrahedral structured mesh]
 	{
 		\label{fig:cantilever-dispErrors-tetstruct}
    	\includegraphics[scale=0.8]{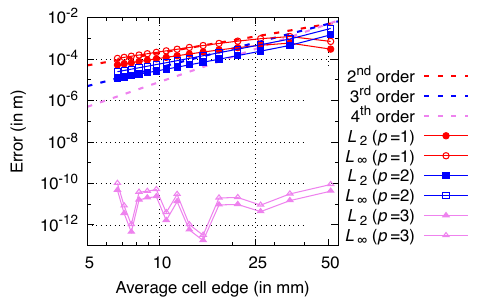}
    }
    \subfigure[Tetrahedral unstructured mesh]
 	{
 		\label{fig:cantilever-dispErrors-tetunstruct}
        \includegraphics[scale=0.8]{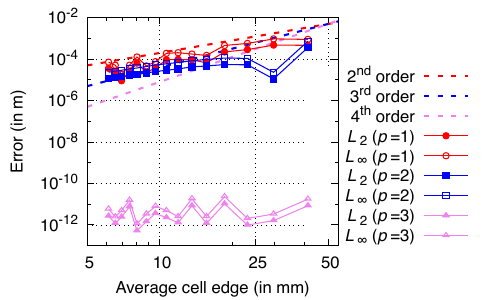}
    }
 	\caption{Cantilever beam case: displacement magnitude discretisation errors.}
 	\label{fig:cantilever-dispErrors}
 \end{figure}

The analytical displacement field is a third-order polynomial, implying that a fourth-order discretisation should reproduce the solution exactly. In this sense, the problem serves as a patch test for the fourth-order scheme. Fig.~\ref{fig:cantilever-dispErrors} presents the $L_2$ and $L_\infty$ convergence of the displacement error. As expected, the fourth-order discretisation ($p=3$) achieves the exact displacement field over the entire domain, reaching machine precision, across all meshes. The second and third-order schemes both converge at a second-order rate.

\begin{figure}[H]
 	\centering
 	\subfigure[Tetrahedral structured mesh]
 	{
 		\label{fig:cantilever-stressErrors-tetstruct_a}
    	\includegraphics[scale=0.8]{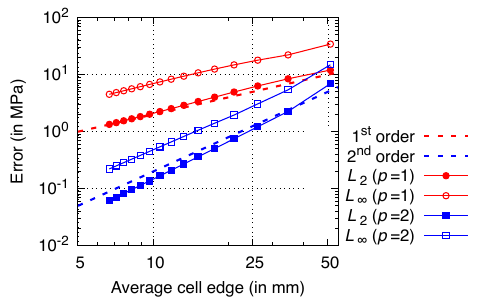}
    }
    \subfigure[Tetrahedral unstructured mesh]
 	{
 		\label{fig:cantilever-stressErrors-tetstruct_b}
    	\includegraphics[scale=0.8]{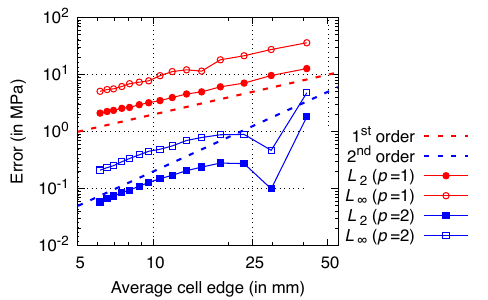}
    }
 	\caption{Cantilever beam case: error convergence for stress magnitude. Results for $p=3$ are omitted for brevity, i.e. stress error for $p=3$ is below $\sim$ $10^{-6}$ MPa.}
 	\label{fig:cantilever-stressErrors}
 \end{figure}

Fig.~\ref{fig:cantilever-stressErrors} shows the corresponding convergence of the stress error norms. As expected, the $p=1$ and $p=2$ discretisation converge with first- and second-order accuracy, respectively.
 %
%
%------------------------------------------------------------------------------
\subsubsection{Case 3: Hole in a Infinite Plate Subjected to Remote Stress}
\label{case:plateHole}
\begin{figure}[H]
 	\centering
 	\subfigure[Case geometry (in m)]
 	{
 		\label{fig:plateHole:a}
    		\includegraphics[scale=0.95]{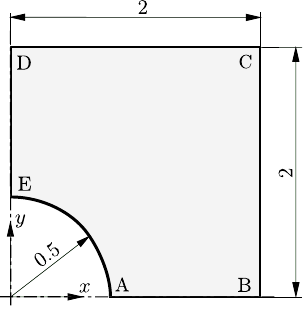}
    }
 	\subfigure[Tetrahedral unstructured mesh]
 	{
 		\label{fig:plateHole:b}
    	\includegraphics[scale=0.051]{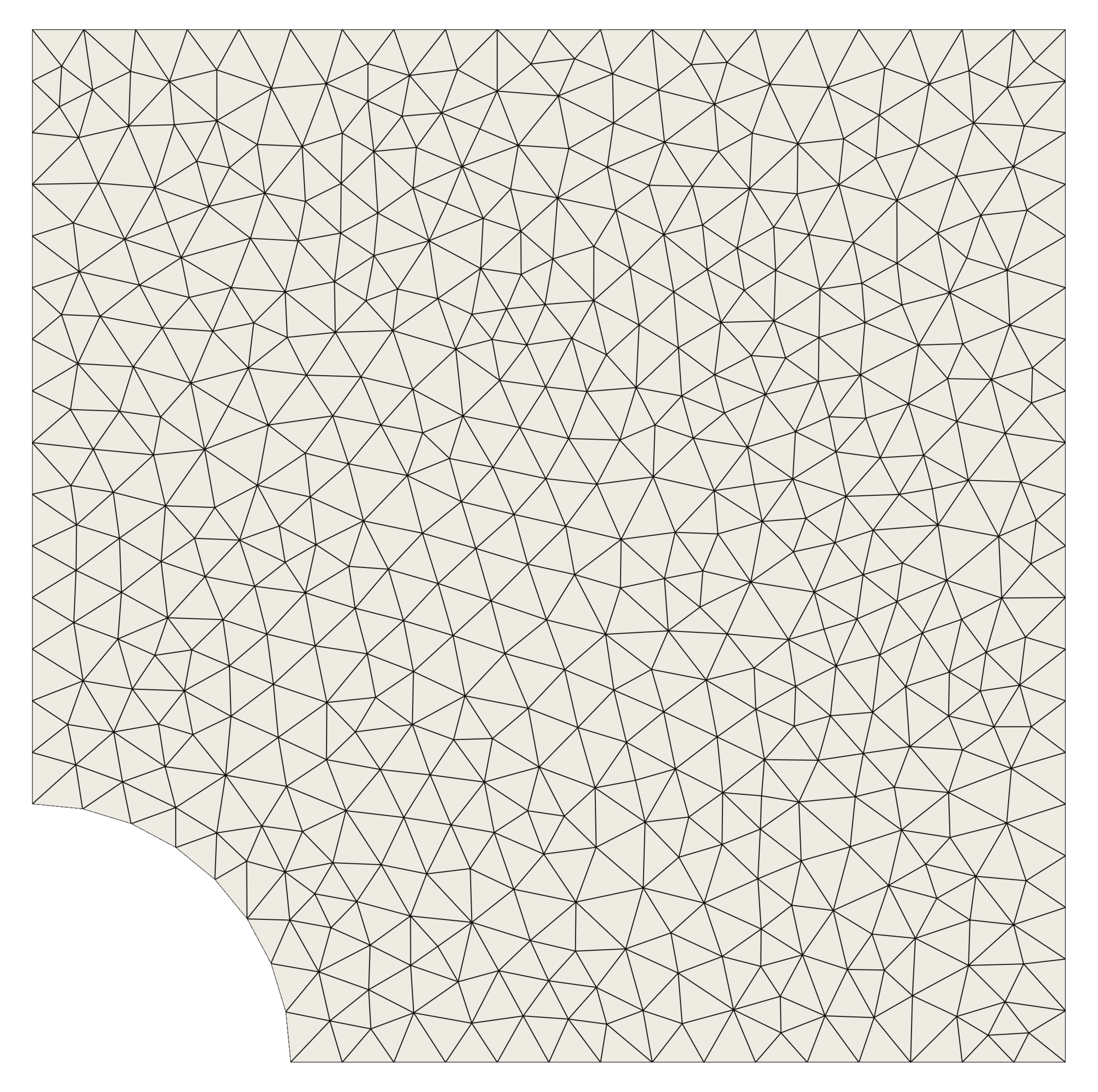}
    }
    \subfigure[Von Mises stress $\sigma_\text{eq}$  distribution (finest mesh, $p=3$)]
 	{
 		\label{fig:plateHole:c}
    	\includegraphics[scale=0.79]{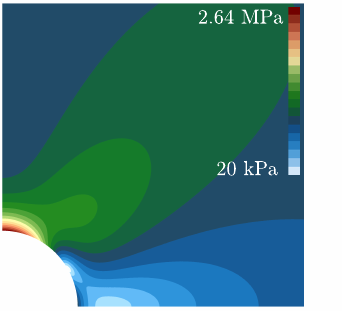}
    }
 	\caption{Plate hole case: geometry, mesh and resulting stress distribution.}
 	\label{fig:plateHole}
\end{figure}

%------------------------------------------------------------------------------
%
% Case can be also found in:
% Finite volume method for stress analysis in complex domain
% Application of the nite volume method and unstructured meshes to linear elasticity
%
This benchmark problem consists of a thin, infinitely large plate with a circular hole, subjected to uniaxial tension of  $\sigma_{xx} = T = 1$ MPa, see Fig.~\ref{fig:plateHole:a}. Owing to the symmetry of the geometry and loading, only one quarter of the plate is modelled. To minimise the influence of the finite computational boundaries, the exact tractions obtained from the analytical solution \cite{Demirdzic1997} are prescribed on the outer edges BC and CD. Symmetry boundary conditions are applied on boundaries AB and DE, while zero traction is specified on the hole boundary. The material properties are defined by a Young’s modulus of $E = 200 \,\text{GPa}$ and a Poisson’s ratio of $\nu = 0.3$.
Five successively refined irregular triangular meshes are used: 978 (Fig.~\ref{fig:plateHole:b}), $4\;050$, $16\;296$, $64\;257$ and $256\;512$ cells.

\begin{figure}[H]
 	\centering
	\subfigure[Displacement magnitude discretisation errors]
 	{
 		\label{fig:plateHole-disp}
    	\includegraphics[scale=0.8]{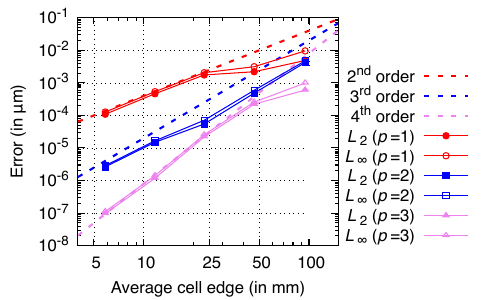}
    }
 	\subfigure[Stress magnitude discretisation errors]
 	{
 		\label{fig:plateHole-stress}
    	\includegraphics[scale=0.8]{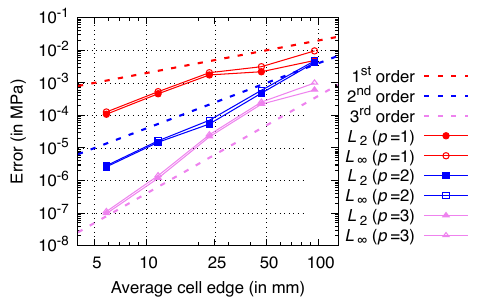}
    }
 	\caption{Plate hole case: displacement and stress error convergence.}
 	\label{fig:plateHole-convergence}
 \end{figure}

Fig.~\ref{fig:plateHole:c} shows the predicted equivalent stress distribution, highlighting the steep stress gradients and stress concentration around the hole. The convergence of the $L_2$ and $L_\infty$ norms of the displacement and stress discretisation errors is presented in Fig.~\ref{fig:plateHole-convergence}. As before, the expected convergence rates for both displacement and stress are achieved. Notably, the presence of the curved hole boundary does not impair the convergence behaviour, despite the absence of any special treatment for this curved edge.
%.
%
%
%
%------------------------------------------------------------------------------
\subsubsection{Case 4: Internally Pressurised Thick-walled Cylinder}
\begin{figure}[H]
 	\centering
 	\subfigure[Case geometry]
 	{
 		\label{fig:pressurisedCylinder:a}
    		\includegraphics[scale=1]{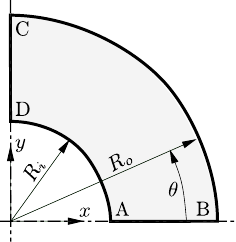}
    }
	\subfigure[Irregular triangular mesh (M1 mesh)]
 	{
 		\label{fig:pressurisedCylinder:b}
    		\includegraphics[scale=0.052]{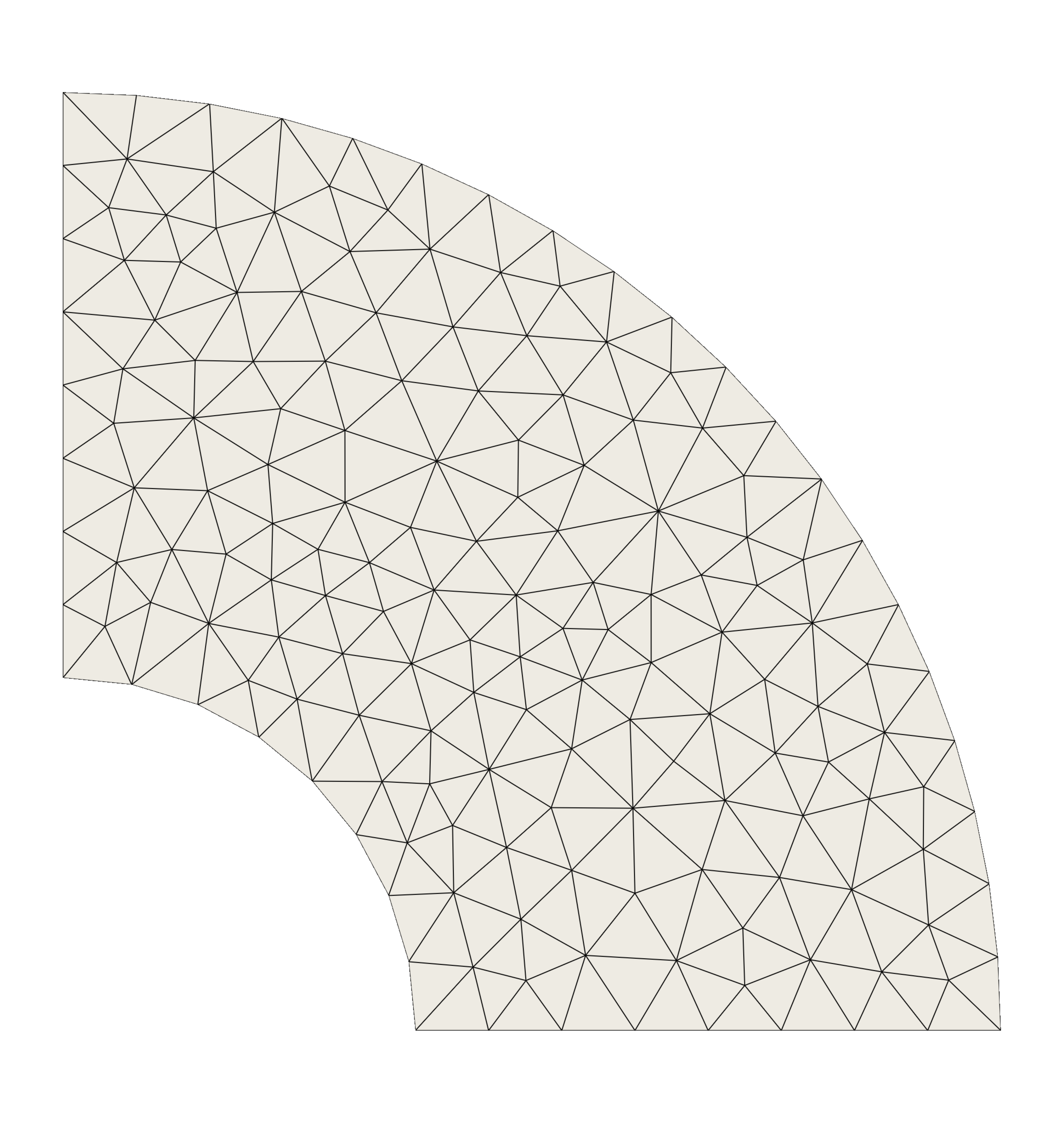}
    }\quad
    \subfigure[Reference configuration (red line) and deformed configuration, coloured by the displacement field (in m, M3 mesh, $p=3$)]
 	{
 		\label{fig:pressurisedCylinder:c}
    		\includegraphics[scale=0.74]{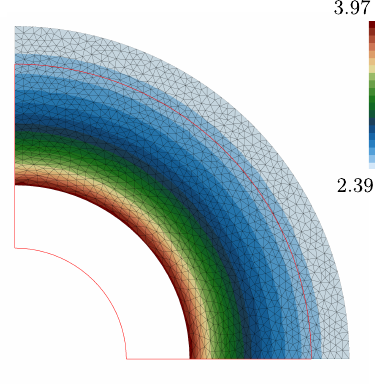}
    }
 	\caption{Pressurised cylinder case: geometry and the coarsest mesh.}
 	\label{fig:pressurisedCylinder}
\end{figure}

%------------------------------------------------------------------------------
%
In this case, a homogeneous thick-wall cylindrical pressure vessel with an inner radius $R_i = 7$ m, outer radius $R_o = 18.625$ m, and loaded internally with pressure $p = 100$ MPa is analysed. Two types of material are considered:
\begin{itemize}
\item[i.] Small strain, linear-elastic \cite{Bijelonja2006}: $E=10$ GPa and $\nu=0.3$.
\item[ii.] Finite strain, Mooney-Rivlin law \cite{Bijelonja2005a}: $c_{10} = 80$ MPa, $c_{01} = 20$ MPa, $c_{11} = 0.0$ MPa and penalty parameter (bulk modulus) $\kappa=400$ MPa, see Appendix \ref{app:mechLaws:MR}.
\end{itemize}
The problem is considered under plane-stress assumptions, with the 2D computational domain corresponding to a quarter of the cylindrical geometry, as shown in Fig.~\ref{fig:pressurisedCylinder:a}. Symmetry boundary conditions are applied on edges AB and CD. Gravitational and inertial effects are neglected. The linear-elastic case is solved using a single loading increment, whereas the hyperelastic case is solved using ten equal loading increments. Analytical solutions are available in \cite{Timoshenko1970} for the small-strain and in \cite{Green1992} for the finite-strain formulation.

\begin{figure}[H]
 	\centering
	\subfigure[Displacement convergence]
 	{
 		\label{fig:pressurisedCylindere-error:a}
    	\includegraphics[scale=0.8]{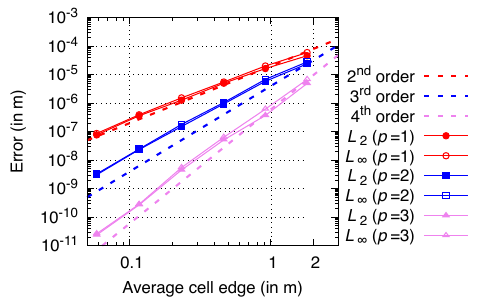}
    }
 	\subfigure[Stress convergence]
 	{
 		\label{fig:pressurisedCylinder-error:b}
    	\includegraphics[scale=0.8]{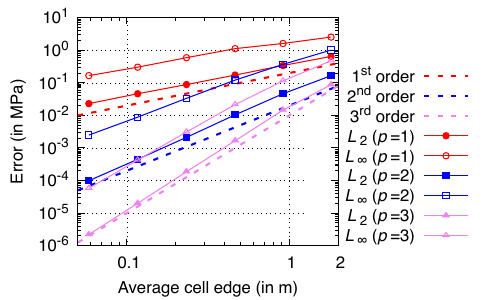}
    }
 	\caption{Pressurised cylinder (small strain case): displacement and stress magnitude discretisation errors for linear elastic material.}
 	\label{fig:pressurisedCylinder-error}
 \end{figure}

\begin{figure}[H]
 	\centering
 	\subfigure[Radial displacement]
 	{
 	 	\label{fig:pressurisedCylinder:displacement}
    	\includegraphics[scale=0.8]{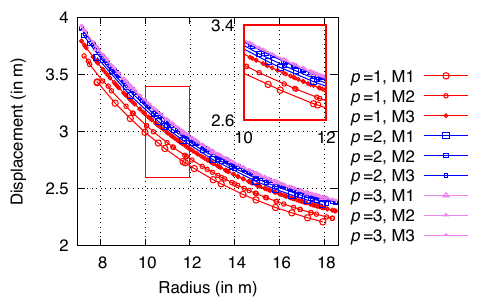}
    }
 	\subfigure[Radial and hoop stress profiles (M3 mesh)]
 	{
 	 	\label{fig:pressurisedCylinder:stress}	
    	\includegraphics[scale=0.8]{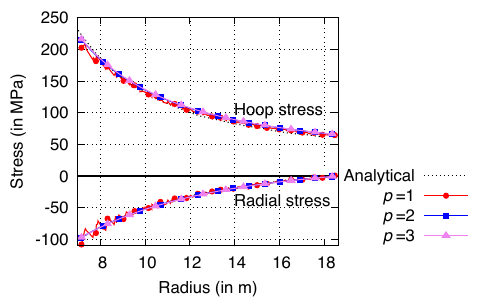}
    }
 	\caption{Pressurised cylinder (finite strain case): radial displacement and stress profiles.}
 	\label{fig:pressurisedCylinder:finiteStrain}
\end{figure}

The linear-elastic case is solved using a series of successively refined irregular triangular meshes. Results are presented in terms of the $L_2$ and $L_\infty$ error norms for displacement and stress, as shown in Fig.~\ref{fig:pressurisedCylinder-error}. Interestingly, the absence of a specialised boundary curvature treatment does not affect the order of convergence, and the expected rates are achieved. This behaviour can be explained by the fact that the problem effectively reduces to a one-dimensional \cite{Bijelonja2005a}, and by the uniform stress state along both the inner and outer boundaries.
The finite-strain case employs a Mooney–Rivlin material model with penalty parameter that yields a Poisson ratio of $\nu=0.4975$, close to the incompressibility limit. In general, accurately approaching the incompressible limit requires a mixed pressure–displacement formulation \cite{Bijelonja2005a, Horvat2025}; alternatively, in second-order discretisations, pressure oscillations may be mitigated using pressure-smoothing techniques \cite{oliveira2020}. However, no such stabilisation or smoothing procedure is applied here. The used Poisson ratio is sufficiently close to incompressibility, and the analytical  solution for incompressible material can be used to validate the resulting stress field. The case is solved on three meshes: M1 with 292 cells (Fig.~\ref{fig:pressurisedCylinder:b}), M2 with 645 cells, and M3 with 2507 cells. Mesh M3 in the deformed configuration, coloured by the $\sigma_\text{eq}$ field, is shown in Fig.~\ref{fig:pressurisedCylinder:c}. Fig.~\ref{fig:pressurisedCylinder:displacement} presents the radial displacement profiles for all three meshes and interpolation orders. It is evident that $p=3$ produces practically identical displacement profiles on all meshes, while $p=1$ converges toward the same solution more slowly than $p=2$, as expected. Fig.~\ref{fig:pressurisedCylinder:stress} shows the radial and hoop stress profiles, which agree closely with the analytical solution. Oscillations are observed in the stress field for $p=1$, which is expected for linear interpolation in a displacement-based formulation when approaching the incompressibility limit \cite{Horvat2025, oliveira2020}. The profiles in Fig.~\ref{fig:pressurisedCylinder:finiteStrain} are obtained by sampling a band of cells near the line $\theta = 45^\circ$ and plotting their corresponding values.
%
%
%------------------------------------------------------------------------------
\subsubsection{Case 5: Cook's Membrane}
\begin{figure}[H]
 	\centering
 	\subfigure[Case geometry and dimensions (in mm). Taken from \cite{Cardiff2025}.]
 	{
 		\label{fig:CooksMembrane:a}
    	\includegraphics[scale=1]{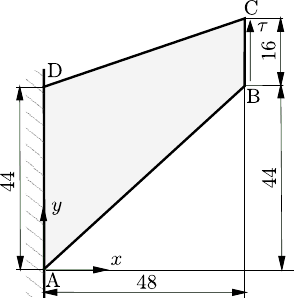}
    }
    \quad
	\subfigure[Irregular triangular mesh (second coarsest, 35 cells)]
 	{
 		\label{fig:CooksMembrane:b}
    	\includegraphics[scale=0.055]{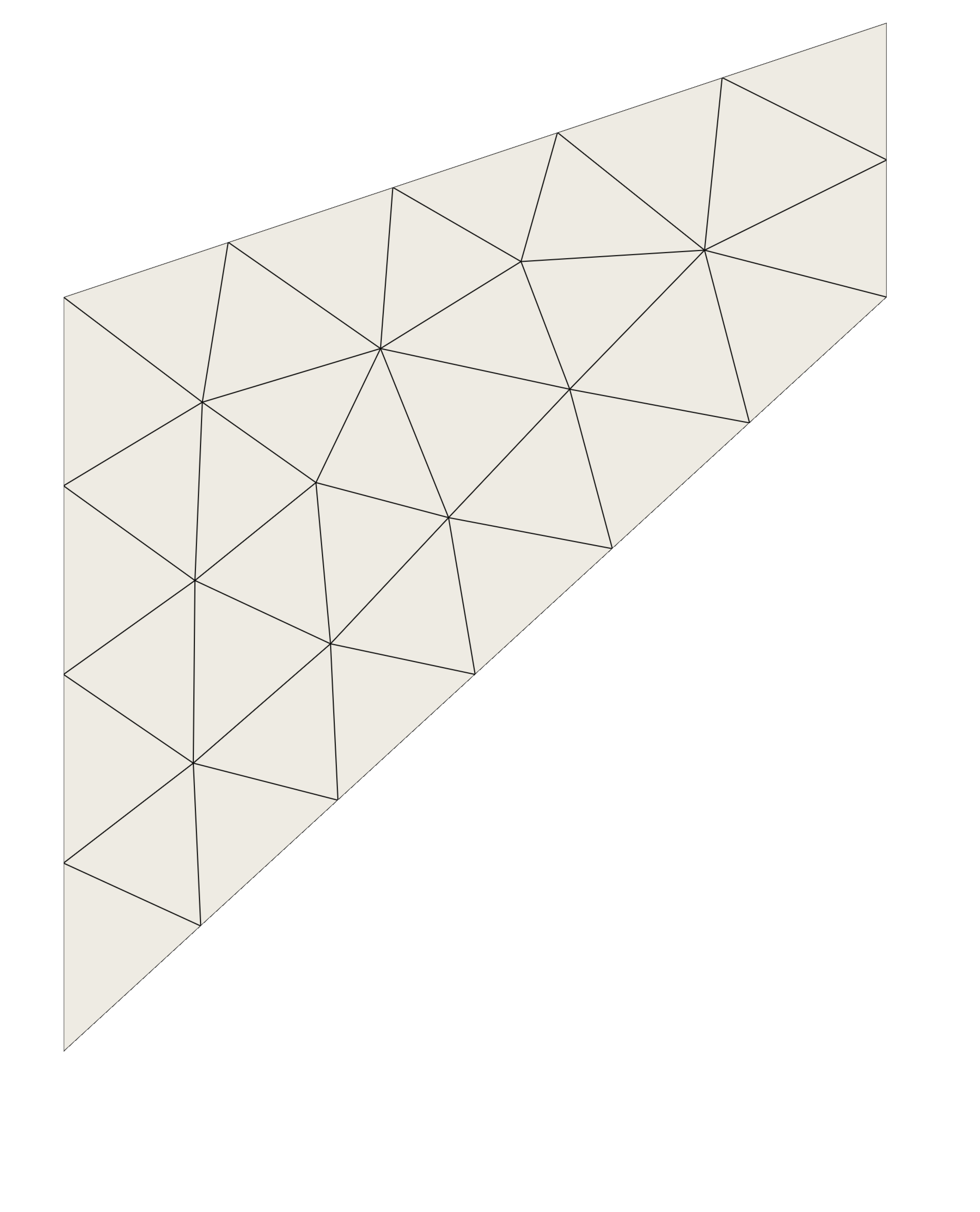}
    }
    \quad
    \subfigure[Reference configuration (red line) and deformed configuration for the hyperelastic case, coloured by the $\sigma_\text{eq}$ field (finest mesh, $p=3$)]
 	{
 		\label{fig:CooksMembrane:c}
    	\includegraphics[scale=0.34]{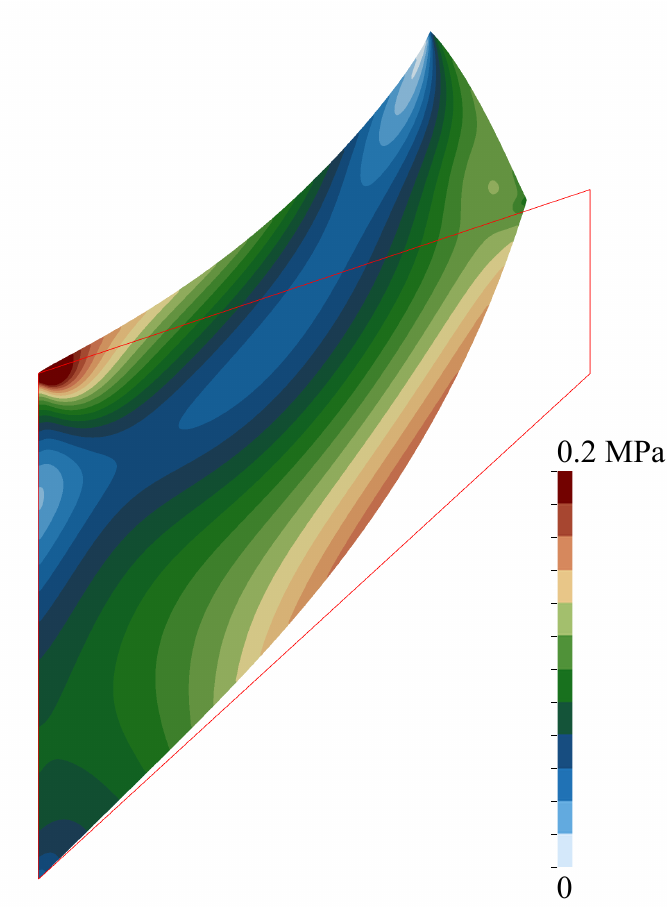}
    }
 	\caption{Cook's membrane case: geometry, mesh and deformation in the hyperelastic case.}
 	\label{fig:CooksMembrane}
\end{figure}

%------------------------------------------------------------------------------
%
The Cook’s membrane configuration (Fig.~\ref{fig:CooksMembrane:a}) is a classic benchmark problem dominated by bending effects.
It consists of a two-dimensional trapezoidal panel (plane-strain condition) that is fully constrained along edge AD and loaded by a uniform shear traction applied to BC edge. Edges AB and CD are traction-free.
%The coordinates of the four vertices of the trapezoid, given in millimetres, are $(0,0)$, $(48,44)$, $(48,60)$, and $(0,44)$.
In this study, two variants of the problem are analysed:
\vspace{0.5em}
\begin{enumerate}[label=\roman*.]
\item Small-strain, linear elastic \cite{Zienkiewicz2000,Simplas}: $E = 70$ MPa, $\nu = 1/3$, and $\tau = 6.25$ kPa.
\item Finite-strain, neo-Hookean hyperelastic \cite{Pelteret2018}: $E = 1.0985$ MPa, $\nu = 0.3$, and $\tau = 62.5$ kPa, see Appendix \ref{app:mechLaws:NH}.
\end{enumerate}
\vspace{0.5em}
\noindent
where $E$ denotes the Young’s modulus, $\nu$ the Poisson ratio, and $\tau$ is the imposed shear traction.
The analysis is performed as quasi-static without body forces.
The linear case is solved in a single step, whereas the hyperelastic case is solved using 30 uniform load increments.
Eight successively refined irregular tetrahedral meshes are used 11, 35 (Fig.~\ref{fig:CooksMembrane:b}), 164, 177, 263, 389, 1~457, 5~637, 22~119 and 88~012 cells.

\begin{figure}[H]
 	\centering
	\subfigure[Linear elastic case]
 	{
 		\label{fig:CooksMembrane-tip:a}
    	\includegraphics[scale=0.8]{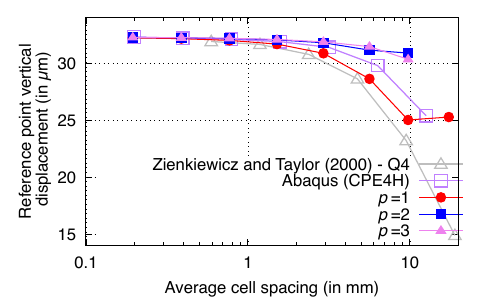}
    }
 	\subfigure[Hyperelastic case]
 	{
 		\label{fig:CooksMembrane-tip:b}
    	\includegraphics[scale=0.8]{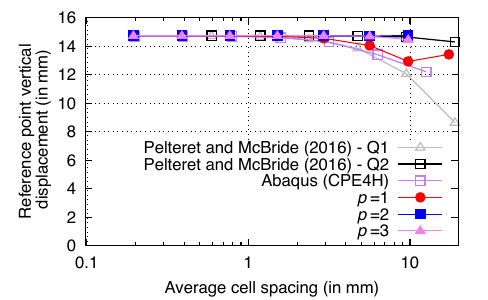}
    }
 	\caption{Cook’s membrane vertical displacement predictions at the reference point as a function of the average cell spacing. Comparisons are given with the results from the literature  \cite{Zienkiewicz2000, Simplas, Pelteret2018}. Case with $p=1$ uses $n+=7$ instead of $n+=10$ to get results on the first mesh with 11 cells.}
 	\label{fig:CooksMembrane-tipConvergence}
 \end{figure}

Fig.~\ref{fig:CooksMembrane:c} shows the predicted equivalent stress distribution on the deformed geometry (using the mesh with 88012 cells) for the hyperelastic cases.  To verify the results, Fig.~\ref{fig:CooksMembrane-tipConvergence} plots the vertical displacement of a key reference point against the average cell spacing, comparing the predictions to literature data and Abaqus (C3D8 element) results. 
The location of this reference point differs by analysis type: it is the top-right corner (48, 60) mm for the elastic case but the midpoint of the loaded edge (48, 52) mm for the hyperelastic case.  
%
%
%------------------------------------------------------------------------------
\subsubsection{Case 6: Spherical Cavity in an Infinite Solid Subjected to Remote Stress}
\label{case:sphericalCavity}
\begin{figure}[H]
 	\centering
 	\subfigure[Unstructured tetrahedral mesh (in m)]
 	{
 		\label{fig:sphericalCavity:a}
    	\includegraphics[scale=1]{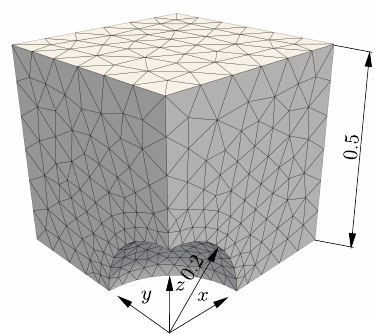}
    }
	\subfigure[Von Mises stress $\sigma_\text{eq}$ distribution (finest mesh, $p=3$)]
 	{
 		\label{fig:sphericalCavity:b}
    	\includegraphics[scale=1]{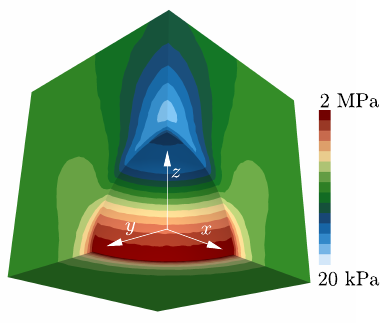}
    }
 	\caption{Spherical cavity case: geometry, mesh and resulting von Mises stress distribution.}
 	\label{fig:sphericalCavity}
\end{figure}

%------------------------------------------------------------------------------
%
This problem considers a spherical cavity embedded within an infinite and isotropic medium with linear elastic material properties ($E = 200$ GPa, $\nu = 0.3$).
The solid is subjected to a far-field uniaxial tensile stress of $\sigma_{zz} = T = 1$ MPa, with all other stress components being zero.
A key feature of this problem is the pronounced stress concentration that forms near the cavity on the plane perpendicular to the applied load; 
this stress rapidly diminishes further from the cavity. 
Although the problem is inherently axisymmetric, it is solved here as a 3D case with a mesh graded towards the cavity.
The computational domain is taken as one-eighth of a $1\,\text{m}~\times~1\,\text{m}~\times~1\,\text{m}$ cube aligned with the Cartesian axes, with one corner at the centre of the sphere, see Fig.~\ref{fig:sphericalCavity:a}.
The analytical expressions for stress \cite{Southwell1926} and displacement \cite{Goodier1933} can be found in \cite{Cardiff2025}.
The analytical tractions are applied to the far-field boundaries to mitigate finite geometry effects. 
The same analytical expressions are also used on the cavity surface, which theoretically corresponds to a zero-traction boundary; however, the quadrature points on the discretised faces do not lie exactly on the true zero-traction surface. Imposing a zero-traction condition directly on this curved boundary, without specialised boundary treatment, would therefore reduce the overall error convergence rate to second order, as it would be dominated by the second-order accurate geometrical approximation \cite{costa2021, Gooch2002}.  
Curved-boundary treatment is beyond the scope of the present work and is left for future investigation.
The analysis is performed using six different irregular tetrahedral meshes with cell counts of $3\;802$, $7\;030$, $14\;891$, $32\;864$, $72\;945$ and $158\;870$, with the coarsest mesh shown in Fig.~\ref{fig:sphericalCavity:a}.

\begin{figure}[H]
 	\centering
 	\subfigure[Displacement convergence]
 	{
 		\label{fig:sphericalCavity:disp}
    	\includegraphics[scale=0.8]{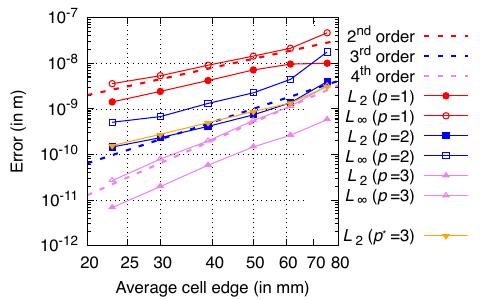}
    }
	\subfigure[Stress convergence]
 	{
 		\label{fig:sphericalCavity:stress}
    	\includegraphics[scale=0.8]{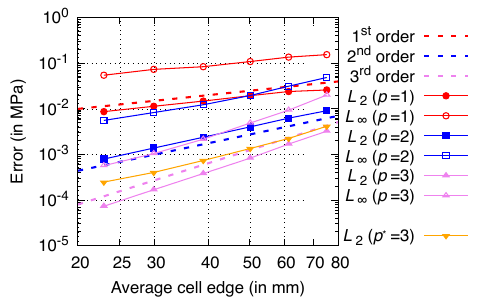}
    }
 	\caption{Spherical cavity case: displacement and stress error convergence.}
 	\label{fig:sphericalCavity:convergence}
\end{figure}

Fig.~\ref{fig:sphericalCavity:b} shows the predicted von Mises stress field on the finest mesh. The convergence of the $L_2$ and $L_\infty$ norms of the displacement and stress discretisation errors is presented in Fig.~\ref{fig:sphericalCavity:convergence}. This Fig.~also includes convergence of the $L_2$ norm computed by prescribing a zero-traction condition on the cavity boundary for the $p=3$ interpolation (denoted by the yellow $p^*$ curves). Corresponding results for $p=1$ and $p=2$ are omitted for clarity, as they closely match existing solutions. Both displacement and stress convergence associated with the zero-traction boundary condition exhibit second-order accuracy, consistent with the geometric approximation of the curved cavity surface.

%
%
%------------------------------------------------------------------------------
\subsubsection{Case 7: Vibration of a 3D Cantilevered Beam}
%------------------------------------------------------------------------------
%
In this case, compressible bending of a 3D thick column under dynamic response is considered \cite{Tukovic2007, Cardiff2025}.
The beam is clamped on its bottom face, and its upper face is subjected to a sudden traction $T = (50, 50, 0)$ kPa.
A neo-Hookean hyperelastic material is assumed with $E = 15.293\ \text{MPa}$, $\nu = 0.3$, and density $\rho = 1000$ $\text{kg/m}^3$.
The beam dimensions are $0.2\,\text{m} \times 0.2\, \text{m}$ and $6\,\text{m}$ in length.
$p$-refinement is performed on the mesh containing $1250$ cells, using a $1$ ms time step and a total simulation time of $1$ s.
The inertia term in Eq.~\eqref{eq:inertiaTerm} is integrated with the second-order backward Euler (BDF2) scheme and one integration point per cell, see Appendix \ref{app:time}.
%The main objective of this problem is to illustrate the robustness of the solution procedure when large rotations and strains are present.

\begin{figure}[H]
 	\centering
 	\subfigure[Displacement magnitude of the centre of the upper surface of the beam vs. time.]
 	{
 		\label{fig:vibratingCantilever:a}
    	\includegraphics[scale=0.8]{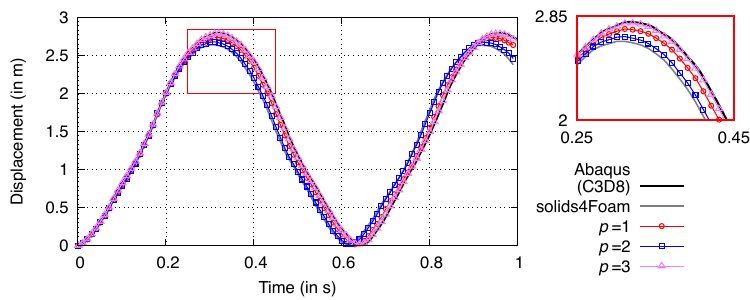}
    }
    \quad
 	\subfigure[Displacement magnitude for $t=0.3$ s and $p=3$]
 	{
 		\label{fig:vibratingCantilever:b}
    	\includegraphics[scale=0.8]{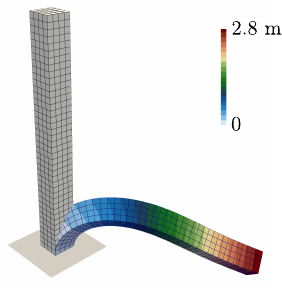}
    }
 	\caption{Vibration of a 3D cantilevered beam.}
 	\label{fig:vibratingCantilever}
\end{figure}

The corresponding displacement magnitude at the centre of the upper surface as a function of time is shown in Fig.~\ref{fig:vibratingCantilever:a}.
A visualisation of the maximum displacement magnitude is provided in Fig.~\ref{fig:vibratingCantilever:b} for $t = 0.3$ s, at which point the upper surface has deflected below the fixed end.
From Fig.~\ref{fig:vibratingCantilever:a}, it can be seen that the predictions converge toward those obtained using Abaqus with the C3D8 element on a mesh containing $138\;240$ elements.
The results for $p = 3$ lie on top of the Abaqus reference solution, while the $p = 2$ and $p = 1$ results show slight offset, with the $p = 1$ solution closer to the $p = 3$ result.
Results obtained using \texttt{solids4foam} are included for reference; they are comparable to the $p = 2$ predictions and also converge toward the Abaqus solution upon mesh refinement \cite{Cardiff2025}.

\subsection{Error-Cost Relationship}
%------------------------------------------------------------------------------
%
To properly assess efficiency, the relevant metric is the computational cost required to achieve a given error level \cite{wang2013}. 
In other words, comparing execution time or memory usage on a fixed mesh is not meaningful, since different method orders produce different error levels on the same mesh. Accordingly, efficiency is evaluated by comparing the $L_2$ error against the total CPU time required for a complete simulation.
All reported results are obtained from serial executions and measured using the GNU time utility on macOS, running on an Apple M4 Pro CPU with the Clang 17.0.0 compiler.

\begin{figure}[H]
 	\centering
	\subfigure[Displacement error]
 	{
 		\label{fig:error-time-MMS-struct-hex:a}
    	\includegraphics[scale=0.8]{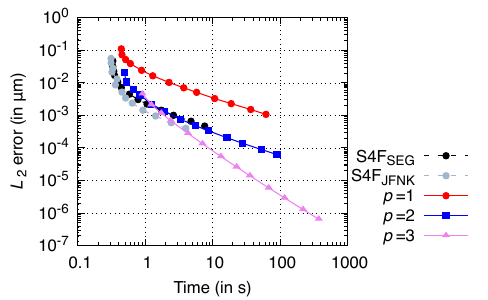}
    }
 	\subfigure[Stress error]
 	{
 		\label{fig:error-time-MMS-struct-hex:b}
    	\includegraphics[scale=0.8]{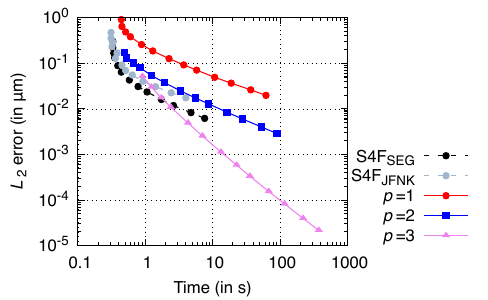}
    }
 	\caption{Error-time relationship for 3D MMS case, regular hexahedral mesh.}
 	\label{fig:error-time-MMS-struct-hex}
\end{figure}

\begin{figure}[H]
 	\centering
	\subfigure[Displacement error]
 	{
 		\label{fig:error-time-MMS-struct-tet:a}
    	\includegraphics[scale=0.8]{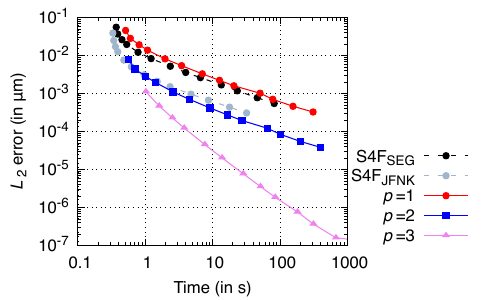}
    }
 	\subfigure[Stress error]
 	{
 		\label{fig:error-time-MMS-struct-tet:b}
    	\includegraphics[scale=0.8]{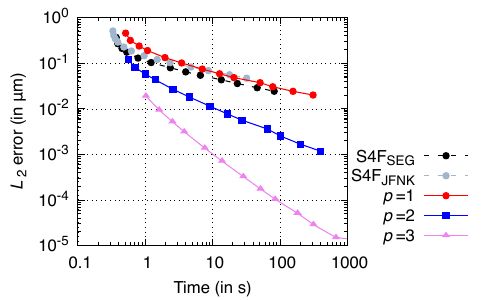}
    }
 	\caption{Error-time relationship for 3D MMS case, regular tetrahedral mesh.}
 	\label{fig:error-time-MMS-struct-tet}
\end{figure}

\begin{figure}[H]
 	\centering
	\subfigure[Displacement error]
 	{
 		\label{fig:error-time-MMS-struct-poly:a}
    	\includegraphics[scale=0.8]{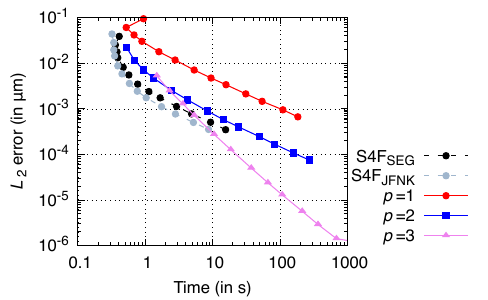}
    }
 	\subfigure[Stress error]
 	{
 		\label{fig:error-time-MMS-struct-poly:b}
    	\includegraphics[scale=0.8]{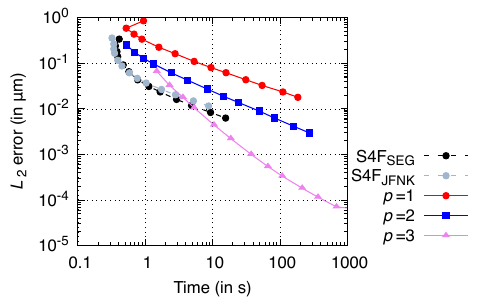}
    }
 	\caption{Error-time relationship for 3D MMS case, regular polyhedral mesh.}
 	\label{fig:error-time-MMS-struct-poly}
\end{figure}

Figs.~\ref{fig:error-time-MMS-struct-hex},~\ref{fig:error-time-MMS-struct-tet},~and \ref{fig:error-time-MMS-struct-poly} present the error–cost (error-time) relationship for the three-dimensional manufactured solution described in Section~\ref{case:mms}, considering three mesh types: hexahedral, tetrahedral, and polyhedral. 
Each Fig.~shows error-cost convergence for both the displacement (a) and the stress (b). Results are reported for the conventional solids4foam (S4F) segregated solver (denoted 
S4F$_{\text{SEG}}$) and the solids4foam Jacobian-Free Newton–Krylov solver (denoted S4F$_{\text{JFNK}}$). Both solids4foam solvers are second-order accurate; their accuracy is therefore comparable to that of the high-order discretisation developed in this work with $p=1$. The timings obtained from these solvers provide a meaningful baseline for comparison, as the segregated solver is built upon a mature and highly optimised OpenFOAM infrastructure.

\begin{figure}[H]
 	\centering
	\subfigure[Regular hexahedral mesh]
 	{
 		\label{fig:error-time-MMS-struct-memory:a}
    	\includegraphics[scale=0.8]{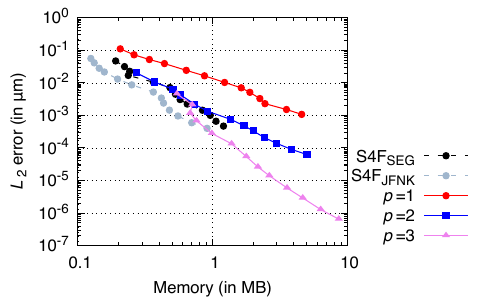}
    }
 	\subfigure[Regular tetrahedral mesh]
 	{
 		\label{fig:error-time-MMS-struct-memory:b}
    	\includegraphics[scale=0.8]{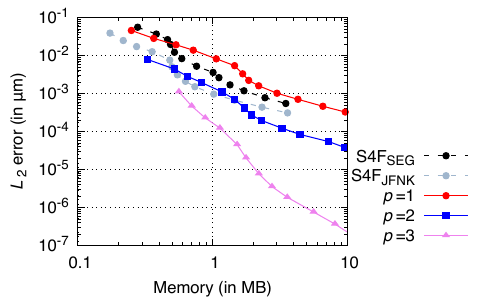}
    }
 	\caption{Error-memory relationship for 3D MMS case.}
 	\label{fig:error-memory}
\end{figure}

From Fig.~\ref{fig:error-time-MMS-struct-hex}, it can be observed that, for the hexahedral mesh, the $p=2$ interpolation reaches a level of efficiency comparable to that of the solids4foam solvers. It is important to note that, on regular hexahedral meshes, stress calculations using the solids4foam solvers exhibit convergence rates that are slightly higher than expected \cite{Syrakos2023, Cardiff2025}. This behaviour is not observed for tetrahedral meshes (Fig.~\ref{fig:error-time-MMS-struct-tet}), where $p=1$ is sufficient to achieve comparable efficiency for stress, while $p=2$ is required to reach similar efficiency in displacement for the JFNK solver. Results for the polyhedral mesh, shown in Fig.~
\ref{fig:error-time-MMS-struct-poly}, indicate that $p=2$ is insufficient to match the efficiency achieved for tetrahedral and hexahedral meshes. This is attributed to the larger number of quadrature points required per face in polyhedral meshes due to their geometric complexity, which increases the computational cost. As a result, face–edge connectivity has a direct impact on efficiency, making tetrahedral meshes the most efficient choice among those considered, owing to their minimal face–edge connectivity.
In all cases, the wall-clock time increases exponentially with decreasing error, resulting in an approximately linear trend on a log–log plot. Similar slopes are observed for 
$p=1$, $p=2$, and the S4F solvers, indicating comparable asymptotic efficiency. An exception is observed for the stress results on tetrahedral meshes, where distinct slopes arise for different interpolation orders, reflecting differences in their asymptotic convergence behaviour. Across all mesh types, the $p=3$ discretisation demonstrates a clear efficiency advantage once the mesh is sufficiently refined. Beyond the crossover point at which the efficiency curves of $p=2$ (or the S4F solvers) intersect that of $p=3$, higher-order discretisation yields performance that cannot be matched by lower-order methods.
Fig.~\ref{fig:error-memory} shows the maximal memory usage for tetrahedral and hexahedral meshes. The observed memory consumption follows trends consistent with the error–cost results presented in Figs.~\ref{fig:error-time-MMS-struct-hex}--\ref{fig:error-time-MMS-struct-poly}, indicating that trends in memory efficiency mirror those observed for time efficiency. Consequently, consistent conclusions can be drawn regarding the performance in terms of memory usage.
The reported timings should be regarded as indicative and may evolve as the implementation develops. Nevertheless, based on the current formulation, it is expected that the efficiency of the S4F solvers will lie between those of the $p=1$ and $p=2$ interpolation. In particular, the $p=1$ method is not expected to outperform S4F, as it requires the storage and computation of interpolation coefficients.
%
%------------------------------------------------------------------------------
\subsection{Stabilisation Scheme Effects}
%------------------------------------------------------------------------------
%
This section analyses the effect of the scaling factor, $\alpha$, which is a user-defined parameter in the stabilisation scheme (Eq.~\eqref{eq:alpha}) used to suppress spurious zero-energy modes. While $\alpha=0.1$ is used as the default value for results in all examples, the investigation here focuses on how a different choice of $\alpha$ impacts solution accuracy and computational timing.
The effect on timing is indirect: while the stabilisation calculation itself has a constant computational cost, different $\alpha$ values alter the properties of the system matrix, in particular the matrix condition number. 
This change in system behavior directly affects the linear solver's performance. 
Accordingly, selecting $\alpha$ involves a balance between introducing numerical dissipation (which impacts accuracy) and maintaining solver efficiency and robustness. 
To quantify this trade-off, the Method of Manufactured Solutions (case \ref{case:mms}) is used with structured hexahedral and unstructured tetrahedral meshes.

\begin{figure}[H]
 	\centering
	\subfigure[Structured hexahedral mesh]
 	{
 		\label{fig:alpha-disp:a}
    	\includegraphics[scale=0.8]{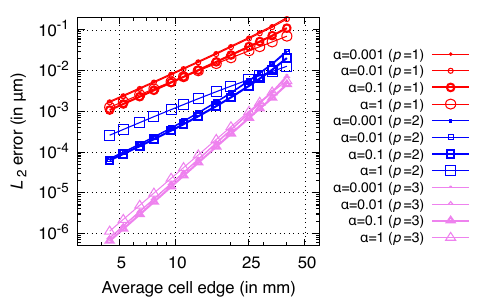}
    }
 	\subfigure[Unstructured tetrahedral mesh]
 	{
 		\label{fig:alpha-disp:b}
    	\includegraphics[scale=0.8]{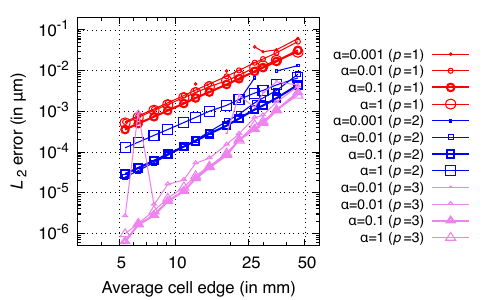}
    }
 	\caption{Convergence of the $L_2$ displacement error norm for case \ref{case:mms} with varying stabilisation parameter: $\alpha = 0.001, 0.01, 0.1,$ and $1.0$.}
 	\label{fig:alpha:disp-mms}
 \end{figure}

\begin{figure}[H]
 	\centering
	\subfigure[Structured hexahedral mesh]
 	{
 		\label{fig:alpha-timing:a}
    	\includegraphics[scale=0.8]{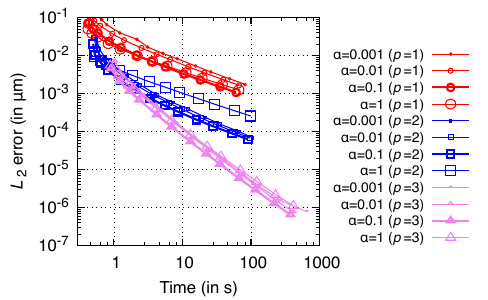}
    }
 	\subfigure[Unstructured tetrahedral mesh]
 	{
 		\label{fig:alpha-timing:b}
    	\includegraphics[scale=0.8]{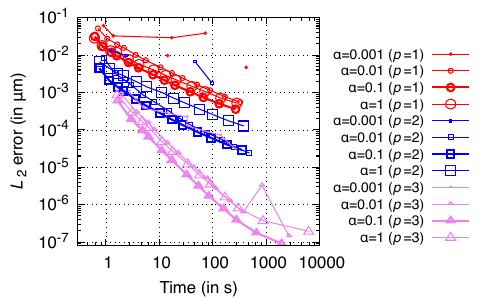}
    }
 	\caption{$L_2$ displacement error norm versus computational time for case \ref{case:mms}, comparing the effect of $\alpha = 0.001, 0.01, 0.1,$ and $1.0$.}
 	\label{fig:alpha-timing}
 \end{figure}

Four values for the scaling factor, 
each differing by an order of magnitude, are investigated: $\alpha = 0.001, 0.01, 0.1,$ and $1.0$. 
Fig.~\ref{fig:alpha:disp-mms} presents the $L_2$ error norm convergence for displacement, while Fig.~\ref{fig:alpha-timing} presents computational timings. 
Different values of $\alpha$ do not affect the convergence order; the exception is for $p=2$ (hexahedral mesh), where $\alpha=1$ forces an earlier transition to second-order. 
In terms of robustness, $\alpha=0.001$ is too low to stabilise the unstructured (irregular) mesh, and the solver converges only on some meshes. 
From these diagrams, one can see that $\alpha=0.1$ and $\alpha=0.01$ are the best values in terms of accuracy. However, when looking at Fig.~\ref{fig:alpha-timing}, 
it is clear that $\alpha=0.1$ is the optimal choice for efficiency. 
It should be emphasised that, $\alpha=0.1$, is the same as the one for Rhie-Chow stabilisation in standard second-order solvers \cite{Cardiff2025}.This agreement arises because, for $p=1$, the Rhie–Chow stabilisation coincide to the $\alpha$-stabilisation.

As an alternative to the proposed $\alpha$-stabilisation, a comparable stabilising effect can be achieved by increasing the interpolation stencil and employing more integration points than the minimum required, an approach often referred to as \emph{over-integration} or \emph{antialiasing}. 
This strategy has been adopted in previous studies \citep{Castrillo2022, Castrillo2024}. 
However, its main drawbacks are the higher computational cost, arising from both the enlarged stencil and the increased number of integration points. 
Moreover, while this approach can be effective enough for regular meshes, 
its impact is significantly reduced for irregularu meshes, where the present $\alpha$-stabilisation scheme demonstrates clear advantages (see Appendix \ref{app:stabilisation}).
%
%
%------------------------------------------------------------------------------
\subsection{Effect of Poor Interpolation Conditioning}
%------------------------------------------------------------------------------
%
The conditioning of the least-squares system used to compute interpolation coefficients can influence both the accuracy and robustness of the numerical solution. 
When the matrix ${\mathbf{M}}$ is highly ill-conditioned, the resulting coefficients may lead to an \emph{artificially stiff discretisation},
which can cause the linear solver to stall or converge very slowly. 
In this context, \emph{artificially stiff discretisation} means that the residual vector $\bb{R}(\bb{u})$ becomes highly sensitive to displacements of some neighbouring cells within the interpolation stencil.

The problem of ill-conditioning is of particular importance for high-order discretisations and is often regarded as one of their key challenges. While irregular meshes naturally help to alleviate ill-conditioning, the issue becomes especially severe on regular meshes, particularly on stretched regular grids, where the structured node distribution tends to amplify ill-conditioning \cite{jalali2013, Gooch2002}.

\begin{figure}[H]
 	\centering
	\subfigure[Hexahedral structured mesh usinng $n^+ = 10$ and $p=3$, coloured with matrix condition number. Red coloured cells have condition number $> 1\cdot 10^{12}$]
 	{
 		\label{fig:plateHole:hex-mesh}
    	\includegraphics[scale=0.052]{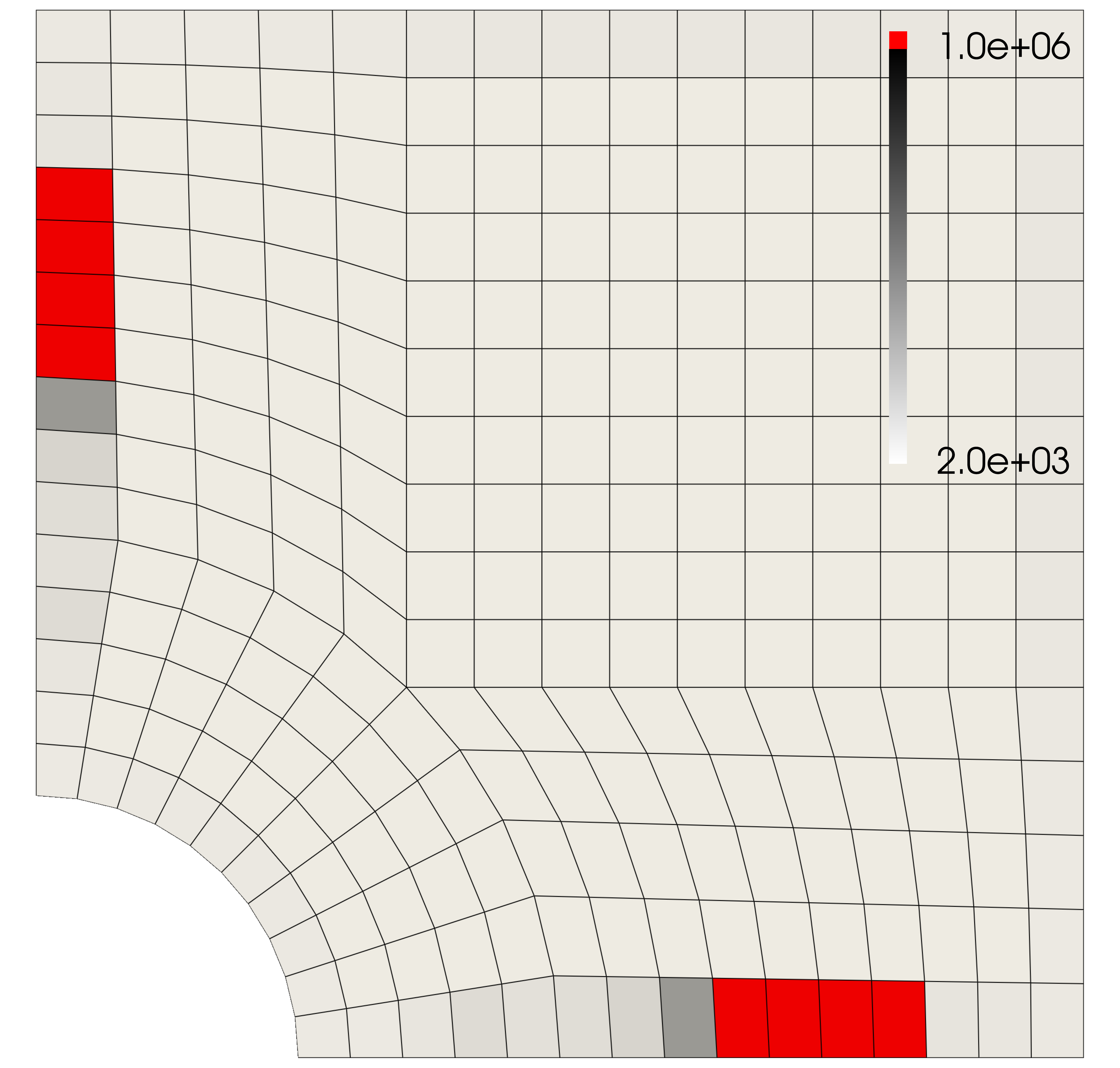}
    }
    \qquad
    \subfigure[Displacement $L_2$ norm convergence for $n^+ = 10$ and $n^+ = 20$]
 	{
 		\label{fig:plateHole-conditioning:b}
    	\includegraphics[scale=0.8]{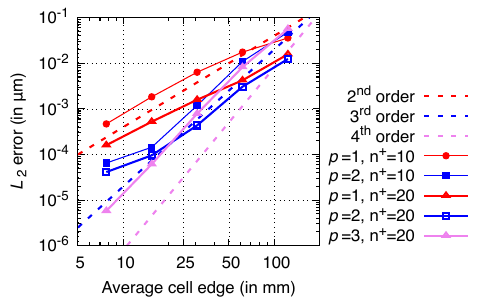}
    }
 	\caption{Plate hole case: effects of least square ill-conditioning when using regular quadrilateral mesh.}
 	\label{fig:plateHole-conditioning}
\end{figure}

The objective of this section is to illustrate the impact of this phenomenon and to demonstrate a simple remedy to mitigate it. Numerical tests reveal that ill-conditioning is most pronounced in structured two-dimensional meshes, where the absence of a third dimension prevents natural improvement of the condition number. To investigate this effect, the plate-hole benchmark from Section~\ref{case:plateHole} is considered, discretised with quadrilateral control volumes. As shown in Fig.~\ref{fig:plateHole:hex-mesh}, the condition number of the least-squares system rises sharply for the boundary cells, revealing that stencil node arrangement leads to severe ill-conditioning

Fig.~\ref{fig:plateHole-conditioning:b} presents the $L_2$-error convergence of the displacement field for two values of the stencil-expansion parameter: $n^+ = 10$, the default used in all previous examples, and $n^+ = 20$. Increasing $n^+$ effectively enlarges the interpolation stencil and reduces the condition number of the least-squares matrix ${\mathbf{M}}$, leading to improved solver robustness. As shown, with $n^+ = 10$ it is not possible to obtain results for $p = 3$ because the linear solver fails to converge, whereas for $n^+ = 20$ convergence is achieved, albeit with a reduction in the accuracy level.
It should be noted that in the case of $p = 3$ and $n^+ = 10$, where the JFNK linear solver stalls, a fully implicit discretisation combined with a direct linear solver (as in \citep{Castrillo2022, Castrillo2024}) can still produce a solution; however, the accuracy of these results is typically affected.

%Hexahedral mesh with N3 needs 20 neighbours to avoid ill conditioning. This results in lower slope for displacement but not affecting the stresses. I'm not sure why S4F results flatten on finer meshes
%
%------------------------------------------------------------------------------
\subsection{Different Choice for Approximate Jacobian}
%------------------------------------------------------------------------------
%
It has previously been demonstrated that a compact-stencil approximate Jacobian, when used as a preconditioner within the Jacobian-Free Newton–Krylov (JFNK) algorithm, is highly effective for second-order finite-volume discretisations in solid mechanics \cite{Cardiff2025}. In the present work, the same preconditioning strategy is employed in conjunction with high-order residual evaluation.

To assess the impact of the preconditioner, the computational performance of a \emph{near-exact} Jacobian (NEJ) is compared with that of the standard compact-stencil \emph{approximate} Jacobian (AJ) for the linear-elastic case, for which the NEJ can be assembled considerably more easily than in the finite-strain case.
Here, the NEJ refers to the Jacobian obtained by linearising the discretised momentum equations while excluding the stabilisation contribution
introduced in Section~\ref{sec:alpha}. As such, it does not represent the \emph{exact} Jacobian of the full residual operator, but provides a significantly closer approximation than the compact-stencil AJ and results in a substantially larger stencil.
Although the availability of such a Jacobian would, in
principle, allow the use of a conventional Newton method, Newton iterations are retained here in
order to maintain a consistent JFNK framework, with the NEJ used solely as a preconditioner. The
assembly procedure for the NEJ is described in Appendix \ref{app:P}.

Table \ref{tab:AJvsTJ-MMS} reports execution times for the manufactured solution case presented in Section~\ref{case:mms}.
For $p=1$, only minor differences in execution time are observed between the two preconditioning strategies. 
For $p=2$, the use of the NEJ leads to a performance improvement of approximately 10–15\%, while for $p=3$ a speed-up of 20–25\% is achieved. 
This behaviour is consistent with expectations, as higher-order discretisations lead to increasingly complex operators for which a more accurate Jacobian representation improves preconditioning effectiveness. These results further confirm that the compact-stencil Laplacian provides an effective and efficient preconditioning strategy for $p=1$, while remaining a reasonable choice for $p=2$ and $p=3$.

\begin{table}[h]
\centering
\begin{tabular}{ccccccc}
\multicolumn{1}{l}{} & \multicolumn{2}{c}{$p=1$} & \multicolumn{2}{c}{$p=2$} & \multicolumn{2}{c}{$p=3$} \\
Mesh size             & AJ          & NEJ           & AJ          & NEJ           & AJ          & NEJ        \\ \hline
$5\;832$                  & 3.72           & 3.29            &  5.56          & 4.63          & 22.67         &  17.47        \\
$9\;261$                 & 5.78           & 5.21            &  8.69         & 7.51         & 35.69      &  28.11            \\
$17\;576$                  & 10.71          & 10.10            &  16.36          & 14.01       & 68.02         &  53.74           \\
$29\;791$                  & 18.84           & 17.68            &  27.85          & 24.15        & 116.89    &   92.05              \\
$54\;872$                 & 36.77           & 34.70            &  52.77          & 45.91         & 219.74  &   180.33                 \\
$91\;125$                 & 61.85          & 61.86          &  89.99        & 80.12       & 376.08 &  304.09              \\ \hline
\end{tabular}
\caption{3D manufactured solution case (regular hexahedral mesh): computational time (in s) for different preconditioning matrices, comparing the AJ with the NEJ.}
\label{tab:AJvsTJ-MMS}
\end{table}

Table \ref{tab:AJvsTJ-SC} presents execution times for the spherical cavity case described in Section~\ref{case:sphericalCavity}.
In this case, the use of the NEJ as a preconditioner results in solve times that are between 2 and 92 times larger than those obtained with the AJ.
This degradation in performance is attributed to the absence of stabilisation terms in the NEJ, which leads to unfavourable numerical properties of the resulting system matrix on irregular meshes (see also Appendix \ref{app:stabilisation}). Consequently, JFNK solver performance is significantly impaired. These results demonstrate the importance of stabilisation and highlight that increased Jacobian fidelity does not necessarily translate into improved computational efficiency when the preconditioning matrix exhibits unfavourable numerical properties.

\begin{table}[h]
\centering
\begin{tabular}{ccccccc}
\multicolumn{1}{l}{} & \multicolumn{2}{c}{$p=1$} & \multicolumn{2}{c}{$p=2$} & \multicolumn{2}{c}{$p=3$} \\
Mesh size             & AJ          & NEJ           & AJ          & NEJ           & AJ          & NEJ        \\ \hline
$2\;068$                  & 1.34           & 2.82            & 1.78           & 4.18          & 4.55         &  418.33        \\
$3\;802$                 & 2.19           & 5.48            & 2.89           & 9.46         & 8.02      &  24.07             \\
$7\;030$                  & 3.94           & 11.43            & 5.23           & 16.82       & 14.80         &  48.75           \\
$14\;891$                  & 8.15           & 30.89            & 10.96           & 44.00        & 32.00    &   91.12              \\
$32\;864$                 & 18.40           & 90.45            & 24.95           & 1896.12         & 75.64  &   246.67                 \\
$72\;945$                 & 43.48          & 1183.30          & 58.71         & 3841.79        & 189.71 &  1340.57              \\ \hline
\end{tabular}
\caption{Spherical cavity case: computational time (in s) for different preconditioning matrices, comparing the AJ with the NEJ.}
\label{tab:AJvsTJ-SC}
\end{table}
%
%
%------------------------------------------------------------------------------
%\subsection{Code parallelisation}
%\label{sec:codeParallelisation}
%------------------------------------------------------------------------------
%
%\hl{To work on after meeting with Philip}
%
%%%%%%%%%%%%%%%%%%%%%%%%%%%%%%%%%%%%%%%%%%%%%%%%%%%%%%%%%%%%%%%%%%
\section{Conclusions} \label{sec:conclusion}
%%%%%%%%%%%%%%%%%%%%%%%%%%%%%%%%%%%%%%%%%%%%%%%%%%%%%%%%%%%%%%%%%%
This work presents a high-order cell-centred finite volume formulation for solid mechanics combined with a Jacobian-free Newton–Krylov nonlinear solution strategy. 
Third- and fourth-order spatial accuracy is achieved through Gaussian quadrature and a least-squares (LRE) gradient reconstruction procedure.
To ensure stability on unstructured meshes, an $\alpha$-stabilisation scheme is introduced to suppress high-frequency error modes. 
Solver efficiency is maintained through the use of a compact-stencil second-order approximate Jacobian as a preconditioner within the JFNK solver.
The formulation does not require additional boundary unknowns, supports arbitrary convex polyhedral cells and is demonstrated across a range of linear and nonlinear benchmark problems.
\\
The principal findings of the study are:
\begin{itemize}
\item \textbf{Accuracy}: The proposed discretisation delivers consistent high-order accuracy across diverse test cases. Linear reconstruction yields second-order displacement and first-order stress accuracy. Quadratic reconstruction achieves second-order accuracy for both displacement and stress; the expected third-order displacement behaviour is not observed in all cases, consistent with prior reports in the literature \cite{Nishikawa2025a}. Cubic reconstruction attains fourth-order displacement and third-order stress accuracy, confirming the effectiveness of the underlying high-order flux integration and least squares interpolation.
\item \textbf{Stabilisation scheme}: The introduced $\alpha$-stabilisation scheme, adapted from high-order CFD, effectively damps spurious high-frequency errors on unstructured meshes. The stabilisation term scales with the order of the method, unlike the Rhie–Chow correction used in second-order finite volume solvers. An optimal value of $\alpha \approx 0.1$ is observed, which is also consistent with recommendations for Rhie-Chow stabilisation \cite{Cardiff2025}.
\item \textbf{Preconditioning choice}: The compact-stencil approximate Jacobian, originating from a second-order semi-implicit discretisation, remains effective when used as a preconditioner in conjunction with high-order residual evaluation. Comparison against a fully assembled \emph{near-exact} Jacobian (linear elastic tests) confirms that compact \emph{approximate} Jacobian offers a reliable choice.
\item \textbf{Robustness}:
The method is robust across linear, nonlinear, structured, and unstructured meshes. The main challenge arises from the ill-conditioning of the local least-squares reconstruction, particularly on structured 2D grids. When conditioning is improved (e.g., via enlarged stencils), robustness is restored.
\item \textbf{Efficiency}: The OpenFOAM implementation enabled direct comparison with conventional second-order finite volume solvers. Tetrahedral meshes provide the best performance between accuracy and computational cost due to fewer quadrature points. Efficiency decreases with increasing face-edge count because this enlarges the quadrature stencil and increases storage requirements for interpolation coefficients. Despite a non-optimised implementation, the method shows competitive performance and substantial potential.
\item \textbf{Integration into existing JFNK frameworks}: The high-order formulation requires only modifications to the residual evaluation and therefore integrates seamlessly into existing second-order JFNK infrastructure. The complete implementation is released within \texttt{solids4foam}, enabling community use, evaluation, and extension.
\end{itemize}
Future work will be directed in several key areas. First, the treatment of curved boundaries in the context of high-order finite volume for solid mechanics is a critical and unresolved challenge that warrants special attention. 
Second, more efficient high-order implementations should be explored to further reduce computational cost. Efficiency gains may be achieved by reducing the number of least-squares systems, for example by forming a single reconstruction per face, per cell, or at mesh vertices, as in the economical approach proposed by Nishikawa \cite{Nishikawa2025}. This strategy has the potential to significantly improve both time and memory efficiency while preserving the overall order of accuracy.
Finally, the successful application of the high-order JFNK algorithm to solid mechanics opens a promising avenue for its extension to a monolithic, high-order fluid–solid interaction framework.
%
%
%
%%%%%%%%%%%%%%%%%%%%%%%%%%%%%%%%%%%%%%%%%%%%%%%%%%%%%%%%%%%%%%%%%%%%%%%%%%%%%%%
\backmatter
%
%%%%%%%%%%%%%%%%%%%%%%%%%%%%%%%%%%%%%%%%%%%%%%%%%%%%%%%%%%%%%%%%%%%%%%%%%%%%%%%
\bmhead{Data Availability}
%%%%%%%%%%%%%%%%%%%%%%%%%%%%%%%%%%%%%%%%%%%%%%%%%%%%%%%%%%%%%%%%%%%%%%%%%%%%%%%
%
The codes presented are publicly available at \url{https://github.com/solids4foam/solids4foam}, and the cases and plotting scripts are available at \url{https://github.com/solids4foam/solid-benchmarks}.
%
%%%%%%%%%%%%%%%%%%%%%%%%%%%%%%%%%%%%%%%%%%%%%%%%%%%%%%%%%%%%%%%%%%%%%%%%%%%%%%%
\bmhead{Declaration of generative AI and AI-assisted technologies in the manuscript preparation process.}
%%%%%%%%%%%%%%%%%%%%%%%%%%%%%%%%%%%%%%%%%%%%%%%%%%%%%%%%%%%%%%%%%%%%%%%%%%%%%%%
%
During the preparation of this work, the authors used ChatGPT and Grammarly as writing assistants.
After using these tools, the authors reviewed and edited the content as needed and takes full responsibility for the content of the publication.
%
%%%%%%%%%%%%%%%%%%%%%%%%%%%%%%%%%%%%%%%%%%%%%%%%%%%%%%%%%%%%%%%%%%%%%%%%%%%%%%%
\bmhead{Acknowledgments}
%%%%%%%%%%%%%%%%%%%%%%%%%%%%%%%%%%%%%%%%%%%%%%%%%%%%%%%%%%%%%%%%%%%%%%%%%%%%%%%
%
This project has received funding from the European Research Council (ERC) under the European Union’s Horizon 2020 research and innovation programme (Grant Agreement No. 101088740).
Financial support is gratefully acknowledged from I-Form, funded by Research Ireland (formerly Science Foundation Ireland)
Grant Number {21/RC/10295\_P2}, co-funded under European Regional Development Fund and by I-Form industry partners, and from NexSys, funded by Research Ireland Grant Number 21/SPP/3756.
Additionally, the authors wish to acknowledge the DJEI/DES/SFI/HEA Irish Centre for High-End Computing (ICHEC) for the provision of computational facilities and support (www.ichec.ie), and part of this work has been carried out using the UCD ResearchIT Sonic cluster which was funded by UCD IT Services and the UCD Research Office. None of the funding sources mentioned above had any role in study design, data collection and analysis, decision to publish, or preparation of the manuscript.
\newpage
\setcounter{figure}{0}
\begin{appendices}
%
%%%%%%%%%%%%%%%%%%%%%%%%%%%%%%%%%%%%%%%%%%%%%%%%%%%%%%%%%%%%%%%%%%%%%%%%%%%%%%%
%
\section{Mechanical Laws}
\label{app:mechLaws}
%
%%%%%%%%%%%%%%%%%%%%%%%%%%%%%%%%%%%%%%%%%%%%%%%%%%%%%%%%%%%%%%%%%%%%%%%%%%%%%%%
%
\subsection{Linear Elasticity}
The definition of engineering stress $\bb{\sigma}_\text{s}$ for linear elasticity can be given as:
\begin{eqnarray}  \label{eq:linearElastic}
	\bb{\sigma}_\text{s}
	&=& 2 \mu \bb{\varepsilon} + \lambda \, \text{tr} \left( \bb{\varepsilon} \right) \bb{I} \notag, \\
	&=& \mu \bb{\nabla} \bb{u} + \mu \left( \bb{\nabla} \bb{u}\right)^{\text{T}} + \lambda \left(\bb{\nabla} \cdot \bb{u} \right) \bb{I},
\end{eqnarray}
where $\lambda$ is the first Lam\'{e} parameter, and $\mu$ is the second Lam\'{e} parameter, synonymous with the shear modulus.
The Lam\'{e} parameters can be expressed in term of the Young's modulus ($E$) and Poisson's ratio $\nu$ as:
\begin{eqnarray}
	\mu = \frac{E}{2(1 + \nu)}, \quad \lambda = \frac{E \nu}{(1+\nu)(1 - 2\nu)}.
\end{eqnarray}
%
%------------------------------------------------------------------------------
%\subsection{St.\,Venant-Kirchoff Hyperelasticity}
%%------------------------------------------------------------------------------
%%
%The St.\,Venant-Kirchoff model defines the second Piola–Kirchhoff stress $\textbf{S}$ as
%\begin{eqnarray}
%	\bb{S} &=& 2 \mu \bb{E} + \lambda \, \text{tr} \left( \bb{E} \right) \textbf{I}
%\end{eqnarray}
%where, as before, $\lambda$ is the first Lam\'{e} parameter, and $\mu$ is the second Lam\'{e} parameter.
%The Lagrangian Green strain $\textbf{E}$ is defined as
%\begin{eqnarray}
%	\bb{E} &=& \frac{1}{2} \left( \bb{\nabla} \bb{u} + \bb{\nabla} \bb{u}^{\text{T}} + \bb{\nabla} \bb{u} \cdot \bb{\nabla} \bb{u}^{\text{T}}  \right)
%\end{eqnarray}
%
%The true stress can be calculated from the second Piola–Kirchhoff stress as
%\begin{eqnarray} \label{eq:S2sigma}
%	\bb{\sigma} &=& \frac{1}{J} \bb{F} \cdot \bb{S} \cdot \bb{F}^{\text{T}}
%\end{eqnarray}
%
%------------------------------------------------------------------------------
\subsection{Neo-Hookean Hyperelasticity}
\label{app:mechLaws:NH}
%------------------------------------------------------------------------------
%
The definition of true (Cauchy) stress $\bb{\sigma}$ for neo-Hookean hyperelasticity can be given as:
\begin{eqnarray} \label{eq:neoHook}
	\bb{\sigma}
	&=& \frac{\mu}{J} \, \text{dev} \left( \bar{\bb{b}} \right) + \frac{\kappa}{2} \frac{J^2 - 1}{J} \bb{I},
\end{eqnarray}
where, once again, $\mu$ is the shear modulus, and $\kappa$ is the bulk modulus.
The bulk modulus can be expressed in terms of the Young's modulus ($E$) and Poisson's ratio $\nu$ as:
\begin{eqnarray}
	\kappa = \frac{E}{3(1 - 2\nu)}.
\end{eqnarray}
The volume-preserving component of the elastic left Cauchy--Green deformation tensor $\boldsymbol{b}$ is given as
\begin{eqnarray}
	\bar{\bb{b}} = J^{-2/3} \bb{b} = J^{-2/3} \bb{F} \cdot \bb{F}^{\text{T}}.
\end{eqnarray}
In the limit of small deformations $\lVert \nabla \mathbf{u} \rVert  \ll 1$, neo-Hookean hyperelasticity (Eq.~\eqref{eq:neoHook}) reduces to linear elasticity (Eq.~\eqref{eq:linearElastic}).
\subsection{Mooney-Rivlin Hyperelasticity} 
\label{app:mechLaws:MR}
%------------------------------------------------------------------------------
%
Strain-energy function in three-parameter Mooney-Rivlin hyperelastic material \cite{Rivlin1951}:
\begin{eqnarray}
\Psi = c_{10}(I_1-3) + c_{01}(I_2-3) + c_{11}(I_1-3)(I_2-3),
\end{eqnarray}
where $c_{10}$, $c_{01}$ and $c_{11}$ are material parameters, $I_1$ and $I_2$ are first and second principal invariant of the left Cauchy-Green deformation.

The Cauchy stress tensor for the loosely incompressible case is then given by \cite{oliveira2020}:
\begin{eqnarray}
\bb{\sigma} = \frac{1}{J} \left[ \dfrac{1}{2}\kappa(J^2-1)\bb{I} + 2(c_{10}+c_{11}(I_2-3)*\text{dev}(\bb{B}_{\mathrm{iso}})) -2(c_{01}+c_{11}(I_1-3)\text{dev}(\bb{B}_{\mathrm{iso}}^{-1}))\right].
\end{eqnarray}
where $\bb{B}_{\mathrm{iso}}=J^{-2/3}\bb{F}\cdot\bb{F}^{\mathrm{T}}$ is isochoric left Cauchy-Green deformation tensor and $\kappa$ is \emph{penalty parameter} which value is typically taken to be equal to the material bulk modulus.
%
%%%%%%%%%%%%%%%%%%%%%%%%%%%%%%%%%%%%%%%%%%%%%%%%%%%%%%%%%%%%%%%%%%%%%%%%%%%%%%%
%
\section{High-order Derivatives}
\label{app:derivatives}
%
%%%%%%%%%%%%%%%%%%%%%%%%%%%%%%%%%%%%%%%%%%%%%%%%%%%%%%%%%%%%%%%%%%%%%%%%%%%%%%%
%
The reconstruction of higher-order derivatives is performed independently for each displacement component $\alpha$ of the vector field
$\bb{u} = [u_x,\,u_y,\,u_z]^{\mathrm{T}}$.
For a scalar field $u_{\alpha}$ representing a single displacement component, the second derivative at the cell centre $P$ is reconstructed as a weighted sum of neighbouring cell values contained in the stencil $\mathcal{S}_P$:
\begin{eqnarray}
(\nabla \nabla {u_{\alpha}})_{P} = \left( \dfrac{\partial^2 u_{\alpha}}{\partial x_i \partial x_j} \right)_P
= \sum_{N \in \mathcal{S}_P} \bb{c}_{{\bb{x}\bb{x}}_N} ({u}_{\alpha})_N,
\label{eq:hessianComponent}
\end{eqnarray}
where $\bb{c}_{{\bb{x}\bb{x}}_N}$ is a symmetric second-order coefficient tensor:
\begin{eqnarray}
\bb{c}_{{\bb{x}\bb{x}}_N} =
\begin{bmatrix}
c_{xx_N} & c_{xy_N} & c_{xz_N} \\
c_{xy_N} & c_{yy_N} & c_{yz_N} \\
c_{xz_N} & c_{yz_N} & c_{zz_N}
\end{bmatrix}.
\label{eq:secondTensor}
\end{eqnarray}
Coefficients in tensor $\bb{c}_{{\bb{x}\bb{x}}_N}$ are derived from the rows of the reconstruction matrix $\mathbf{\bar{A}}$ corresponding to the second-order monomials.
The second-derivative displacement tensor $(\nabla \nabla \bb{u})$ is assembled from the derivatives of the three scalar displacement components.
For $p = 3$, the third-derivative tensor $(\nabla \nabla \nabla \boldsymbol{u})_P$ is reconstructed for each displacement component as:
\begin{eqnarray}
(\nabla \nabla \nabla {u}_{\alpha})_P = \left( \dfrac{\partial^3 u_{\alpha}}{\partial x_i \partial x_j \partial x_k} \right)_P
= \sum_{N \in \mathcal{S}_P} \mathbf{c}_{\bb{x}\bb{x}\bb{x}_N}({u}_{\alpha})_N,
\label{eq:thirdTensor}
\end{eqnarray}
where $\mathbf{c}_{\bb{x}\bb{x}\bb{x}_N}$ is a symmetric third-order coefficient tensor:
\begin{eqnarray}
\mathbf{c}_{\bb{x}\bb{x}\bb{x}_N} =\big[c_{xxx_N},c_{xxy_N},c_{xxz_N},c_{xyy_N},c_{xyz_N},c_{xzz_N},c_{yyy_N},c_{yyz_N},c_{yzz_N},c_{zzz_N}\big]^{\mathrm{T}},
\end{eqnarray}
derived from the rows of the matrix $\mathbf{\bar{A}}$.

The second- and third-order terms appearing in Eq.~\eqref{eq:alphaJump}, namely $\bb{d}_{Nf}^2 \!:\! (\nabla \nabla \bb{u})_N$ and $\bb{d}_{Nf}^3 \!::\! (\nabla \nabla \nabla \bb{u})_N$, are vectors whose components are obtainer through successive tensor contraction performed independently for each displacement component:
\begin{eqnarray}
\begin{split}
 &\left(\bb{d}^2 \!:\! \nabla \nabla \bb{u} \right)_{\alpha} = \sum_{i,j} d_i d_j \dfrac{\partial^2 u_{\alpha}}{\partial x_i \partial x_j}, \\
 & \left(\bb{d}^3 \!::\! \nabla \nabla \nabla \bb{u} \right)_{\alpha} = \sum_{i,j,k} d_i d_j d_k \dfrac{\partial^3 u_{\alpha}}{\partial x_i \partial x_j \partial x_k},
\end{split}
\end{eqnarray}
where indices $i,j,k$ denote spatial coordinates.

To mitigate ill-conditioning in the reconstruction procedure, the polynomial basis $\mathbf{q}$ in Eq.~\eqref{eq:qAnda} is scaled using the characteristic stencil size $h = 2r_s$:
\begin{eqnarray}
\mathbf{q}^T(\bb{x}-\tilde{\bb{x}}) = \left[1, \frac{(x-\tilde{x})}{h}, \frac{(y-\tilde{y})}{h}, \frac{(z-\tilde{z})}{h}, \frac{1}{2}\frac{(x-\tilde{x})^2}{h^2}, \dots \right].
\end{eqnarray}
Accordingly, the monomials and their corresponding coefficients in Eqs.~\eqref{eq:hessianComponent}~and~\eqref{eq:thirdTensor} are scaled by $h^{-2}$ for the second derivative terms and by $h^{-3}$ for the third derivative terms.
%
%
%%%%%%%%%%%%%%%%%%%%%%%%%%%%%%%%%%%%%%%%%%%%%%%%%%%%%%%%%%%%%%%%%%%%%%%%%%%%%%%
%
\section{Preconditioner Based on Near-exact Jacobian}
\label{app:P}
%
%%%%%%%%%%%%%%%%%%%%%%%%%%%%%%%%%%%%%%%%%%%%%%%%%%%%%%%%%%%%%%%%%%%%%%%%%%%%%%%
%
In the case of a linear elastic material, the near-exact Jacobian is constructed as follows:
\begin{eqnarray} \label{app:tj}
\begin{split}
 \sum_{f \,\in\, \mathcal{F}_P} \bb{n}_{f} \cdot &\left[  \sum_{q=1}^{q=N_{f,q}}\alpha_q \; (\bb{\sigma}_\text{s})_{f,q} \right] \Gamma_f = \\
 &\sum_{f \,\in\, \mathcal{F}_P^{\text{int}}} \bb{n}_{f} \cdot \left[  \sum_{q=1}^{q=N_{f,q}}\alpha_q \; (\bb{\sigma}_\text{s})_{f,q} \right] \Gamma_f
 +  \sum_{b \,\in\, \mathcal{F}_P^{\text{non-trac}}} \bb{n}_{b} \cdot \left[  \sum_{q=1}^{q=N_{f,q}}\alpha_q \; (\bb{\sigma}_\text{s})_{b,q} \right] \Gamma_b
\end{split}
\end{eqnarray}
The first term in Eq.~\eqref{app:tj},  corresponding to internal faces, is further discretised following \cite{Castrillo2022, Castrillo2024}:
\begin{eqnarray} 
\begin{split}
&\sum_{f \,\in\, \mathcal{F}_P^{\text{int}}} \bb{n}_{f} \cdot \left[  \sum_{q=1}^{q=N_{f,q}}\alpha_q \; (\bb{\sigma}_\text{s})_{f,q} \right] \Gamma_f\\
& = \sum_{f \,\in\, \mathcal{F}_P^{\text{int}}} \bb{n}_{f} \cdot \left [ \sum_{q=1}^{q=N_{f,q}}\alpha_q \; \left(\mu (\nabla \bb{u})_{f,q} + \mu (\nabla \bb{u})^T_{f,q} +\lambda\, \text{tr}((\nabla \bb{u})_{f,q} )\bb{I}\right) \right]\Gamma_f, \\
=&
\sum_{f \,\in\, \mathcal{F}_P^{\text{int}}} \left [ \sum_{q=1}^{q=N_{f,q}}\alpha_q \; \left(\mu \sum_{N\,\in\,\mathcal{S}_f}(\bb{c}_{x_N}\cdot \bb{n}_f)\bb{I} \bb{u}_N 
 + \mu \sum_{N\,\in\,\mathcal{S}_f}(\bb{c}_{x_N} \otimes \bb{n}_f) \bb{u}_N 
 +\lambda  \sum_{N\,\in\,\mathcal{S}_f}(\bb{n}_f \otimes \bb{c}_{x_N}) \bb{u}_N\right) \right]\Gamma_f.
\end{split}
\end{eqnarray}
The second term corresponds to symmetry and displacement boundaries (traction boundaries are omitted, as their contribution is moved to the right-hand-side source vector):
\begin{eqnarray}
\begin{split}
\sum_{b \,\in\, \mathcal{F}_P^{\text{non-trac}}}& \bb{n}_{b} \cdot \left[  \sum_{q=1}^{q=N_{f,q}}\alpha_q \; (\bb{\sigma}_\text{s})_{b,q} \right] \Gamma_b = \\
&=\sum_{b \,\in\, \mathcal{F}_P^{\text{symm}}} \bb{n}_{b} \cdot \left[  \sum_{q=1}^{q=N_{f,q}}\alpha_q \; (\bb{\sigma}_\text{s})_{b,q} \right] \Gamma_b 
+ \sum_{b \,\in\, \mathcal{F}_P^{\text{disp}}} \bb{n}_{b} \cdot \left[  \sum_{q=1}^{q=N_{f,q}}\alpha_q \; (\bb{\sigma}_\text{s})_{b,q} \right] \Gamma_b.\\
\end{split}
\end{eqnarray}
Faces on the displacement boundary $\mathcal{S}_f^{\text{disp}}$ are discretised in the same way as internal faces.
The ghost node contribution on displacement boundaries is included to the right-hand side vector.
For symmetry faces, the stencil is split into the physical part $\mathcal{F}_P^{\text{symm,p}}$ and the ghost part $\mathcal{F}_P^{\text{symm,g}}$.
The physical part is again discretised as internal faces, while the ghost contribution includes the reflection tensor $\bb{R}$:
\begin{eqnarray}
\begin{split}
\sum_{b \,\in\, \mathcal{F}_P^{\text{symm,g}}} \bb{n}_{f} \cdot &\left[  \sum_{q=1}^{q=N_{f,q}}\alpha_q \; (\bb{\sigma}_\text{s})_{b,q} \right] \Gamma_b 
= \\
& \sum_{b \,\in\, \mathcal{F}_P^{\text{symm,g}}} \left [ \sum_{q=1}^{q=N_{f,q}}\alpha_q \; \left(\mu \sum_{N\,\in\,\mathcal{S}_f^{\text{symm,g}}}(\bb{c}_{x_N}\cdot \bb{n}_b)\bb{I}\cdot \bb{R} \cdot \bb{u}_N  \right) \right] 
\\
&  
 +\sum_{b \,\in\, \mathcal{F}_P^{\text{symm,g}}} \left [ \sum_{q=1}^{q=N_{f,q}}\alpha_q \; \left(\mu \sum_{N\,\in\,\mathcal{S}_f^{\text{symm,g}}}(\bb{c}_{x_N} \otimes \bb{n}_b) \cdot \bb{R} \cdot \bb{u}_N  \right) \right] \\
& 
+\sum_{b \,\in\, \mathcal{F}_P^{\text{symm,g}}} \left [ \sum_{q=1}^{q=N_{f,q}}\alpha_q \; \left(\lambda \sum_{N\,\in\,\mathcal{S}_f^{\text{symm,g}}}(\bb{n}_b \otimes \bb{c}_{x_N}) \cdot \bb{R} \cdot \bb{u}_N \right) \right] .
\end{split}
\end{eqnarray}
The resulting Jacobian can be used within an implicit solution procedure in which the solution is obtained in a single Newton iteration, as demonstrated in \cite{Castrillo2022}.
The Jacobian constructed here excludes the stabilisation term, which is required for irregular meshes. 
%
%
%%%%%%%%%%%%%%%%%%%%%%%%%%%%%%%%%%%%%%%%%%%%%%%%%%%%%%%%%%%%%%%%%%%%%%%%%%%%%%%
%
\section{Time integration scheme}
\label{app:time}
%
%%%%%%%%%%%%%%%%%%%%%%%%%%%%%%%%%%%%%%%%%%%%%%%%%%%%%%%%%%%%%%%%%%%%%%%%%%%%%%%
%
For the unsteady problem considered in this study, the second-order backward Euler differentiation formula (BDF2) is employed for time integration. The inertia term in Eq.~\eqref{eq:inertiaTerm} is calculated by summing contribution of cell quadrature points $\bb{x}_{\Omega, q}$ values, with acceleration computed as:
\begin{eqnarray} \label{eq:inertia2}. 	\left(\dfrac{\partial^2 \boldsymbol{u}}{\partial t^2}\right)_{\Omega, q}
 	\approx
 	\dfrac{3\boldsymbol{v}_{\Omega, q}^{[t+1]} - 4\boldsymbol{v}_{\Omega, q}^{[t]} + \boldsymbol{v}_{\Omega, q}^{[t-1]}}{2\Delta t}
 	\approx
 	\frac{3\left(
 		\dfrac{3\boldsymbol{u}_{\Omega, q}^{[t+1]} - 4\boldsymbol{u}_{\Omega, q}^{[t]} + \boldsymbol{u}_{\Omega, q}^{[t-1]}}{2\Delta t}
 		\right)
 	- 4\boldsymbol{v}_{\Omega, q}^{[t]} + \boldsymbol{v}_{\Omega, q}^{[t-1]}}{2\Delta t}
\end{eqnarray}
where $\Delta t$ is the time increment -- assumed constant here.
Superscript $[t]$ indicates the time level, with $\bb{u}_{\Omega, q}^{[t+1]}$ corresponding to the unknown displacement at the current time step.
The velocity vector $\bb{v} = \partial \bb{u}/\partial t$ at the current time step is also updated using the BDF2 scheme as:
\begin{eqnarray}
 	\boldsymbol{v}_{\Omega, q}^{[t+1]}	\approx
 		\dfrac{3\boldsymbol{u}_{\Omega, q}^{[t+1]} - 4\boldsymbol{u}_{\Omega, q}^{[t]} + \boldsymbol{u}_{\Omega, q}^{[t-1]}}{2\Delta t}.
\end{eqnarray}
Consequently, the displacement and velocity at the two previous time steps must be stored, or alternatively, the displacement at the previous  four preceding time levels.
The displacement value at a cell quadrature point, $\boldsymbol{u}_{\Omega,q}$, is obtained by extrapolation from the cell centre. A Taylor-based extrapolation, identical to that used in Eq.~\eqref{eq:alphaJump}, is employed using derivatives evaluated at the cell centre, which are already available and stored as part of the $\alpha$ -stabilisation scheme. 
The described approach integrates the temporal term in space with the same accuracy as the remaining terms in the governing equations. The BDF2 scheme is limited to second-order accuracy; however, higher-order temporal schemes or alternative time-integration operators \cite{Costa2022b} can be incorporated.
%
%%%%%%%%%%%%%%%%%%%%%%%%%%%%%%%%%%%%%%%%%%%%%%%%%%%%%%%%%%%%%%%%%%%%%%%%%%%%%%%
%
\section{Stabilisation using \emph{over-integration}}
\label{app:stabilisation}
\begin{figure}[H]
 	\centering
	\subfigure[Regular hexahedral mesh]
 	{
 		\label{fig:noalpha-disp:a}
    	\includegraphics[scale=0.8]{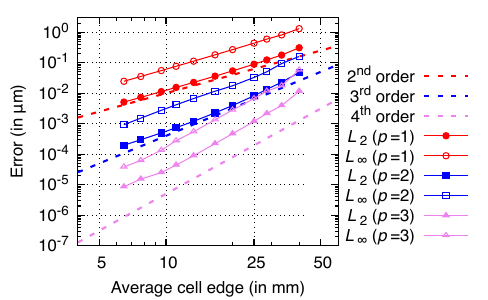}
    }
 	\subfigure[Irregular tetrahedral mesh]
 	{
 		\label{fig:noalpha-disp:b}
    	\includegraphics[scale=0.8]{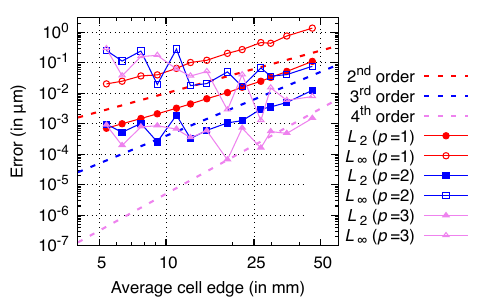}
    }
 	\caption{Convergence of the $L_2$ and $L_\infty$ displacement error norm for case \ref{case:mms} with over-integration as stabilisation.}
 	\label{fig:noalpha-disp}
 \end{figure}

%
%%%%%%%%%%%%%%%%%%%%%%%%%%%%%%%%%%%%%%%%%%%%%%%%%%%%%%%%%%%%%%%%%%%%%%%%%%%%%%%
%
Stabilisation by over-integration has been employed in previous studies \cite{Castrillo2022,Castrillo2024} to alleviate high-frequency numerical errors. For regular meshes, over-integration can be sufficient; however, on irregular meshes the stabilisation it provides is insufficient to adequately suppress such errors. This limitation is demonstrated in Figs.~\ref{fig:noalpha-disp:a}~and~\ref{fig:noalpha-disp:b} for the three-dimensional manufactured solution case presented in Section~\ref{case:mms}. Fig.~\ref{fig:noalpha-disp:a} shows results obtained on a regular hexahedral mesh, for which the solution remains stable, whereas Fig.~\ref{fig:noalpha-disp:b} presents results for an irregular tetrahedral mesh, where high-frequency oscillations persist despite the use of over-integration. This behaviour manifests as abrupt increases in the error norms and prevents the observation of a clear convergence trend, as the $L_{\infty}$ error exhibits sudden localised spikes at random locations within the computational domain.
All results were generated using 12 quadrature points per triangular element for all interpolation orders, which is sufficient to exactly integrate up to sixth-order polynomial. The stencil employed is larger than that used with the $\alpha$-stabilisation approach, with $n^{+}=60$ for $p=1$, $n^{+}=70$ for $p=2$, and $n^{+}=80$ for $p=3$. When stabilisation is insufficient, the linear solver within the JFNK solution procedure stalls and fails to converge. For this reason, the results presented in Fig.~\ref{fig:noalpha-disp} were obtained without JFNK, i.e. by solving directly the near-exact Jacobian, whose construction is described in  Appendix \ref{app:P}.
%
%%%%%%%%%%%%%%%%%%%%%%%%%%%%%%%%%%%%%%%%%%%%%%%%%%%%%%%%%%%%%%%%%%%%%%%%%%%%%%%
%
\iffalse

\section{List of things to check/discuss/improve in the paper}
%
%%%%%%%%%%%%%%%%%%%%%%%%%%%%%%%%%%%%%%%%%%%%%%%%%%%%%%%%%%%%%%%%%%%%%%%%%%%%%%%
%
\begin{itemize}
%\item[•] \textbf{Add comment:} We use openFoam normal in momentum balance when integrating over quadrature points. This assumes that face is flat. When face is not flat this affects convergence behaviour and fourth order method is affected. This is not commented in the text.
\item[•] Finish idealised ventricle case, consider adding it later.
\item[•] Consider using coarser mesh in vibration beam.
\item[•] Time scheme integration - check with Philip is it ok to leave as this
\item[•] Add other papers from 
Coppeans2026
\end{itemize}
\fi
%
%
\end{appendices}
\bibliography{bibliography}

@article{wang2013,
	author = {Wang, Z.J. and Fidkowski, Krzysztof and Abgrall, Rémi and Bassi, Francesco and Caraeni, Doru and Cary, Andrew and Deconinck, Herman and Hartmann, Ralf and Hillewaert, Koen and Huynh, H.T. and Kroll, Norbert and May, Georg and Persson, Per-Olof and van Leer, Bram and Visbal, Miguel},
	title = {High-order {CFD} methods: current status and perspective},
	journal = {International Journal for Numerical Methods in Fluids},
	volume = {72},
	number = {8},
	pages = {811-845},
	doi = {https://doi.org/10.1002/fld.3767},
	url = {https://onlinelibrary.wiley.com/doi/abs/10.1002/fld.3767},
	year = {2013},
}

@article{Tsoutsanis2025,
	title = {Arbitrary-order unstructured finite-volume methods for implicit large eddy simulation of turbulent flows with adaptive dissipation/dispersion adjustment (ADDA)},
	journal = {Journal of Computational Physics},
	volume = {523},
	pages = {113653},
	year = {2025},
	issn = {0021-9991},
	doi = {https://doi.org/10.1016/j.jcp.2024.113653},
	url = {https://www.sciencedirect.com/science/article/pii/S002199912400901X},
	author = {Panagiotis Tsoutsanis and Xesus Nogueira}
}

@article{Nishikawa2025,
	title = {An economical third-order nodal-gradient cell-centered finite-volume method for mixed-element grids},
	journal = {Journal of Computational Physics},
	volume = {540},
	pages = {114292},
	year = {2025},
	issn = {0021-9991},
	doi = {https://doi.org/10.1016/j.jcp.2025.114292},
	url = {https://www.sciencedirect.com/science/article/pii/S0021999125005753},
	author = {Hiroaki Nishikawa and Jeffery A. White}
}

@article{Costa2022,
	title = {Very high-order accurate finite volume scheme for the steady-state incompressible Navier–Stokes equations with polygonal meshes on arbitrary curved boundaries},
	journal = {Computer Methods in Applied Mechanics and Engineering},
	volume = {396},
	pages = {115064},
	year = {2022},
	issn = {0045-7825},
	doi = {https://doi.org/10.1016/j.cma.2022.115064},
	url = {https://www.sciencedirect.com/science/article/pii/S0045782522002729},
	author = {Ricardo Costa and Stéphane Clain and Gaspar J. Machado and João M. Nóbrega}
}

@article{Costa2023,
	title = {Imposing slip conditions on curved boundaries for {3D} incompressible flows with a very high-order accurate finite volume scheme on polygonal meshes},
	journal = {Computer Methods in Applied Mechanics and Engineering},
	volume = {415},
	pages = {116274},
	year = {2023},
	issn = {0045-7825},
	doi = {https://doi.org/10.1016/j.cma.2023.116274},
	url = {https://www.sciencedirect.com/science/article/pii/S0045782523003985},
	author = {Ricardo Costa and Stéphane Clain and Gaspar J. Machado and João M. Nóbrega and Hugo {Beirão da Veiga} and Francesca Crispo}
}

@article{Antoniadis2022,
	title = {{UCNS3D}: An open-source high-order finite-volume unstructured {CFD} solver},
	journal = {Computer Physics Communications},
	volume = {279},
	pages = {108453},
	year = {2022},
	issn = {0010-4655},
	doi = {https://doi.org/10.1016/j.cpc.2022.108453},
	url = {https://www.sciencedirect.com/science/article/pii/S0010465522001722},
	author = {Antonis F. Antoniadis and Dimitris Drikakis and Pericles S. Farmakis and Lin Fu and Ioannis Kokkinakis and Xesús Nogueira and Paulo A.S.F. Silva and Martin Skote and Vladimir Titarev and Panagiotis Tsoutsanis}
}

@article{Demirdzic2016,
	author =       {I. Demird{\v{z}}i\'{c}},
	title =        {A fourth-order finite volume method for structural
	analysis},
	journal =      {Applied Mathematical Modelling},
	year =         {2016},
	volume =       {40},
	number = {4},
	pages =        {3104-3114},
	doi =          {10.1016/j.apm.2015.09.098}
}

@article{Castrillo2022,
	title = {High-order finite volume method for linear elasticity on unstructured meshes},
	author = {Pablo Castrillo and Alfredo Canelas and Eugenio Schillaci and Joaquim Rigola and Asensio Oliva},
	journal = {Computers \& Structures},
	volume = {268},
	pages = {106829},
	year = {2022},
	issn = {0045-7949},
	doi = {https://doi.org/10.1016/j.compstruc.2022.106829},
	url = {https://www.sciencedirect.com/science/article/pii/S004579492200089X}
}

@article{Castrillo2024,
	title = {High-order cell-centered finite volume method for solid dynamics on unstructured meshes},
	journal = {Computers \& Structures},
	volume = {295},
	pages = {107288},
	year = {2024},
	issn = {0045-7949},
	doi = {https://doi.org/10.1016/j.compstruc.2024.107288},
	url = {https://www.sciencedirect.com/science/article/pii/S0045794924000178},
	author = {Pablo Castrillo and Eugenio Schillaci and Joaquim Rigola}
}

@article{Tzanio2021,
	author = {Tzanio Kolev and Paul Fischer and Misun Min and Jack Dongarra and Jed Brown and Veselin Dobrev and Tim Warburton and Stanimire Tomov and Mark S Shephard and Ahmad Abdelfattah and Valeria Barra and Natalie Beams and Jean-Sylvain Camier and Noel Chalmers and Yohann Dudouit and Ali Karakus and Ian Karlin and Stefan Kerkemeier and Yu-Hsiang Lan and David Medina and Elia Merzari and Aleksandr Obabko and Will Pazner and Thilina Rathnayake and Cameron W Smith and Lukas Spies and Kasia Swirydowicz and Jeremy Thompson and Ananias Tomboulides and Vladimir Tomov},
	title ={Efficient exascale discretizations: High-order finite element methods},
	journal = {The International Journal of High Performance Computing Applications},
	volume = {35},
	number = {6},
	pages = {527-552},
	year = {2021},
	doi = {10.1177/10943420211020803},
	URL = {https://doi.org/10.1177/10943420211020803},
}

@book{Bathe1996,
	author    = {K. J. Bathe},
	title     = {Finite Element Procedures},
	publisher = {Prentice Hall},
	address   = {New Jersey},
	year      = {1996},
	isbn      = {9780133014587},
	url       = {https://books.google.com/books/about/Finite_Element_Procedures.html?id=TKtVAAAAMAAJ}
}

@article{Demirdzic1988,
	author =       {I. Demird{\v{z}}i\'{c} and D. Martinovi\'{c} and
	A. Ivankovi\'{c}},
	title =        {Numerical simulation of thermal deformation in
	welded workpiece},
	journal =      {Zavarivanje},
	year =         {1988},
	volume =       {31},
	number =       {},
	pages =        {209-219},
	note =         {In Croatian. English translation available at \url{https://www.researchgate.net/publication/296148474_Numerical_simulation_of_thermal_deformation_in_welded_workpiece}},
}

@article{Cardiff2021,
	title={Thirty years of the finite volume method for solid mechanics},
	author={Cardiff, Philip and Demird{\v{z}}i{\'c}, Ismet},
	journal={Archives of Computational Methods in Engineering},
	volume={28},
	number={5},
	pages={3721--3780},
	year={2021},
	publisher={Springer},
	doi = {10.1007/s11831-020-09523-0}
}

@article{Felgueroso2007,
	title = {Finite volume solvers and Moving Least-Squares approximations for the compressible Navier–Stokes equations on unstructured grids},
	journal = {Computer Methods in Applied Mechanics and Engineering},
	volume = {196},
	number = {45},
	pages = {4712-4736},
	year = {2007},
	issn = {0045-7825},
	doi = {https://doi.org/10.1016/j.cma.2007.06.003},
	url = {https://www.sciencedirect.com/science/article/pii/S0045782507002484},
	author = {Luis Cueto-Felgueroso and Ignasi Colominas and Xesús Nogueira and Fermín Navarrina and Manuel Casteleiro}
}

@article{Ramirez2014,
	title = {A new higher-order finite volume method based on Moving Least Squares for the resolution of the incompressible Navier–Stokes equations on unstructured grids},
	journal = {Computer Methods in Applied Mechanics and Engineering},
	volume = {278},
	pages = {883-901},
	year = {2014},
	issn = {0045-7825},
	doi = {https://doi.org/10.1016/j.cma.2014.06.028},
	url = {https://www.sciencedirect.com/science/article/pii/S004578251400214X},
	author = {Luis Ramírez and Xesús Nogueira and Sofiane Khelladi and Jean-Camille Chassaing and Ignasi Colominas}
}

@phdthesis{Castrillo2023,
	title={High-order finite volume method for solid dynamics in fluid-structure interaction applications},
	author={Pablo Castrillo},
	year={2023},
	school={Universitat Polit{\`e}cnica de Catalunya},
	doi = {10.5821/dissertation-2117-396563}
}

@article{Gear1983,
	author  = {C. W. Gear and Y. Saad},
	title   = {Iterative solution of linear equations in {ODE} codes},
	journal = {SIAM Journal on Scientific and Statistical Computing},
	volume  = {4},
	number  = {4},
	pages   = {583--601},
	year    = {1983},
	doi     = {10.1137/0904040}
}

@article{Chan1984,
	author  = {T. F. Chan and K. R. Jackson},
	title   = {Nonlinearly preconditioned {Krylov} subspace methods for discrete {Newton} algorithms},
	journal = {SIAM Journal on Scientific and Statistical Computing},
	volume  = {5},
	number  = {3},
	pages   = {533--542},
	year    = {1984},
	doi     = {10.1137/0905039}
}

@article{Brown1986,
	author  = {P. N. Brown and A. C. Hindmarsh},
	title   = {Matrix-free methods for stiff systems of {ODE}'s},
	journal = {SIAM Journal on Numerical Analysis},
	volume  = {23},
	number  = {3},
	pages   = {610--638},
	year    = {1986},
	doi     = {10.1137/0723039}
}

@article{Brown1990,
	author  = {P. N. Brown and Y. Saad},
	title   = {Hybrid {Krylov} methods for nonlinear systems of equations},
	journal = {SIAM Journal on Scientific and Statistical Computing},
	volume  = {11},
	number  = {3},
	pages   = {450--481},
	year    = {1990},
	doi     = {10.1137/0911026}
}

@article{Qin2000,
	author    = {Ning Qin and David K. Ludlow and Scott T. Shaw},
	title     = {A matrix-free preconditioned {Newton}/{GMRES} method for unsteady {Navier}--{Stokes} solutions},
	journal   = {International Journal for Numerical Methods in Fluids},
	volume    = {33},
	number    = {2},
	pages     = {223--248},
	year      = {2000},
	month     = {May},
	doi       = {10.1002/(SICI)1097-0363(20000530)33:2<223::AID-FLD10>3.0.CO;2-V}
}

@article{Geuzaine2001,
	author = {Geuzaine, Philippe},
	title = {{Newton}-{Krylov} Strategy for Compressible Turbulent Flows on Unstructured Meshes},
	journal = {AIAA Journal},
	volume = {39},
	number = {3},
	pages = {528-531},
	year = {2001},
	doi = {10.2514/2.1339}
}

@article{Pernice2001,
	author  = {Michael Pernice and Michael D. Tocci},
	title   = {A Multigrid-Preconditioned {Newton}--{Krylov} Method for the Incompressible {Navier}--{Stokes} Equations},
	journal = {SIAM Journal on Scientific Computing},
	volume  = {23},
	number  = {2},
	pages   = {398--418},
	year    = {2001},
	doi     = {10.1137/S1064827500372250}
}

@article{Nejat2008,
	author  = {Amir Nejat and Carl Ollivier-Gooch},
	title   = {A high-order accurate unstructured finite volume {Newton}--{Krylov} algorithm for inviscid compressible flows},
	journal = {Journal of Computational Physics},
	volume  = {227},
	number  = {4},
	pages   = {2582--2609},
	year    = {2008},
	month   = {February},
	doi     = {10.1016/j.jcp.2007.11.011},
	publisher = {Elsevier}
}

@inproceedings{Nakashima2025,
	author = {Yoshitaka Nakashima and Yusuke Higo and Hiroaki Nishikawa},
	title = {Recent Algorithmic Advances in the scFLOW Unstructured-Polyhedral-Grid CFD Solver},
	booktitle = {AIAA SCITECH 2025 Forum},
	publisher = {American Institute of Aeronautics and Astronautics (AIAA)},
	doi = {10.2514/6.2025-0075},
	pages= {Paper No. 2025-0075},
	year ={2025}
}

@article{Mchugh1994,
	author    = {Paul R. McHugh and Dana A. Knoll},
	title     = {Comparison of standard and matrix-free implementations of several {Newton}--{Krylov} solvers},
	journal   = {AIAA Journal},
	volume    = {32},
	number    = {12},
	pages     = {2394--2400},
	year      = {1994},
	doi       = {10.2514/3.12305} 
}

@article{Knoll2004,
	author = {D.A. Knoll and D.E. Keyes},
	title = {Jacobian-free {Newton}–{Krylov} methods: a survey of approaches and applications},
	journal = {Journal of Computational Physics},
	volume = {193},
	number = {2},
	pages = {357-397},
	year = {2004},
	issn = {0021-9991},
	doi = {https://doi.org/10.1016/j.jcp.2003.08.010},
	url = {https://www.sciencedirect.com/science/article/pii/S0021999103004340}
}

@article{Vaassen2008,
	author    = {J.-M. Vaassen and D. Vigneron and J.-A. Essers},
	title     = {An implicit high order finite volume scheme for the solution of 3D {Navier}--{Stokes} equations with new discretization of diffusive terms},
	journal   = {Journal of Computational and Applied Mathematics},
	volume    = {215},
	number    = {2},
	pages     = {595--601},
	year      = {2008},
	month     = {June},
	publisher = {Elsevier},
	doi       = {10.1016/j.cam.2006.04.066}
}

@article{Lucas2010,
	author    = {Peter Lucas and Alexander H. van Zuijlen and Hester Bijl},
	title     = {Fast unsteady flow computations with a Jacobian-free {Newton}--{Krylov} algorithm},
	journal   = {Journal of Computational Physics},
	volume    = {229},
	number    = {24},
	pages     = {9201--9215},
	year      = {2010},
	month     = {December},
	publisher = {Elsevier},
	doi       = {10.1016/j.jcp.2010.08.033}
}

@article{Nishikawa2020,
	title = {A hyperbolic Poisson solver for tetrahedral grids},
	journal = {Journal of Computational Physics},
	volume = {409},
	pages = {109358},
	year = {2020},
	issn = {0021-9991},
	doi = {https://doi.org/10.1016/j.jcp.2020.109358},
	url = {https://www.sciencedirect.com/science/article/pii/S0021999120301327},
	author = {Hiroaki Nishikawa}
}

@misc{Cardiff2025,
	title={A Jacobian-free Newton-Krylov method for cell-centred finite volume solid mechanics}, 
	author={Philip Cardiff and Dylan Armfield and Željko Tuković and Ivan Batistić},
	year={2025},
	eprint={2502.17217},
	archivePrefix={arXiv},
	primaryClass={math.NA},
	note={Preprint available at: \url{https://arxiv.org/abs/2502.17217}},
}

@article{Weller1998,
	author =       {H. G. Weller and G. Tabor and H. Jasak and
	C. Fureby},
	title =        {A tensorial approach to computational continuum
	mechanics using object orientated techniques},
	journal =      {Computers in Physics},
	year =         {1998},
	volume =       {12},
	pages =        {620-631},
	doi = {10.1063/1.168744}
}

@article{Cardiff2025a,
	title={solids4foam: A toolbox for performing solid mechanics and fluid-solid interaction simulations in OpenFOAM},
	author={Philip Cardiff and Ivan Batisti{\'c} and Željko Tukovi{\'c}},
	journal={J. Open Source Softw.},
	year={2025},
	volume={10},
	pages={7407},
	doi = {10.21105/joss.07407}
}

@inproceedings{Nishikawa2010,
	author    = {Hiroaki Nishikawa},
	title     = {Beyond Interface Gradient: A General Principle for Constructing Diffusion Schemes},
	booktitle = {40th Fluid Dynamics Conference and Exhibit},
	year      = {2010},
	address   = {Chicago, Illinois, USA},
	month     = {June 28 -- July 1},
	doi       = {10.2514/6.2010-5093},
	pages = {Paper No. 2010-5093},
	publisher = {AIAA}
}

@article{Nishikawa2011,
	title = {Robust and accurate viscous discretization via upwind scheme – I: Basic principle},
	journal = {Computers \& Fluids},
	volume = {49},
	number = {1},
	pages = {62-86},
	year = {2011},
	issn = {0045-7930},
	doi = {https://doi.org/10.1016/j.compfluid.2011.04.014},
	url = {https://www.sciencedirect.com/science/article/pii/S0045793011001484},
	author = {Hiroaki Nishikawa}
}

@inproceedings{Nishikawa2025a,
	author = {Hiroaki Nishikawa},
	title = {On Second- and Higher-Order Finite-Volume Schemes for Viscous Terms on Unstructured Grids},
	booktitle = {AIAA AVIATION FORUM AND ASCEND 2025},
	pages = {Paper No. 2025-3674},
	publisher = {AIAA},
	doi = {10.2514/6.2025-3674},
	year = {2025}
}

@article{Nishikawa2021,
	author = {Nishikawa, Hiroaki},
	title = {Economically high-order unstructured-grid methods: Clarification and efficient FSR schemes},
	journal = {International Journal for Numerical Methods in Fluids},
	volume = {93},
	number = {11},
	pages = {3187-3214},
	doi = {https://doi.org/10.1002/fld.5028},
	year = {2021}
}

@article{Dunavant1985,
	author = {Dunavant, D. A.},
	title = {High degree efficient symmetrical Gaussian quadrature rules for the triangle},
	journal = {International Journal for Numerical Methods in Engineering},
	volume = {21},
	number = {6},
	pages = {1129–1148},
	year = {1985},
	doi = {10.1002/nme.1620210612}
}

@article{Shunn2012,
	author = {Shunn, L. and Ham, F. E.},
	title = {Symmetric quadrature rules for tetrahedra based on a cubic close-packed lattice arrangement},
	journal = {Journal of Computational and Applied Mathematics},
	volume = {236},
	number = {17},
	pages = {4348–4364},
	year = {2012},
	doi = {10.1016/j.cam.2012.03.032}
}

@article{Cardiff2017,
	author = {Cardiff, P. and {\v{Z}}. Tukovi{\'c} and Jaeger, P. De and Clancy, M. and Ivankovi{\'c}, A.},
	title = {A {L}agrangian cell-centred finite volume method for metal forming simulation},
	journal = {International {J}ournal for {N}umerical {M}ethods in {E}ngineering},
	volume = {109},
	number = {13},
	pages = {1777-1803},
	year = {2017},
	doi = {https://doi.org/10.1002/nme.5345},
}

@article{Cardiff2018,
	author =       {P. Cardiff and A. Kara\v{c} and P. De Jaeger and H. Jasak and J. Nagy and A. Ivankovi\'{c} and {\v{Z}}. Tukovi\'{c}},
	title =        {An open-source finite volume toolbox for solid mechanics and fluid-solid interaction simulations},
	journal = {arXiv preprint arXiv:1808.10736},
	year =         {2018},
	doi =          {10.48550/arXiv.1808.10736}
}

@article{Rhie1983,
	author = {Rhie, C. M. and Chow, W. L.},
	title = {Numerical study of the turbulent flow past an airfoil with trailing edge separation},
	journal = {AIAA Journal},
	volume = {21},
	number = {11},
	pages = {1525-1532},
	year = {1983},
	doi = {10.2514/3.8284}
}

@article{Tukovic2013,
	author = {Tukovi\'{c}, {\v{Z}} and Ivankovi\'{c}, A. and Kara\v{c}, A.},
	title = {Finite{-}volume stress analysis in multi{-}material linear elastic body},
	journal = {International Journal for Numerical Methods in Engineering},
	volume = {93},
	number = {4},
	pages = {400-419},
	year = {2013},
	doi = {https://doi.org/10.1002/nme.4390}
}

@article{Khelladi2011,
	title = {Toward a higher order unsteady finite volume solver based on reproducing kernel methods},
	journal = {Computer Methods in Applied Mechanics and Engineering},
	volume = {200},
	number = {29},
	pages = {2348-2362},
	year = {2011},
	issn = {0045-7825},
	doi = {https://doi.org/10.1016/j.cma.2011.04.001},
	url = {https://www.sciencedirect.com/science/article/pii/S0045782511001344},
	author = {Sofiane Khelladi and Xesús Nogueira and Farid Bakir and Ignasi Colominas}
}

@book{strang2012,
	title={Linear algebra and its applications},
	author={Strang, Gilbert},
	year={2012},
	publisher={Academic Press},
	address={Belmont},
	edition = {2},
	url = {https://openlibrary.org/books/OL4428412M/Linear_algebra_and_its_applications}
}

@article{Tsoutsanis2023,
	title = {Stencil selection algorithms for WENO schemes on unstructured meshes},
	journal = {Journal of Computational Physics},
	volume = {475},
	pages = {108840},
	year = {2023},
	issn = {0021-9991},
	doi = {https://doi.org/10.1016/j.jcp.2019.07.039},
	url = {https://www.sciencedirect.com/science/article/pii/S0021999119305248},
	author = {Panagiotis Tsoutsanis}
}

@article{Coppeans2026,
	title = {Anisotropic mesh adaptation for high-Order meshes in two dimensions},
	journal = {Journal of Computational Physics},
	volume = {545},
	pages = {114506},
	year = {2026},
	issn = {0021-9991},
	doi = {https://doi.org/10.1016/j.jcp.2025.114506},
	url = {https://www.sciencedirect.com/science/article/pii/S0021999125007880},
	author = {Alexander W.C. Coppeans and Krzysztof J. Fidkowski and Joaquim R.R.A. Martins}
}

@article{Costa2021,
	title = {Efficient very high-order accurate polyhedral mesh finite volume scheme for 3{D} conjugate heat transfer problems in curved domains},
	journal = {Journal of Computational Physics},
	volume = {445},
	pages = {110604},
	year = {2021},
	issn = {0021-9991},
	doi = {https://doi.org/10.1016/j.jcp.2021.110604},
	url = {https://www.sciencedirect.com/science/article/pii/S002199912100499X},
	author = {Ricardo Costa and João M. Nóbrega and Stéphane Clain and Gaspar J. Machado}
}

@article{Gooch2002,
	title = {A High-Order-Accurate Unstructured Mesh Finite-Volume Scheme for the Advection–Diffusion Equation},
	journal = {Journal of Computational Physics},
	volume = {181},
	number = {2},
	pages = {729-752},
	year = {2002},
	issn = {0021-9991},
	doi = {https://doi.org/10.1006/jcph.2002.7159},
	url = {https://www.sciencedirect.com/science/article/pii/S0021999102971597},
	author = {Carl Ollivier-Gooch and Michael {Van Altena}}
}

@article{Demirdzic2022,
	author = {I. Demirdžić and P. Cardiff},
	title = {Symmetry plane boundary conditions for cell-centered finite-volume continuum mechanics},
	journal = {Numerical Heat Transfer, Part B: Fundamentals},
	year = {2022},
	doi = {10.1080/10407790.2022.2105073}
}

@article{Tukovic2018,
	author =       {{\v{Z}}. Tukovi\'{c} and A. Kara\v{c} and P. Cardiff
	and H. Jasak and A. Ivankovi\'{c}},
	title =        {{OpenFOAM} finite volume solver for fluid-solid interaction},
	year =         {2018},
	volume =       {42},
	number =        {3},
	pages =        {1-31},
	journal =      {Transactions of {FAMENA}},
	doi =         {10.21278/TOF.42301},
}

@article{Demirdzic2015,
	author =       {I. Demird{\v{z}}i\'{c}},
	year =         {2015},
	title =        {On the Discretization of the Diffusion Term in
	Finite-Volume Continuum Mechanics},
	journal =      {Numerical Heat Transfer, Part B: Fundamentals: An
	International Journal of Computation and
	Methodology},
	volume =       {68:1},
	pages =        {1-10},
	doi =          {10.1080/10407790.2014.985992}
}

@Misc{PETSc,
	author = {Satish Balay and Shrirang Abhyankar and Mark~F. Adams and Jed Brown and Peter Brune
	and Kris Buschelman and Lisandro Dalcin and Victor Eijkhout and William~D. Gropp
	and Dinesh Kaushik and Matthew~G. Knepley
	and Lois Curfman McInnes and Karl Rupp and Barry~F. Smith
	and Stefano Zampini and Hong Zhang},
	title =  {{PETS}c. {W}eb page},
	url =    {http://www.mcs.anl.gov/petsc},
	year = {2015}
}

@phdthesis{Mazzanti2024,
	title={Coupled Vertex-Centred Finite Volume Methods for Large-Strain Elastoplasticity},
	author={Mazzanti, Federico},
	year={2024},
	school={University College Dublin},
	url = {http://hdl.handle.net/10197/29209}
}

@article{augarde2008,
	title = {The use of Timoshenko's exact solution for a cantilever beam in adaptive analysis},
	journal = {Finite Elements in Analysis and Design},
	volume = {44},
	number = {9},
	pages = {595-601},
	year = {2008},
	issn = {0168-874X},
	doi = {https://doi.org/10.1016/j.finel.2008.01.010},
	url = {https://www.sciencedirect.com/science/article/pii/S0168874X08000140},
	author = {Charles E. Augarde and Andrew J. Deeks}
}

@article{Demirdzic1997,
	title={Benchmark solutions of some structural analysis problems using finite-volume method and multigrid acceleration},
	author={Demird{\v{z}}i{\'c}, I and Muzaferija, S and Peri{\'c}, M},
	journal={International journal for numerical methods in engineering},
	volume={40},
	number={10},
	pages={1893--1908},
	year={1997},
	publisher={Wiley Online Library},
	doi = {10.1002/(SICI)1097-0207(19970530)40:10<1893::AID-NME146>3.0.CO;2-L}
}

@article{Bijelonja2006,
	author =       {I. Bijelonja and I. Demird{\v{z}}i\'{c} and
	S. Muzaferija},
	title =        {A finite volume method for incompressible linear
	elasticity},
	journal =      {Computer Methods in Applied Mechanics and Engineering},
	year =         {2006},
	volume =       {195},
	pages =        {6378-6390},
	doi = {10.1016/j.cma.2006.01.005}
}

@article{Bijelonja2005a,
	author =       {I. Bijelonja and I. Demird{\v{z}}i\'{c} and
	S. Muzaferija},
	title =        {A finite volume method for large strain analysis of
	incompressible hyperelastic materials},
	journal =      {International Journal for Numerical Methods in Engineering},
	year =         {2005},
	volume =       {64},
	pages =        {1594-1609},
	doi = {10.1002/nme.1413}
}

@book{Timoshenko1970,
	title={Theory of Elasticity},
	author={Timoshenko, S and Goodier, JN},
	edition   = {3},
	publisher={McGraw--Hill Book Company},
	address={India},
	year={1970},
	isbn      = {9780070858053},
	url = {https://books.google.com.bo/books/about/Theory_of_elasticity.html?id=HSfto8If1moC}
}

@book{Green1992,
	title={Theoretical elasticity},
	author={Green, Albert Edward and Zerna, Wolfgang},
	year={1992},
	publisher={Dover Publications},
	edition =      {2},
	address =      {New York},
	isbn      = {9780486670768},
	url = {https://books.google.com.uy/books/about/Theoretical_Elasticity.html?id=Jg6oAAAAIAAJ&redir_esc=y}
}

@Article{Horvat2025,
	AUTHOR = {Horvat, Anja and Milovi\'{c}, Philipp and Kar\v{s}aj, Igor and Tukovi\'{c}, {\v{Z}}eljko},
	TITLE = {A Block-Coupled Finite Volume Method for Incompressible Hyperelastic Solids},
	JOURNAL = {Applied Sciences},
	VOLUME = {15},
	YEAR = {2025},
	NUMBER = {23},
	ARTICLE-NUMBER = {12660},
	URL = {https://www.mdpi.com/2076-3417/15/23/12660},
	ISSN = {2076-3417},
	doi = {10.3390/app152312660}
}

@inproceedings{oliveira2020,
	title={Implementation and numerical verification of an incompressible three-parameter Mooney-Rivlin model for large deformation of soft tissues},
	author={Oliveira, I and Gasche, JL and Cardiff, P},
	booktitle={The 15th OpenFOAM Workshop},
	year={2020},
	url ={https://openfoam-extend.sourceforge.net/OpenFOAM_Workshops/OFW15_2020_Arlington/program/s106.html}
}

@book{Zienkiewicz2000,
	title={The finite element method},
	author={Zienkiewicz, Olgierd Cecil and Taylor, Robert Leroy},
	volume={2},
	year={2000},
	publisher={Butterworth-heinemann},
	address={Berlin, Germany},
	url={https://books.google.com.uy/books/about/The_Finite_Element_Method_Solid_mechanic.html?id=MhgBfMWFVHUC&redir_esc=y}
}

@misc{Simplas,
	author       = {P. Areias},
	title        = {{Simplas}},
	publisher    = {Portuguese Software Association (ASSOFT)},
	url          = {http://www.simplassoftware.com.},
	year = {1996}
}

@misc{Pelteret2018,
	author       = {Pelteret, Jean-Paul and McBride, Andrew},
	title        = {{The deal.II code gallery: Quasi-Static Finite-Strain Compressible Elasticity}},
	month        = apr,
	year         = 2018,
	publisher    = {Zenodo},
	doi          = {10.5281/zenodo.1228964},
}

@article{Southwell1926,
	author = {R.V. Southwell and H.J. Gough},
	title = {VI. On the concentration of stress in the neighbourhood of a small spherical flaw; and on the propagation of fatigue fractures in “Statistically Isotropic” materials },
	journal = {The London, Edinburgh, and Dublin Philosophical Magazine and Journal of Science},
	volume = {1},
	number = {1},
	pages = {71--97},
	year = {1926},
	publisher = {Taylor \& Francis},
	doi = {10.1080/14786442608633614}
}

@article{Goodier1933,
	title={Concentration of stress around spherical and cylindrical inclusions and flaws},
	author={Goodier, J. N.},
	journal={Journal of Applied Mechanics},
	volume={1},
	number={2},
	pages={39--44},
	year={1933},
	doi = {10.1115/1.4012173}
}

@article{Tukovic2007,
	author  = {Tukovi{\'c}, {\v{Z}}eljko and Jasak, Hrvoje},
	title   = {Updated Lagrangian finite volume solver for large deformation dynamic response of elastic body},
	journal = {Transactions of FAMENA},
	volume  = {31},
	number  = {1},
	pages   = {55--70},
	year    = {2007},
	note = {\url{https://www.croris.hr/crosbi/publikacija/prilog-casopis/141726}}
}

@article{Syrakos2023,
	title={A unification of least-squares and Green--Gauss gradients under a common projection-based gradient reconstruction framework},
	author={Syrakos, Alexandros and Oxtoby, Oliver and de Villiers, Eugene and Varchanis, Stylianos and Dimakopoulos, Yannis and Tsamopoulos, John},
	journal={Mathematics and Computers in Simulation},
	volume={205},
	pages={108--141},
	year={2023},
	month={March},
	doi={10.1016/j.matcom.2022.09.008},
	publisher={Elsevier}
}

@inproceedings{Jalali2013,
	title={Higher-order finite volume solution reconstruction on highly anisotropic meshes},
	author={Jalali, Alireza and Ollivier Gooch, Carl F},
	booktitle={21st AIAA Computational Fluid Dynamics Conference},
	pages={2565},
	year={2013},
	doi= {https://doi.org/10.2514/6.2013-2565}
}

@article{Rivlin1951,
	author = {Rivlin, R. S.  and Saunders, D. W. },
	title = {Large elastic deformations of isotropic materials VII. Experiments on the deformation of rubber},
	journal = {Philosophical Transactions of the Royal Society of London. Series A, Mathematical and Physical Sciences},
	volume = {243},
	number = {865},
	pages = {251-288},
	year = {1951},
	doi = {10.1098/rsta.1951.0004}
}

@article{Costa2022b,
	title = {A novel approach for temporal simulations with very high-order finite volume schemes on polyhedral unstructured grids},
	journal = {Journal of Computational Physics},
	volume = {453},
	pages = {110960},
	year = {2022},
	issn = {0021-9991},
	doi = {https://doi.org/10.1016/j.jcp.2022.110960},
	url = {https://www.sciencedirect.com/science/article/pii/S0021999122000225},
	author = {Pedro M.P. Costa and Duarte M.S. Albuquerque}
}

@article{geuzaine2009gmsh,
	title={Gmsh: A 3-D finite element mesh generator with built-in pre- and post-processing facilities},
	author={Geuzaine, Christophe and Remacle, Jean-Fran{\c{c}}ois},
	journal={International journal for numerical methods in engineering},
	volume={79},
	number={11},
	pages={1309--1331},
	year={2009},
	publisher={Wiley Online Library},
	doi={10.1002/nme.2579}
}

@book{Gurtin2010,
	author = {Gurtin, Morton E. and Fried, Eliot. and Anand, Lallit.},
	isbn = {9780521405980},
	pages = {694},
	publisher = {Cambridge University Press},
	title = {The Mechanics and Thermodynamics of Continua},
	year = {2010},
	edition   = {1},
	address   = {Cambridge}, 
	url = {https://www.cambridge.org/highereducation/books/the-mechanics-and-thermodynamics-of-continua/D098B8CCA2EC48608BC2444FBD03D43F},
}

@book{Gurtin1981,
	author = {Gurtin, Morton E.},
	title = {An Introduction to Continuum Mechanics},
	series = {Mathematics in Science and Engineering}, 
	volume = {158}, 
	publisher = {Academic Press, Inc.},
	year = {1981},
	address = {New York},
	url = {https://www.sciencedirect.com/bookseries/mathematics-in-science-and-engineering/vol/158/}
}
\end{document}